%% file: master.tex
\documentclass[leqno]{amsart}
\usepackage[utf8]{inputenc}
\usepackage[margin=3cm]{geometry}
\usepackage{extdash}
\usepackage{amsthm,amsfonts,amssymb,amsmath, microtype}
\usepackage[all]{xy}
\usepackage{xr-hyper}
\usepackage[colorlinks=true]{hyperref}
\usepackage{verbatim}
\usepackage{pdfsync}
\usepackage{stmaryrd}
\usepackage{marvosym}

\usepackage{mathtools}
\usepackage{enumitem}

\usepackage{tikz}
\usetikzlibrary{arrows,calc,positioning}
\usetikzlibrary{circuits.logic.IEC}
\usetikzlibrary{cd}
\usepackage{xcolor}
\usetikzlibrary{shapes}
\usetikzlibrary{decorations.markings}
\input{tikzstyledefs.tex}

\newcommand{\I}{\mathrm{I}}
\newcommand{\Z}{\mathbb{Z}}
\newcommand{\Q}{\mathbb{Q}}
\newcommand{\N}{\mathbb{N}}

\newcommand{\tr}{\mathrm{tr}}
\newcommand{\id}{\mathrm{Id}}
\newcommand{\im}{\mathrm{im}}

\newcommand{\Mod}{\mathrm{Mod}}

\newcommand{\coker}{\mathrm{coker}}

\newcommand{\A}{\mathcal{A}}
\newcommand{\C}{\mathcal{C}}
\newcommand{\T}{\mathcal{T}}
\newcommand{\F}{\mathcal{F}}
\newcommand{\D}{\mathcal{D}}
\newcommand{\Se}{\mathcal{S}}
\newcommand{\unit}{\mathbb{I}}
\newcommand{\II}{\mathrm{I}}
\newcommand{\PP}{\mathrm{P}}

\newcommand{\Fun}{\mathrm{Fun}}
\newcommand{\hproj}{\mathrm{\text{h-}proj}}

\newcommand{\Mor}{\mathrm{Mor}}
\newcommand{\CC}{\mathrm{C}}
\newcommand{\CCC}{\mathbb{C}}
\newcommand{\kk}{\mathsf{k}}
\newcommand{\KK}{\mathrm{K}}
\newcommand{\DD}{\mathrm{D}}
\newcommand{\Aut}{\mathrm{Aut}}
\newcommand{\Gr}{\mathrm{Gr}}

\newcommand{\Cone}{\mathrm{Cone}}

\newcommand{\colim}{\mathrm{colim}}
\newcommand{\Quot}{\mathrm{Quot}}

\newcommand{\nocontentsline}[3]{}
\newcommand{\tocless}[2]{\bgroup\let\addcontentsline=\nocontentsline#1{#2}\egroup}

\newtheorem{theorem}{Theorem}[section]
\newtheorem{proposition}[theorem]{Proposition}
\newtheorem{lemma}[theorem]{Lemma}

\newtheorem{corollary}[theorem]{Corollary}

\theoremstyle{definition}
\newtheorem{definition}[theorem]{Definition}
\newtheorem{example}[theorem]{Example}
\newtheorem{observation}[theorem]{Observation}

\newtheorem{remark}[theorem]{Remark}

\begin{document}

\author{Felix K\"ung}
\date{\today}

\title{Structures on the category of $N$-complexes}

\maketitle

\tocless{\section*{Abstract}}
The theory of $N$-complexes is a generalization of both ordinary chain complexes and graded objects. Hence it yields deeper insight in the structure of these and offers a broader range of applications.
This work generalizes the tensor product of chain complexes and graded objects to the case of $N$-complexes using the structures of $q$-binomial coefficients. We then study different approaches to realize the derived category of $N$-complexes. In particular we realize it as the Verdier quotient of the homotopy category of $N$-complexes, as the $\mathrm{h}$-projective objects and as the homotopy category of a category admitting a Quillen model structure.

\tableofcontents

\newpage

\section*{Introduction}

The theory of chain complexes is a very thoroughly studied field in mathematics. Since a defining property of these is $d^2=0$, an obvious generalization is the notion of an $N$-complex, which has the same structure as a chain complex but instead of $d^2=0$ we have $d^N=0$.

Objects of this form were first studied by W. Mayer \cite{mayer} in 1942 as an alternative homology theory on topological spaces. This approach did not get significant attention since this $N$-homology theory fulfills the Eilenberg-Steenrod axioms and therefore has to be equivalent to ordinary homology. Later in 1991 Kapranov studied them from a purely algebraic point of view \cite{kapranov} which spiked great interest in that topic. In particular a considerable amount of homological framework for these was developed by M. Dubois-Violette \cite{DuboisViolette1998,DuboisViolette1996}. One of the main motivations of these papers was that $N$-complexes can be used to study representations of quantum groups similarly to the way one uses chain complexes to study representations of groups. Furthermore the concept of cohomology of $N$-complexes also got attention in physics where it is referred to as amplitude cohomology, or homology depending on the convention \cite{CSW}. After this initial work there have been papers about the relation of $N$-complexes to other topics such as Koszulity \cite{Berger}, simplicial structures \cite{Djalal1,Djalal2} or dualizing categories \cite{YHNZ}. An especially interesting example of the work done in this field is the paper ``Algebraic Model for $N$-connection'' \cite{VAOL} since it applies the otherwise purely algebraic notions of $N$-complexes to geometry and suggests an application of $N$-complexes to non-flat connections. During the last few years the concept of $N$-complexes has also arrived in the field of  categorification: in particular categorification at roots of unity relies on the theory of $N$-complexes. For example M. Khovanov used hopfological algebra to model $N$-complexes in his paper ``Hopfological algebra and categorification at a root of unity: the first steps'' \cite{khovanov}.

Although there was already quite some theory on the applications and properties of $N$-complexes, a rigid foundation of the homological properties and constructions used in these papers was missing. This has been remedied by Iyama, Kato and Miyachi in \cite{N-complex}. However, some of the notions defined in their paper deserve a more thorough study, and although one has now this solid foundation of the theory involved, there are still some open questions on the structure of $N$-complexes. \\
For example:
\begin{itemize}
\item Can we construct a canonical monoidal structure on the category of $N$-complexes if the underlying category is monoidal?
\item Can one characterize the derived category of $N$-complexes as well over h-projectives?
\item Is the derived category really the Gabriel-Zismann localisation of the category of $N$-complexes along $N$-quasi-isomorphisms?
\end{itemize}
Throughout this work we will try to answer these questions and give very thorough study of the properties of the involved structures. 

In the first three sections we will collect some general properties of $N$-complexes, study $q$-binomial coefficients and apply these to lift bifunctors using a twisted Leibniz rule. This involves the construction of certain standard $N$-complexes, proving that the category of $N$-complexes is complete, cocomplete, abelian or additive if the underlying category is, a weaker but more intuitive form of cyclic sieving and how we can lift which functors and bifunctors.

In the remaining section~\ref{section 4} we then apply these structures of $q$-binomial coefficients and liftable functors to lift monoidal structures in most cases.

We continue in section~\ref{section 5} and \ref{section 6} by studying the homotopy structure of the category of $N$-complexes, we start by constructing its homotopy category as stable category of a Frobenius exact sturcture and prove that it is triangulated.

Using this structure we imitate in section~\ref{section 7} the three classical constructions for realizing the derived category as Verdier-quotient, $\mathrm{h}$-projective objects and prove that these can be identified by considering semifree objects as $\mathrm{h}$-projective resolutions.

Finally in section~\ref{section 8} we finish by defining a Quillen model structure on the category of $N$-complexes realizing the derived category as homotopy category. This shows in particular that the derived category of $N$-complexes indeed is the Gabriel-Zismann localisation along $N$-quasi-isomorphisms.

Although the model structure constructed in section~\ref{section 8} coincides with the model structure constructed by  J. Gillespie and M. Hovey in their paper \cite{Gillespie2}, our work is completely independent and also uses different techniques that are more hands on. Furthermore it is quite natural that both approaches yield the same model structure since it is the obvious generalization of the projective model structure on the category of chain complexes.

\vspace*{\fill} 
\tocless{\section*{Acknowledgement}}

I would like to express my gratitude to my advisor O. Schnürer for pointing me towards this topic and his helpful comments and engagement.
 Furthermore I would like to thank G. Jasso for his lecture on dg-structures and insight on classical constructions of chain complexes.
Finally I want to thank G. Janssens for suggesting to finally publish this work and believing in it.

\input{Category_CN}
\input{interlude_q-binomial}

\input{monoidal_structure}

\input{Frobenius_exact_structure}

\input{derived_category}
\input{model_structure}

\newpage

\section*{Conclusion and Outlook}

We construct, using $q$-binomial coefficients, a theory of lifted bifunctors to the category of $N$-complexes and apply these to construct an induced monoidal structure on the category of $N$-complexes for countable distributive monoidal categories and prove that this monoidal structure is closed if the underlying category is the category of modules over a ring. After that we study homotopy and cohomology of $N$-complexes to get a characterization of the derived category of $N$-complexes over a ring as a Verdier quotient or h-projective complexes. Finally we establish a cofibrantly generated model structure on the category of $N$-complexes of modules over a ring which yields as homotopy category again the derived category of $N$-complexes and in particular shows that the derived category is indeed equivalent to the Gabriel-Zismann localisation along the $N$-quasi-isomorphisms.

The approaches presented here may be generalized in several ways. For instance one could consider the category of modules over a small category in order to show that this category still permits the arguments used for constructing the derived category and the model structure. Another approach would be to consider monoid objects in our monoidal structure, which would lead to the theory of $N$-dg algebras.

Especially, since both of these possible generalizations point in the direction of a possible generalization to a theory of dg-$N$-categories, we believe that a detailed study of this topic merits further attention and we intend to address it in future work.

\bibliographystyle{abbrv}
\bibliography{master.bib}

\end{document}

%% file: tikzstyledefs.tex
\tikzstyle heightone=[scale=.7,shift={(0,-.3)}]
\tikzstyle heightones=[scale=.8,xscale=.35,shift={(0,.1)}]
\tikzstyle heightoneonehalf=[scale=.9,shift={(0,-.2)}]
\tikzstyle heighttwo=[scale=.9,shift={(0,-.4)}]
\tikzstyle heighttwos=[scale=.5,xscale=.6,shift={(0,-.1)}]
\tikzstyle heightthree=[scale=.6,shift={(0,-.9)}]
\tikzstyle heightthrees=[scale=.4,xscale=.7,shift={(0,-.2)}]

\tikzstyle arrowstyle=[blue,semitransparent,scale=2]

\tikzstyle basiclabel=[draw=none,fill=none,shape=rectangle,inner sep=2pt,scale=.8]
\tikzstyle leftlabel=[basiclabel,anchor=east]
\tikzstyle rightlabel=[basiclabel,anchor=west]
\tikzstyle bottomlabel=[basiclabel,anchor=north]
\tikzstyle toplabel=[basiclabel,anchor=south]

\tikzstyle vertex=[circle,draw,fill=black,inner sep=1pt]
\tikzstyle ciliation=[circle,draw=none,fill=red,inner sep=1pt,semitransparent]
\tikzstyle ciliatednode=[vertex,pin={[pin distance=1mm,pin edge={semitransparent,red},ciliation]#1:{}}]

\tikzstyle vector=[black,thick,rectangle,draw=gray!50!yellow,top color=yellow!30,bottom color=black!10,scale=.8,inner sep=2pt]
\tikzstyle small vector=[vector,scale=.8]
\tikzstyle plain vector=[rectangle,draw=none,fill=white,scale=.7]

\tikzstyle my signal=[black,thick,signal,signal pointer angle=120,draw=blue!50,top color=blue!20,bottom color=black!10,scale=.8,inner sep=2pt]
\tikzstyle matrix=[my signal,signal from=south,signal to=north]
\tikzstyle reverse matrix=[my signal,signal from=north,signal to=south]

\tikzstyle small matrix=[matrix,scale=.7]
\tikzstyle reverse small matrix=[reverse matrix,scale=.7]
\tikzstyle matrix on edge=[small matrix,sloped,rotate=-90]
\tikzstyle reverse matrix on edge=[small matrix,sloped,rotate=90]

\tikzstyle trivalent=[very thick]
\tikzstyle dotdotdot=[decorate,decoration={markings,
    mark=at position .3 with{\node{.};},
    mark=at position .5 with {\node{.};},
    mark=at position .7 with {\node{.};}}]

\tikzstyle wavyup=[out=90,in=-90]
\tikzstyle wavydown=[out=-90,in=90]

\tikzstyle symmetrizer=[rectangle,fill=gray!10,draw=black]
\tikzstyle permutation=[symmetrizer]
\tikzstyle antisymmetrizer=[rectangle,fill=black,draw=black]
\tikzstyle symlabel=[draw=none,fill=none,black,scale=.8]
\tikzstyle asymlabel=[draw=none,fill=none,white,scale=.8]

%% file: Category_CN.tex
\section{The category $\CC_N(\A)$}

The aim of this section is to introduce the concept of $N$-complexes. In particular we will give the definition and first properties of the category of $N$-complexes in a pointed category. We will use a general definition of diagrams that look like $N$-complexes, without restriction on the differential, which we will refer to as $\infty$-complexes. This allows the use of some abstract arguments for the structure and existence of this category and later also simplifies certain proofs about properties of the category of $N$-complexes. After that we will define the category of $N$-complexes as a full subcategory and discuss certain ``standard'' $N$-complexes which resemble the standard projective chain complexes and admit similar adjunctions as in the chain-complex case. Later we will prove that the category of $N$-complexes and $\infty$-complexes inherits properties from their base category. In particular we will show that the category of $N$-complexes is abelian if the ground category is, and that it admits colimits and limits if the ground category does.

Let us begin by defining the category of $\infty$-complexes as the category of diagrams of a certain shape in an arbitrary category.
\begin{definition}
Let $\C$ be a category.
\begin{enumerate}
\item An \textbf{$\infty$-complex in $\C$} is a diagram
\begin{equation*}
\tikz[heighttwo,xscale=2,yscale=2,baseline]{

\node (P0) at (0,0){$\cdot \cdot \cdot$};
\node (X1) at (1,0){$X^{n-1}$};
\node (X2) at (2,0){$X^n$};
\node (X3) at (3,0){$X^{n+1}$};

\node (P2) at (4,0){$\cdot \cdot \cdot$};
    
\draw[->]
(P0) edge node[above] {$d_X$} (X1)
(X1) edge node[above] {$d_X$} (X2)
(X2) edge node[above] {$d_X$} (X3)
(X3) edge node[above] {$d_X$} (P2);
}
\end{equation*}
in $\C$.

We will sometimes drop the $\C$ if it is clear from the context.
\item A \textbf{morphism $g:X \to Y$ of $\infty$-complexes in $\C$} is a family of morphisms $\left(g^i:X^i \to Y^i\right)_{i \in \Z}$, such that the following diagram in $\C$ commutes:
\begin{equation*}
\tikz[heighttwo,xscale=2,yscale=2,baseline]{

\node (P0) at (0,1){$\cdot \cdot \cdot$};
\node (X1) at (1,1){$X^{i-1}$};
\node (X2) at (2,1){$X^i$};
\node (X3) at (3,1){$X^{i+1}$};

\node (P2) at (4,1){$\cdot \cdot \cdot$};

\node (P1) at (0,0){$\cdot \cdot \cdot$};
\node (Y1) at (1,0){$Y^{i-1}$};
\node (Y2) at (2,0){$Y^i$};
\node (Y3) at (3,0){$Y^{i+1}$};
\node (P3) at (4,0){$\cdot \cdot \cdot$};
\node[overlay] (Point) at (4.5,0){.};

\draw[->]
(P0) edge node[above] {$d_X$} (X1)
(X1) edge node[above] {$d_X$} (X2)
(X2) edge node[above] {$d_X$} (X3)
(X3) edge node[above] {$d_X$} (P2)

(P1) edge node[above] {$d_Y$} (Y1)
(Y1) edge node[above] {$d_Y$} (Y2)
(Y2) edge node[above] {$d_Y$} (Y3)
(Y3) edge node[above] {$d_Y$} (P3)

(X1) edge node[right] {$g^{i-1}$} (Y1)
(X2) edge node[right] {$g^i$} (Y2)
(X3) edge node[right] {$g^{i+1}$} (Y3);
}
\end{equation*}

\end{enumerate}
\end{definition}

\begin{lemma}\label{1.2}

The collection of $\infty$-complexes in $\C$ together with morphisms of $\infty$-complexes in $\C$ defines a category, which we will refer to as $\CC_\infty (\C)$.

\begin{proof}
Let $\I$ be the category with objects $x \in \Z$ and $\I(x,y):=\{*\}$ if $x\le y$ and $\I(x,y):=\emptyset$ otherwise. Then $\CC_\infty (\C)$ is the category of functors from $\I$ to $\C$.
\end{proof}
\end{lemma}

\begin{corollary}\label{CN additive}
If $\C$ is  pre-additive, additive or abelian, so is $\CC_\infty(\C)$.
\begin{proof}
This is clear by the identification of $\CC_\infty(\C)$ as functors from  $\I$ to $\C$.
\end{proof}
\end{corollary}

Now we can define the category of $N$-complexes as a full subcategory of the category of $\infty$-complexes. 
\begin{definition}
Let $\C$ be a pointed category, i.\,e.\ a category with a zero object, and $N\in \N$.
\begin{enumerate}
\item An $\infty$-complex $X$ in $\C$ is an \textbf{$N$-complex in $\C$} or just an \textbf{$N$-complex}, if $d_X^N=0$.
\item We denote the full subcategory of $\CC_\infty(\C)$ consisting of the $N$-complexes by \textbf{$\CC_N(\C)$}.
\end{enumerate}
\end{definition}

It is particularly interesting to observe that the notion of $N$-complexes generalizes the following cases.

\begin{observation}
There are three well known cases for the category of $N$-complexes in a pointed category $\C$:
\begin{itemize}
\item[$N=0$:] In this case we have for every $n \in \Z$ that $d_X^0=\id_{X^n}=0$, hence $X^n\cong 0$ and so $\CC_0(\C)$ is the subcategory consisting of zero objects in $\CC_\infty(\C)$.  
\item[$N=1$:] Here the above equation yields $d_X=0$. Since morphisms in a pointed category always commute with the zero map we get that the category $\CC_1(\C)$ is isomorphic to the category of graded objects in $\C$.
\item[$N=2$:] This case yields by definition the case of chain complexes in $\C$.
\end{itemize}
\end{observation}

Since we will later imitate constructions from the theory of chain complexes, we will need to define the following objects introduced in \cite{N-complex} and prove that they give rise to certain adjunctions. This will be one of the main ingredients for our constructions regarding the homotopy category of $N$-complexes and the projective model structure.

\begin{definition}
 Let $\C$ be a pointed category.
\begin{enumerate}

\item Let $X \in \C$ and $i \in \Z$. The diagram 
\begin{equation*}
\tikz[heighttwo,xscale=1.7,yscale=2,baseline]{

\node (P0) at (0,0){$\cdot \cdot \cdot$};
\node (X1) at (1,0){0};
\node (X2) at (2,0){$X$};
\node (X3) at (3,0){$X$};
\node (P1) at (4,0){$\cdot \cdot \cdot$}; 
\node (XN-2) at (5,0){$X$};
\node (XN-1) at (6,0){$X$};
\node (XN) at (7,0){0};
\node (P2) at (8,0){$\cdot \cdot \cdot$};

\node (P3) at (0,0.5){\tiny{$\cdot \cdot \cdot$}};
\node (i-n-1) at (1,0.5){\tiny{ $i-N$}};
\node (i-n) at (2,0.5){\tiny{$i-N+1$}};
\node (i-n+1) at (3,0.5){\tiny{$i-N+2$}};
\node (P4) at (4,0.5){\tiny{$\cdot \cdot \cdot$}};
\node (i-1) at (5,0.5){\tiny{$i-1$}};
\node (i) at (6,0.5){\tiny{$i$}};
\node (i+1) at (7,0.5){\tiny{$i+1$}};
\node (P5) at (8,0.5){\tiny{$\cdot \cdot \cdot$}};

\node[overlay] (Point) at (8.5,0){.};
    
\draw[->]
(P0) edge node[above] {0} (X1)
(X1) edge node[above] {0} (X2)
(X2) edge node[above] {$\id_{X}$} (X3)
(X3) edge node[above] {$\id_{X}$} (P1)
(P1) edge node[above] {$\id_{X}$} (XN-2)
(XN-2) edge node[above] {$\id_{X}$} (XN-1)
(XN-1) edge node[above] {0} (XN)
(XN) edge node[above] {0} (P2);
}
\end{equation*}
defines an $N$-complex that we denote by \textbf{$\mu^i_N(X)$}. Here the numbers in the top row denote the degree of the corresponding object.
\item Let $g: M \to N$ be a morphism in $\C$, then $\mu_N^i(g):\mu_N^i(M) \to \mu_N^i (N)$ is the morphism of $N$-complexes in $\C$ which is levelwise $(\mu_N^i(g))^j:= g$ for $i-N < j\le i$ and the zero morphism otherwise.


\end{enumerate}
\end{definition}

\begin{example}
Consider the $\CCC$-vectorspace  $\CCC$. Then we have that $\mu_4^7(\CCC)$ is the following $4$-complex
\begin{equation*}
\tikz[heighttwo,xscale=1.7,yscale=2,baseline]{

\node (P0) at (0,0){$\cdot \cdot \cdot$};
\node (X3) at (1,0){0};
\node (X4) at (2,0){$\CCC$};
\node (X5) at (3,0){$\CCC$};
\node (X6) at (4,0){$\CCC$};
\node (X7) at (5,0){$\CCC$};
\node (X8) at (6,0){$0$};
\node (P1) at (7,0){$\cdot \cdot \cdot$};

\node (P3) at (0,0.5){\tiny{$\cdot \cdot \cdot$}};
\node (i-n-1) at (1,0.5){\tiny{ $3$}};
\node (i-n) at (2,0.5){\tiny{$4$}};
\node (i-n+1) at (3,0.5){\tiny{$5$}};
\node (P4) at (4,0.5){\tiny{$6$}};
\node (i-1) at (5,0.5){\tiny{$7$}};
\node (i) at (6,0.5){\tiny{$8$}};
\node (P5) at (7,0.5){\tiny{$\cdot \cdot \cdot$}};

\draw[->]
(P0) edge node[above] {0} (X3)
(X3) edge node[above] {$0$} (X4)
(X4) edge node[above] {$\id_{X}$} (X5)
(X5) edge node[above] {$\id_{X}$} (X6)
(X6) edge node[above] {$\id_{X}$} (X7)
(X7) edge node[above] {$0$} (X8)
(X8) edge node[above] {$0$} (P1);
}
\end{equation*}
in the category of $\CCC$-vectorspaces.
\end{example}

\begin{proposition}\label{functoriality mu}
Let $\C$ be a pointed category and $N\in \N$. Then the assignments $M \mapsto \mu_N^i(M)$ and $ g \mapsto \mu_N^i(g)$ define functors $\C \to \CC_N(\C)$ for each $i \in \Z$. \\
Furthermore \begin{enumerate}
\item the functor $(\_)^i:\CC_N(\C)\to \C$ is left adjoint to $\mu_N^i:\C \to \CC_N(\C)$,
\item the functor $\mu_N^i : \C \to \CC_N(\C)$ is left adjoint to $(\_)^{i-N+1}:\CC_N(\C)\to \C$.
\end{enumerate}
\begin{proof}

\begin{enumerate}
\item Let $X \in \CC_N(\C)$ and $M \in \C$. If $f:X^i \to M$ is a morphism in $\C$, then 
\begin{equation*}
\tikz[heighttwo,xscale=2,yscale=2,baseline]{

\node (P10) at (0,1){$\cdot \cdot \cdot$};
\node (X1) at (1,1){$X^{i-N}$};
\node (X2) at (2,1){$X^{i-N+1}$};
\node (X3) at (3,1){$X^{i-N+2}$};
\node (P11) at (4,1){$\cdot \cdot \cdot$}; 
\node (XN-2) at (5,1){$X^{i-1}$};
\node (XN-1) at (6,1){$X^i$};
\node (XN) at (7,1){$X^{i+1}$};
\node (P12) at (8,1){$\cdot \cdot \cdot$};

\node (P00) at (0,0){$\cdot \cdot \cdot$};
\node (M1) at (1,0){$0$};
\node (M2) at (2,0){$M$};
\node (M3) at (3,0){$M$};
\node (P01) at (4,0){$\cdot \cdot \cdot$}; 
\node (MN-2) at (5,0){$M$};
\node (MN-1) at (6,0){$M$};
\node (MN) at (7,0){$0$};
\node (P02) at (8,0){$\cdot \cdot \cdot$};

\draw[->]
(P10) edge node[above] {$d_X$} (X1)
(X1) edge node[above] {$d_X$} (X2)
(X2) edge node[above] {$d_X$} (X3)
(X3) edge node[above] {$d_X$} (P11)
(P11) edge node[above] {$d_X$} (XN-2)
(XN-2) edge node[above] {$d_X$} (XN-1)
(XN-1) edge node[above] {$d_X$} (XN)
(XN) edge node[above] {$d_X$} (P12)

(P00) edge node[above] {0} (M1)
(M1) edge node[above] {0} (M2)
(M2) edge node[above] {$\id_{M}$} (M3)
(M3) edge node[above] {$\id_{M}$} (P01)
(P01) edge node[above] {$\id_{M}$} (MN-2)
(MN-2) edge node[above] {$\id_{M}$} (MN-1)
(MN-1) edge node[above] {0} (MN)
(MN) edge node[above] {0} (P02)

(X1) edge node[right] {0} (M1)
(X2) edge node[right] {$f \circ d^{N-1}$} (M2)
(X3) edge node[right] {$f \circ d^{N-2}$} (M3)

(XN-2) edge node[right] {$f\circ d$} (MN-2)
(XN-1) edge node[right] {$f$} (MN-1)
(XN) edge node[right] {0} (MN);
}
\end{equation*}
is a morphism $\widehat{f}:X \to \mu^i_N(M)$ in $\CC_N(\C)$ with $\widehat{f}^i=f$.

Conversely every $g:X  \to \mu^i_N(M)$ in $\CC_N(\C)$ is of this form for a unique $f:X^i \to M$, namely $f=g^i$. This establishes a bijection which is natural in $X$ and $M$.

\item
 Let $X \in \CC_N(\C)$ and $M \in \C$. If $f:M \to X^{i-N+1}$ is a morphism in $\C$, then 
\begin{equation*}
\tikz[heighttwo,xscale=2,yscale=2,baseline]{

\node (P10) at (0,0){$\cdot \cdot \cdot$};
\node (X1) at (1,0){$X^{i-N}$};
\node (X2) at (2,0){$X^{i-N+1}$};
\node (X3) at (3,0){$X^{i-N+2}$};
\node (P11) at (4,0){$\cdot \cdot \cdot$}; 
\node (XN-2) at (5,0){$X^{i-1}$};
\node (XN-1) at (6,0){$X^i$};
\node (XN) at (7,0){$X^{i+1}$};
\node (P12) at (8,0){$\cdot \cdot \cdot$};

\node (P00) at (0,1){$\cdot \cdot \cdot$};
\node (M1) at (1,1){$0$};
\node (M2) at (2,1){$M$};
\node (M3) at (3,1){$M$};
\node (P01) at (4,1){$\cdot \cdot \cdot$}; 
\node (MN-2) at (5,1){$M$};
\node (MN-1) at (6,1){$M$};
\node (MN) at (7,1){$0$};
\node (P02) at (8,1){$\cdot \cdot \cdot$};

\draw[->]
(P10) edge node[above] {$d_X$} (X1)
(X1) edge node[above] {$d_X$} (X2)
(X2) edge node[above] {$d_X$} (X3)
(X3) edge node[above] {$d_X$} (P11)
(P11) edge node[above] {$d_X$} (XN-2)
(XN-2) edge node[above] {$d_X$} (XN-1)
(XN-1) edge node[above] {$d_X$} (XN)
(XN) edge node[above] {$d_X$} (P12)

(P00) edge node[above] {0} (M1)
(M1) edge node[above] {0} (M2)
(M2) edge node[above] {$\id_{M}$} (M3)
(M3) edge node[above] {$\id_{M}$} (P01)
(P01) edge node[above] {$\id_{M}$} (MN-2)
(MN-2) edge node[above] {$\id_{M}$} (MN-1)
(MN-1) edge node[above] {0} (MN)
(MN) edge node[above] {0} (P02)

(M1) edge node[right] {0} (X1)
(M2) edge node[right] {$ f  $} (X2)
(M3) edge node[right] {$d \circ f $} (X3)

(MN-2) edge node[right] {$d^{N-2} \circ f$} (XN-2)
(MN-1) edge node[right] {$d^{N-1} \circ f$} (XN-1)
(MN) edge node[right] {0} (XN);
}
\end{equation*}
is a morphism $\widehat{f}:\mu_N^{i}(M)\to X$ in $\CC_N(\C)$ with $\widehat{f}^{i-N+1}=f$.

Conversely every $g:\mu_N^{i}(M) \to X$ in $\CC_N(\C)$ is of this form for a unique morphism $f:M \to X^{i-N+1} $, namely $f=g^{i-N+1}$. This establishes a bijection which is natural in $X$ and $M$.\qedhere
\end{enumerate}
\end{proof}
\end{proposition}

\begin{proposition}
If $\A$ is an additive category, then $\mu^i_N:\A \to \CC_N(\A)$ defines an additive functor, and all the above adjunctions become adjunctions of additive functors.
\begin{proof}
For $f,g \in \A(X,Y)$ it is clear by the definition of $\mu^i_N$ that $\mu_N^i(f+g)=\mu_N^i(f)+ \mu_N^i(g)$. Moreover the bijections from Proposition~\ref{functoriality mu} are compatible with this addition.
\end{proof}
\end{proposition}

After having defined these standard complexes we will inspect which properties the category of $N$-complexes inherits from the ground category. We are especially intrested in limits and colimits since we need coproducts for defining our monoidal structure and colimits for the model structure.


\begin{proposition}\label{complete/cocomplete}
Let $\C$ be a category and assume that $\C$ admits all colimits (respectively limits). Then $\CC_\infty(\C)$ admits all colimits (respectively limits). If $\C$ is pointed, then $\CC_N(\C)$ admits colimits and limits if $\C$ does.
\begin{proof}
The claim for $\CC_\infty(\C)$ follows from the fact that $\CC_\infty(\C)$ is the category of functors from the small category $\I$ to the category $\C$ by Lemma~\ref{1.2}. In particular colimits and limits in $\CC_\infty(\C)$ are computed pointwise.

In the pointed case we need to show that the pointwise construction of limits and colimits defines an object in $\CC_N(\C)$. Hence consider $X:=\colim_{j\in J}X_j$. Now $d_X= \colim_{j\in J} d_{X_j}$, but this implies by pointwise calculation that $d_X^N= \colim_{j\in J} d_{X_j}^N$ and by $d_{X_j}^N=0$ for all $j \in J$ we get $d_X^N=0$. So $\colim_{j\in J}X_j$ is again an object in $\CC_N(\C)$. 
An analogous argument shows the claim for limits.
\end{proof}
\end{proposition}

\begin{proposition}
Let $\A$ be an abelian category, then $\CC_\infty(\A)$ and $\CC_N(\A)$ are abelian as well.
\begin{proof}
The category $\CC_\infty(\A)$ is abelian since it can be interpreted by Lemma~\ref{1.2} as the functor category $\Fun(\I,\A)$ where kernels and cokernels are computed pointwise.
Now to show that $\CC_N(\A)$ is abelian it suffices to show that the kernels and cokernels in $\CC_\infty(\A)$ of maps in the subcategory $\CC_N(\A)$ are again objects in $\CC_N(\A)$. Consider for this a morphism $\varphi : X \to Y$ in $\CC_N(\A)$.
\begin{description}
\item[\normalfont kernels]Let $f: K \to X$ be the kernel of $\varphi$ in $\CC_\infty(\A)$. Then $f$ is a pointwise monomorphism. Furthermore we have $f\circ d_K=d_X \circ f$, in particular one gets 
\begin{equation*}
f \circ d_K^N=d_X^N \circ f =0 \circ f =0.
\end{equation*}
Since $f$ is a pointwise monomorphism this implies $d_K^N=0$. So $K$ is an object of $\CC_N(\A)$ as desired.
\item[\normalfont cokernels] Let $g:Y\to C$ be the cokernel of $\varphi$ in $\CC_\infty(\A)$, then $g$ is a pointwise epimorphism and we have again that $d_C \circ g = g \circ d_Y$, which yields
\begin{align*}
d_C ^N \circ g = g \circ d_Y^N =g \circ 0 =0.
\end{align*}
Now we get by the cancelation property that $d_C^N=0$, and so $C$ is an object of $\CC_N(\A)$.
\end{description}
All together we have that $\CC_N(\A)$ is abelian.
\end{proof}
\end{proposition}

%% file: interlude_q-binomial.tex
\section{Interlude on combinatorics of $q$-binomial coefficients}\label{section 2}

\subsection{Background}
Sign rules play a very important role in homological algebra, and hence also in the theory of chain complexes. In order to extend these to $N$-complexes for arbitrary natural numbers $N$ we study $q$-binomial coefficients. These are of interest, since they arise as coefficients if we apply the Leibniz rule with $-1$ replaced by a generic $q$. We also show that these vanish for primitive roots of unity, which is exactly what we desire to define a monoidal structure on $N$-complexes.

A lot of the results in this subsection are known in the general literature. However, we still show the proofs as these show the intimate relation to binomial coefficients and to supply an easy accessible reference for these results usually given as exercises.

\begin{definition}In $\Z[q]$ define:
\begin{enumerate}
\item The \textbf{$q$-bracket} of $n \in \N$ for $1 < n$ as
\begin{align*}[n]_q:=1+q+q^2+...+q^{n-1}\text{ and }[0]_q:=0,[1]_q:=1.
\end{align*}
\item The \textbf{$q$-factorial} of $n\in \N$ as
\begin{align*}[n]!_q:=[1]_q[2]_q...[n]_q \text{ and } [0]!_q:=1.
\end{align*}
\end{enumerate}
And in $\Q(q)$ define:
\begin{enumerate}[resume]
\item The \textbf{$q$-binomial coefficient} or \textbf{Gaussian binomial coefficient} of $n$ over $k$ for $n,k \in \N$ where $k \le n$ as
\begin{align*}\binom{n}{k}_q:=\frac{[n]!_q}{[n-k]!_q[k]!_q}.
\end{align*}
\end{enumerate}
\end{definition}

\begin{remark}
We will prove in Lemma~\ref{q binomial formula}  that actually $\binom{n}{k}_q$ defines an element in $\N[q]$, i.\,e.\ a polynomial with only natural coefficients.
\end{remark}

\begin{example}
Some non-trivial examples for $q$-binomial coefficients are:
\begin{align*}
\binom{4}{2}_q&=q^4+q^3+2q^2+q+1,\\
\binom{5}{2}_q&=q^6 + q^5 + 2 q^4 + 2 q^3 + 2 q^2 + q + 1,\\
\binom{6}{3}_q&=q^9 + q^8 + 2 q^7 + 3 q^6 + 3 q^5 + 3 q^4 + 3 q^3 + 2 q^2 + q + 1.
\end{align*}
\end{example}
 
 \begin{lemma}\label{q binomial formula} 
 The following equations hold for all $n \in \N$ and $1 \le k \le n-1$:
 \begin{align}
 \binom{n}{k}_q&= q^k \binom{n-1}{k}_q + \binom{n-1}{k-1}_q,\\
\binom{n}{k}_q&=  \binom{n-1}{k}_q + q^{n-k}\binom{n-1}{k-1}_q.
 \end{align}
 \begin{proof}
 \begin{enumerate}
\item We calculate explicitly:
 \begin{align*}
q^k \binom{n-1}{k}_q +\binom{n-1}{k-1}_q &= q^k \frac{[n-1]!_q}{[n-k-1]!_q[k]!_q} +\frac{[n-1]!_q}{[n-k]!_q[k-1]!_q} \\
&= q^k \frac{[n-1]!_q}{[k]_q[n-k-1]!_q[k-1]!_q} +\frac{[n-1]!_q}{[n-k]_q[n-k-1]!_q[k-1]!_q}\\
&= \frac{q^k[n-k]_q[n-1]!_q + [k]_q[n-1]!_q}{[k]_q[n-k]_q[n-k-1]!_q[k-1]!_q}\\
&=\frac{\left( q^k[n-k]_q + [k]_q \right)[n-1]!_q}{[n-k]!_q[k]!_q}\\
&=\frac{\left( q^k\sum_{i=0}^{n-k-1} q^i + \sum_{j=0}^{k-1} q^j \right)[n-1]!_q}{[n-k]!_q[k]!_q}\\
&=\frac{\left( \sum_{i=k}^{n-1} q^i + \sum_{j=0}^{k-1} q^j \right)[n-1]!_q}{[n-k]!_q[k]!_q}\\
&=\frac{\left( \sum_{i=0}^{n-1} q^i \right)[n-1]!_q}{[n-k]!_q[k]!_q}\\
&=\frac{[n]_q[n-1]!_q}{[n-k]!_q[k]!_q} =\frac{[n]!_q}{[n-k]!_q[k]!_q}= \binom{n}{k}_q.
 \end{align*}
 \item Again explicit calculation yields:
  \end{enumerate}
 \begin{align*}
 \binom{n-1}{k}_q +q^{n-k} \binom{n-1}{k-1}_q &= \frac{[n-1]!_q}{[n-k-1]!_q[k]!_q} +q^{n-k} \frac{[n-1]!_q}{[n-k]!_q[k-1]!_q} \\
&= \frac{[n-1]!_q}{[k]_q[n-k-1]!_q[k-1]!_q} +q^{n-k}\frac{[n-1]!_q}{[n-k]_q[n-k-1]!_q[k-1]!_q}\\
&= \frac{[n-k]_q[n-1]!_q + q^{n-k}[k]_q[n-1]!_q}{[k]_q[n-k]_q[n-k-1]!_q[k-1]!_q}\\
&=\frac{\left( [n-k]_q + q^{n-k}[k]_q \right)[n-1]!_q}{[n-k]!_q[k]!_q}\\
&=\frac{\left(\sum_{i=0}^{n-k-1} q^i + q^{n-k}\sum_{j=0}^{k-1} q^j \right)[n-1]!_q}{[n-k]!_q[k]!_q}\\
&=\frac{\left( \sum_{i=0}^{n-k-1} q^i + \sum_{j=n-k}^{n-1} q^j \right)[n-1]!_q}{[n-k]!_q[k]!_q}\\
&=\frac{\left( \sum_{i=0}^{n-1} q^i \right)[n-1]!_q}{[n-k]!_q[k]!_q}\\
&=\frac{[n]_q[n-1]!_q}{[n-k]!_q[k]!_q} =\frac{[n]!_q}{[n-k]!_q[k]!_q}= \binom{n}{k}_q. \qedhere
 \end{align*}

 \end{proof}
 \end{lemma}
\begin{remark}
As mentioned before the above lemma in particular proves by induction $\binom{n}{k}_q \in \N[q]$, since we have $\binom{n}{0}_q=\binom{1}{0}_q =1=\binom{1}{1}_q=\binom{n}{n}_q$ and $\binom{n}{1}_q=\binom{n}{n-1}_q=[n]_q$ for all $n \in \N$.
\end{remark}

\begin{definition}
For $R$ a ring and $\xi \in R$, we denote the images of the above elements under the evaluation homomorphism 
\begin{align*}
ev_\xi: \Z[q] &\to R
\\q &\mapsto \xi
\end{align*}
 by $[n]_\xi, [n]!_\xi$ and $\binom{n}{k}_\xi$.
\end{definition} 
 
The following Lemma will be of use in multiple calculations and generalizes the ordinary binomial formula to a $q$-binomial formula:
 
 \begin{lemma}\label{q binom}
The following formula holds in $\Z[x,y,q]$:
$$\prod_{k=0}^{n-1}(x+q^ky)=\sum_{k=0}^{n}q^{\frac{k(k-1)}{2}}\binom{n}{k}_q x^{n-k}y^k .$$
\begin{proof}
We use induction. For the induction base consider the case $n=1$, which leads 
\begin{align*}
\prod_{k=0}^{0}(x+q^ky)&= (x+y)=\sum_{k=0}^{1}q^{\frac{k(k-1)}{2}}\binom{n}{k}_q x^{n-k}y^k .
\intertext{Now the claim follows by the following calculation:}
\prod_{k=0}^{n}\left(x+q^k y \right)&=\left(\prod_{k=0}^{n-1}\left(x+q^k y \right)\right)\left(x+q^ny\right)\\
&=\left(\sum_{k=0}^{n}q^{\frac{k(k-1)}{2}}\binom{n}{k}_q x^{n-k}y^k\right)\left(x+q^n y\right)\\
&=\sum_{k=0}^{n}q^{\frac{k(k-1)}{2}}\binom{n}{k}_q x^{n-k+1}y^k + \sum_{k=0}^{n}q^{\frac{k(k-1)}{2}+n}\binom{n}{k}_q x^{n-k}y^{k+1}\\
&=\sum_{k=0}^{n}q^{\frac{k(k-1)}{2}}\binom{n}{k}_q x^{n-k+1}y^k + \sum_{k=1}^{n+1}q^{\frac{(k-1)(k-2)}{2}+n}\binom{n}{k-1}_q x^{n-k+1}y^{k}\\
&=x^{n+1} + q^{\frac{n(n-1)}{2}+n}y^{n+1} +\sum_{k=1}^{n}\left(q^{\frac{k(k-1)}{2}}\binom{n}{k}_q + q^{\frac{(k-1)(k-2)}{2}+n}\binom{n}{k-1}_q\right) x^{n-k+1}y^{k}\\
&=x^{n+1} + q^\frac{(n+1)n}{2}y^{n+1} +\sum_{k=1}^{n}\left(q^{\frac{k(k-1)}{2}}\binom{n}{k}_q + q^{\frac{k(k-1)}{2}+n-k+1}\binom{n}{k-1}_q\right)  x^{n-k+1}y^{k}\\
&=x^{n+1} + q^\frac{(n+1)n}{2}y^{n+1}+\sum_{k=1}^{n}q^{\frac{k(k-1)}{2}}\left(\binom{n}{k}_q + q^{n-k+1}\binom{n}{k-1}_q\right) x^{n-k+1}y^{k} \\
&=x^{n+1} + q^\frac{(n+1)n}{2}y^{n+1}+\sum_{k=1}^{n}q^{\frac{k(k-1)}{2}}\binom{n+1}{k}_q x^{n-k+1}y^{k} \\
&=\sum_{k=0}^{n+1}q^{\frac{k(k-1)}{2}}\binom{n+1}{k}_q x^{n-k+1}y^{k}. \\
\end{align*}
Here we used in the third to last equality formula (2) from Lemma~\ref{q binomial formula}.
\end{proof}
\end{lemma}
 
\subsection{Vanishing at roots of unity}
Now we will show that the $q$-binomial coefficients vanish at certain primitive roots of unity. A very powerful approach for this is to use cyclic sieving. Since this needs a more involved argument and we only need a vanishing argument for $q$-binomial coefficients, we will present a less technical proof. However our proof is strongly inspired by cyclic sieving and suffices for our purposes. For the cyclic sieving argument we refer to \cite{cyclic}.

We first need to review some standard representation theory notions to apply them in the vanishing theorem.

\begin{definition}
Let $n,k \in \N$ such that $0 \le k \le n$. Then
\begin{enumerate}
\item let the set $\left\{1,2,...,n \right\}$ be denoted by $[n]$,
\item let $\langle(1,2,...,n)\rangle :=C_n\subset S_n$ denote the \textbf{cyclic group of order $n$} identified as the subgroup of the symmetric group generated by $(1,2,...,n)$.
\end{enumerate} 
\end{definition}

\begin{definition}
Let $G$ be a group, $V$ a finite dimensional vector space over a field $\kk$ and $\varphi: G \to \Aut(V)$ an action of $G$ on $V$. Then the \textbf{character} $\chi_\varphi: G \to \kk $ of $\varphi$ is defined by $\chi_\varphi(g)=\tr(\varphi( g))$ for $g \in G$.
\end{definition}

\begin{remark}
This is well defined since the trace is basis independent
\end{remark}

\begin{remark} We will denote for a ring $R$ and a set $I$ the free $R$-module over $I$ by $\langle I \rangle_R$.
\end{remark}

\begin{definition}
Let $n \in \N$. Then define the action $\varphi_{n}$ of $C_n:= \langle (1,2,...,n)\rangle \subset S_n$  on $\langle[n]\rangle_\CCC$ as
$$ \varphi_{n}(1,2,...,n)(i)=(1,2,...,n)(i). $$
\end{definition}

\begin{definition}
Let $V$ be a vector space over a field $\kk$. Then denote for $k \in \N$ the \textbf{$k$-th exterior product of $V$} by $\Lambda^kV$. 

Furthermore define for a group representation $\varphi:G \to \Aut({V,V})$ the induced representation 
\begin{align*}
\Lambda^k\varphi: G &\to \Aut_\kk(\Lambda^kV,\Lambda^kV)\\
g &\mapsto  (\Lambda^k\varphi(g) :\Lambda^kV \to \Lambda^kV)
\end{align*}
\end{definition}

\begin{proposition}\label{2.11}
Let $n,k \in \N$ such that $1\le k <n$. Then we have 
$$0=\chi_{\Lambda^k\varphi_{n}}(1,2,...,n).$$ 
 \begin{proof}
 Consider the basis $B:=\{e_{\{i_1,i_2,...,i_k\}}:= i_1\wedge i_2 \wedge ... \wedge i_k \mid 1 \le i_1 < i_2 <...<i_k \le n\}$ of $\Lambda^k \langle[n]\rangle_\CCC $. Then we get by definition 
 \begin{align*} 
 \Lambda^k \varphi_{n}(1,2,...,n)(e_{\{i_1,i_2,...,i_k\}}) &= \Lambda^k\varphi_{n}(1,2,...,n)(i_1\wedge i_2 \wedge ... \wedge i_k)\\
 &=(1,2,...,n)i_1\wedge (1,2,...,n)i_2 \wedge ... \wedge (1,2,...,n)i_k  \\
 &= \delta e_{\{(1,2,...,n)i_1,(1,2,...,n)i_2,...,(1,2,...,n)i_k\}},
 \end{align*}
 where $\delta = \pm 1$. In particular the matrix of  $\Lambda^k \varphi_{n}(1,2,...,n)$ with respect to $B$ has only zeroes on the diagonal. So we get 
 $$\chi_{\Lambda^k\varphi_{n}}(1,2,...,n)=\tr\left(\Lambda^k \varphi_{n}(1,2,...,n)\right)=0$$
 as claimed.\qedhere
 \end{proof}
 \end{proposition}

\begin{proposition}\label{2.15}
Let $V$ be an $n$-dimensional vector space over a field $\kk$ and $f: V \to V$ a diagonalizable endomorphism with eigenvalues $x_1,x_2,...,x_n$. Then we have for $k \in \N$ that the endomorphism   $$\Lambda^k(f): \Lambda^k(V) \to \Lambda^k(V)$$ induced by $f$ is diagonalizable and we have
\begin{equation*}
\tr(\Lambda^k(f))= \sum_{\mathclap{1 \le i_1 <i_2<...<i_k \le n}} x_{i_1} x_{i_2}...x_{i_k}.
\end{equation*}
\begin{proof}
We have by assumption a basis $B:=\{b_1,b_2,...,b_n\}$ of $V$ consisting of eigenvectors for the respective eigenvalues of $f$. This basis induces a basis $\Lambda^k B$ on  $\Lambda^k(V)$. We may consider $\Lambda^k(f)$ with respect to this basis and get 
\begin{align*}
\Lambda^k(f)\left(b_{i_1}\wedge b_{i_2} \wedge ... \wedge  b_{i_k}\right)=x_{i_1}b_{i_1}\wedge  x_{i_2}b_{i_2}\wedge  ... \wedge  x_{i_k}b_{i_k}=x_{i_1}x_{i_2} ...  x_{i_k}b_{i_1}\wedge  b_{i_2}\wedge  ... \wedge  b_{i_k}.
\end{align*}
So the vectors in $\Lambda^k B$ are eigenvectors of $\Lambda^k(f)$ as well. Using this basis to calculate $\tr(\Lambda^k(f))$ gives 
\begin{align*}
\tr(\Lambda^k(f))= \sum_{\mathclap{1\le i_1 <i_2<...<i_k \le n}} x_{i_1} x_{i_2}...x_{i_k}.
\end{align*}
This finishes the proof. \qedhere
\end{proof}
\end{proposition}

\begin{definition}
Let $n,k \in \N$, such that $0\le k \le n$. We denote by $e_{k,n}(x_1,...,x_n)$ the \textbf{$k$-th elementary symmetric polynomial in  $n$ variables}. This is defined as 
\begin{align*}
e_{k,n}(x_1,...,x_n):=&\sum_{\mathclap{1 \le i_1 <i_2<...<i_k \le n}} x_{i_1} x_{i_2}...x_{i_k}\in \Z[x_1,x_2,...,x_n],
\intertext{for $0<k$ and for $k=0$ as}e_{0,n}(x_1,...,x_n):=&1.
\end{align*}
\end{definition}

After this reminder on some linear algebra and representation theory the only thing missing is a connection to $q$-binomials. For this we consider the following relation between $e_{k,n}(1,...,q^{n-1})$ and $\binom{n}{k}_q$.

\begin{proposition}\label{2.17}
Let $n,k\in \N$ such that $k\le n$ and $1 \le n$. Then the equation 
\begin{equation*}
e_{k,n}(1,q,q^2,...,q^{n-1})=q^{\frac{k(k-1)}{2}}\binom{n}{k}_q
\end{equation*}
holds in $\Z[q]$.
\begin{proof}
This follows from Lemma~\ref{q binom} by considering $x=1$ and $y=t$. More precisely we get the following sequence of identities, using for the second identity Vieta's formula:
\begin{align*}
\sum_{k=0}^{n}q^{\frac{k(k-1)}{2}}\binom{n}{k}_q t^k = \prod_{k=0}^{n-1}\left( 1 + q^k t\right)= \sum_{k=0}^{n}e_{k,n}(1,q,...,q^{n-1})t^k.
\end{align*}
This gives the desired equation by comparing coefficients.
\end{proof}
\end{proposition}

Now we just need to put the collected statements together and get:

\begin{theorem} \label{q-binomial vanishing at root of unity in C}
Let $n,k \in \N$ such that $0<k<n$ and let $\xi\in \CCC$ be a primitive $n$-th root of unity. Then we have in $\CCC$ $$\binom{n}{k}_\xi=0.$$
\begin{proof}
First of all we know $0 = \chi_{\Lambda^k\varphi_{n}}(1,2,...,n)$ by Proposition~\ref{2.11}. Now $\varphi_{n}(1,2,...,n)$ has characteristic polynomial $x^n-1$ and so it has a basis of eigenvectors $\{v_1,...,v_n\}$, where $v_i$ has the eigenvalue $\xi^i$.  In particular it is diagonalizable. Now by Proposition~\ref{2.15} also $\Lambda^k(\varphi_{n})(1,2,...,n)$ is diagonalizable and we get:
\begin{align*}
0 &= \chi_{\Lambda^k(\varphi_{n})}(1,2,...,n)=\sum_{\mathclap{1 \le i_1 <i_2<...<i_k \le n}} \xi^{i_1} \xi^{i_2}...\xi^{i_k}=e_{k,n}(\xi,\xi^2,...,\xi^{n-1},\xi^n).
\intertext{We get by $\xi^n=1=\xi^0$, Proposition~\ref{2.17} and the symmetry of the elementary symmetric polynomials that}
0&=e_{k,n}(\xi,\xi^2,...,\xi^{n-1},\xi^n)=e_{k,n}(1,\xi,\xi^2,...,\xi^{n-1})=\xi^{\frac{k(k-1)}{2}}\binom{n}{k}_\xi.
\end{align*}
Since $\CCC $ is a field and $\xi \neq 0$ we get $\binom{n}{k}_\xi=0$ as claimed.
\end{proof}
\end{theorem}

The previous statement about $q$-binomial coefficients relies heavily on the particular structure of complex numbers.
Importantly, standard arguments, such as the cyclotomic polynomial \cite{lang}, allow to carry it over to arbitrary integral domains:

\begin{lemma}\label{q binomial vanishing at root of unity}
Let $R$ be an integral domain, $n \in \N$ and $\xi \in R$ a primitive $n$-th root of unity i.\,e.\ an element such that $\xi^n=1$ and $\xi^j\neq 1$ for all $0<j<n$. Then we have $\binom{n}{k}_\xi=0$ for all $0 < k < n$.
\begin{proof}
We will prove this statement in two steps, firstly showing the statement for $x \in \Z[x]/\left(\Phi_n\right)$, where $\Phi_n$ is the $n$-th cyclotomic polynomial, and in a second step that the evaluation morphism $\Z[x] \to R$ defined by $x \mapsto \xi$ factors through $\Z[x]/\left(\Phi_n\right)$.

Consider the map $\Z[x]/\left(\Phi_n\right) \to \CCC$ defined by $x \mapsto \xi$ for $\xi$ a primitive $n$-th root of unity. This is a well defined monomorphism since $\Phi_n$ is the $n$-th cyclotomic polynomial. We have $\binom{n}{k}_\xi=0$ for all $0< k < n$ by Theorem~\ref{q-binomial vanishing at root of unity in C} since this is just the image of $\binom{n}{k}_x \in \Z[x]/\left(\Phi_n\right) $ we get the desired formula in 
$\Z[x]/\left(\Phi_n\right)$.

Now consider for an arbitrary integral domain the following composition denoting by $\Quot(R)$ and $\overline{\Quot(R)}$ the quotient field of $R$ and its algebraic closure:

\begin{equation*}
\tikz[heighttwo,xscale=2,yscale=2,baseline]{
\node (ZX) at (0,1){$\Z[x]$};
\node (ZXmodphi) at (0,0){$\Z[x]/\left(\Phi_n\right)$};
\node (R) at (1.1,1){$R$};
\node (QR) at (2.1,1){$ \Quot(R)$};
\node (QRclosure) at (3.45,1){$\overline{\Quot(R)},$};

\draw[->]
(ZX) edge node[above]{$X \mapsto \xi$}(R);
\draw[right hook->]
(R) edge (QR)
(QR) edge (QRclosure);
\draw[->>]
(ZX)edge (ZXmodphi);
\draw[dashed, ->]
(ZXmodphi) edge (R);
}
\end{equation*}
where we only need to show that the dashed arrow exists such that the diagram commutes. By the homomorphism theorem it suffices to check that $\Phi_n(\xi)=0$ in $R$. This is equivalent to $\Phi_n(\xi)=0$ in $\overline{\Quot(R)}$, which holds since $\Phi_n$ is the $n$-th cyclotomic polynomial and $\xi$ is an $n$-th primitive root of unity in $R$ and hence also in $\overline{\Quot(R)}$.
\end{proof}
\end{lemma}

\begin{remark}
This statement heavily uses the cyclotomic polynomial $\Phi_n$ and its properties. In particular that for every field $\kk$ with $\mathrm{char}(\kk)\nmid n$ the $n$-th primitive roots of unity are roots of $\Phi_n$, and that $\Phi_n$ can be defined in $\Z[x]$ and hence is independent of the field. For the precise arguments that show this we refer to Lang's book \cite{lang}.
\end{remark}

%% file: monoidal_structure.tex
\section{Lifting functors}\label{section 3}

In this section we will apply the tools developed in the previous sections to define functors on $\CC_N(\C)$ using functors on $\C$. This is of particular interest since a lot of structures on a category can be encoded using functors. First we will consider ordinary functors, but later we will discuss how to lift bifunctors which will be essential for constructing the monoidal structure on the category of $N$-complexes.

\begin{definition}Let $\C,\D$ be categories and $F: \C \to \D$ a functor. Then define for an $\infty$-complex $X$ in $\C$ the $\infty$-complex $F(X)$ in $\D$ as
\begin{align*}
F(X)^i &:= F(X^i) \\
(d_{F(X)} :F(X^i) \to F(X^{i+1}) ) &:=  F(d_X:X^i \to X^{i+1}). 
\end{align*}
Furthermore define for a morphism of $\infty$-complexes $f:X \to Y$, the family of morphisms $(F(f)^i)_{i \in \Z}$ as
\begin{align*}
F(f)^i= F(f^i).
\end{align*}
\end{definition}  
\begin{proposition}\label{induced functor}
The above assignments define a functor $F: \CC_\infty(\C) \to \CC_\infty(\D)$.
\begin{proof}
By construction $F(X)$ is an $\infty$-complex in $\D$, so it just remains to be checked that for a morphism of $\infty$-complexes $f:X \to Y$, the family $F(f)$ is indeed a morphism of $\infty$-complexes, since functoriality follows by levelwise functoriality. Hence we need to show that the diagram
\begin{equation*}
\tikz[heighttwo,xscale=3,yscale=2,baseline]{
\node (Xi) at (0,1){$ F(X^i)=F(X)^i$};
\node (Xi+1) at (0,0){$F(X^{i+1})=F(X)^{i+1}$};
\node (Yi+1) at (2,0){$F(Y^{i+1})=F(Y)^{i+1}$};
\node (Yi) at (2,1){$F(Y^{i})=F(Y)^i$};

\draw[->]
(Xi) edge node[left]{$d_{F(X)} =F(d_X)$}(Xi+1)
(Yi) edge node[right]{$d_{F(Y)} =F(d_Y)$}(Yi+1)
(Xi) edge node[above]{$ F(f)^i=F(f^i)$}(Yi)
(Xi+1) edge node[below]{$ F(f)^{i+1}=F(f^{i+1}) $}(Yi+1);
}
\end{equation*}commutes for all $i \in \Z$. Since $F: \C \to \D$ is a functor, it suffices to show that 
\begin{equation*}
\tikz[heighttwo,xscale=2,yscale=2,baseline]{
\node (Xi) at (0,1){$ X^i$};
\node (Xi+1) at (0,0){$X^{i+1}$};
\node (Yi+1) at (2,0){$Y^{i+1}$};
\node (Yi) at (2,1){$Y^{i}$};

\draw[->]
(Xi) edge node[left]{$d_X$}(Xi+1)
(Yi) edge node[right]{$d_Y$}(Yi+1)
(Xi) edge node[above]{$f^i$}(Yi)
(Xi+1) edge node[below]{$f^{i+1} $}(Yi+1);
}
\end{equation*} commutes for all $i \in \Z$, but this holds since $f$ is a morphism of $\infty$-complexes in $\C$. 
\end{proof}
\end{proposition}

\begin{corollary}\label{induced pointed functors}
Let $\C,\D$ be pointed categories, $N \in \N$ and $F:\C \to \D$ a pointed functor, i.\,e.\ $F(0)\cong 0$ for $0$ a zero object in $\C$. Then the functor $F: \CC_\infty(\C) \to \CC_\infty(\D)$ induces a functor $\CC_N(\C) \to \CC_N(\D)$.
\begin{proof}
By Proposition~\ref{induced functor} it suffices to check that $F(X)$ is an $N$-complex in $\D$ if $X$ is an $N$-complex in $\C$, but this is clear, since 
$$d_{F(X)}^N=F(d_X^N)=F(0)=0,$$
and so $F(X)$ is an $N$-complex.\qedhere
\end{proof}
\end{corollary}

\begin{proposition}\label{induced adjunction}
Let $\C$,$\D$ be categories and $F:\C \rightleftarrows \D:G$ a pair of adjoint functors. Then $F:\CC_\infty(\C) \rightleftarrows \CC_\infty(\D): G$ is also a pair of adjoint functors. 
\begin{proof}
Since $F,G$ are adjoint functors we have natural transformations $\eta:\id_\C \to G\circ F$ and  $\epsilon: F \circ G \to \id_\D$ such that the following diagrams commute:
\begin{equation*}
\tikz[heighttwo,xscale=2,yscale=2,baseline]{
\node (F1) at (0,1){$F$};
\node (FGF) at (1,1){$FGF,$};
\node (F2) at (1,0){$F$};
\node (G1) at (3,1){$ G$};
\node (GFG) at (4,1){$GFG$};
\node (G2) at (4,0){$G.$};

\draw[->]
(F1) edge node[above]{$F \eta $}(FGF)
(FGF) edge node[right]{$\epsilon F$}(F2)
(G1) edge node[above]{$\eta G$}(GFG)
(GFG) edge node[right]{$G \epsilon $}(G2)
(F1) edge node[below left]{$\id$}(F2)
(G1) edge node[below left]{$\id $}(G2);
}
\end{equation*} Now since $\eta$ and $\epsilon$ are natural transformations, we get that the following squares commute for all $X,Y \in \CC_\infty(\C)$ and $j \in \Z$:
\begin{equation*}
\tikz[heighttwo,xscale=3,yscale=2,baseline]{
\node (X1) at (0,1){$X^j$};
\node (GFX) at (1,1){$G\circ F(X^j)$};
\node (GFY) at (1,0){$G\circ F(X^{j+1})$};
\node (Y1) at (0,0){$X^{j+1}$};
\node (FGY) at (3,0){$ F\circ G (Y^{j+1})$};
\node (FGX) at (3,1){$ F\circ G (Y^j)$};
\node (X2) at (4,1){$Y^j$};
\node (Y2) at (4,0){$Y^{j+1}.$};
\draw[->]
(X1) edge node[above]{$ \eta_{X^j} $}(GFX)
(Y1) edge node[below]{$\eta_{X^{j+1}} $}(GFY)
(X1) edge node[left]{$d_X$}(Y1)
(GFX) edge node[right]{$G\circ F (d_X) $}(GFY)
(FGY) edge node[below]{$ \epsilon_{Y^{j+1}} $}(Y2)
(FGX) edge node[above]{$\epsilon_{Y^{j}} $}(X2)
(FGX) edge node[left]{$F \circ G (d_Y)$}(FGY)
(X2) edge node[right]{$d_Y $}(Y2);
}
\end{equation*}
Hence $(\eta_{X^j}:X^j \to G\circ F (X^j))_{j \in \Z}$ and $(\eta_{Y^j}:F \circ G(Y^j) \to Y^j)_{j \in \Z}$ define morphisms of $\infty$\Hyphdash*complexes. These define by construction natural transformations $\eta_{\CC_\infty}:\id_{\CC_\infty(\C)} \to G \circ F$ and $\epsilon_{\CC_\infty(\C)}: F \circ G \to \id_{\CC_\infty(\D)}$. Furthermore the diagrams
\begin{equation*}
\tikz[heighttwo,xscale=2,yscale=2,baseline]{
\node (F1) at (0,1){$F$};
\node (FGF) at (1,1){$FGF$};
\node (F2) at (1,0){$F$};
\node (G1) at (3,1){$ G$};
\node (GFG) at (4,1){$GFG$};
\node (G2) at (4,0){$G$};
\draw[->]
(F1) edge node[above]{$F \eta_{\CC_\infty} $}(FGF)
(FGF) edge node[right]{$\epsilon_{\CC_\infty} F$}(F2)
(G1) edge node[above]{$\eta_{\CC_\infty} G$}(GFG)
(GFG) edge node[right]{$G \epsilon_{\CC_\infty} $}(G2)
(F1) edge node[below left]{$\id$}(F2)
(G1) edge node[below left]{$\id $}(G2);
}
\end{equation*}
commute levelwise by construction and so they commute in general. This means that $\eta_{\CC_\infty}$ and $\epsilon_{\CC_\infty}$ define a counit and unit for $F$ and $G$ such that the triangle identities hold, which proves the claim. 
\end{proof}
\end{proposition}
\begin{corollary}
Let $\C,\D$ be pointed categories, $N\in \N$ and $F:\C \rightleftarrows \D:G$ an adjunction of functors. Then $F:\CC_N(\C) \rightleftarrows \CC_N(\D): G$ is also a pair of adjoint functors. 
\begin{proof}
By Proposition~\ref{induced adjunction} we get an adjunction $F:\CC_\infty(\C) \rightleftarrows \CC_\infty(\D): G$. Now by Corollary~\ref{induced pointed functors} it suffices to show that $F:\C \to \D$ and $G:\D \to \C$ are actually pointed functors. Consider the following isomorphisms induced by the adjunction for all $X \in \C$ and $Y \in \D$: 
\begin{align*}
\{*\}&\cong \C(0,F(Y))\cong \D(G(0),Y)\\
\{*\}&\cong \D(X,F(0))\cong \D(G(X),0).
\end{align*}
So $G(0)$ is an initial object and $F(0)$ is a final object. Now since $\C$ and $\D$ are pointed we get $G(0)\cong 0$ and $F(0)\cong 0$ which finishes the claim.
\end{proof}
\end{corollary}

The following definition is a general construction which can be used to encode the property of bifunctors respecting an enriched structure in the sense that they commute with the operations coming from this structure. 

\begin{definition}\label{tensorcat}
Let $(\D,\otimes_\D,\unit_\D,\alpha_\D,\lambda_\D,\rho_\D)$ be a symmetric monoidal category and $\C,\C'$ be categories enriched over it. 
 We define \textbf{$\C \otimes_\D \C'$} as the category consisting of objects of the form $(X , Y)$, where $X \in \C,Y \in \C'$, and morphism objects $$\C \otimes_\D \C'\left((X_1,Y_1),(X_2,Y_2)\right):= \C(X_1,X_2)\otimes_\D \C'(Y_1, Y_2).$$

\end{definition}

Now we will need some kind of lifting for bifunctors of the form $\C \times \C' \to \D$  since we want to construct a monoidal structure on the category of $N$-complexes. Similar to the ordinary case of chain complexes we need some kind of twisted Leibniz rule. We say twisted, because the ordinary Leibniz rule just works for $2$-complexes or $N$-complexes over rings with characteristic $2$. This is owing to the fact that $-1$ is a primitive second root of unity, which will be clear by the proof of Proposition~\ref{lifting bifunctor Ncomplex}. Hence we will define our lifted functors using a Leibniz rule which involves arbitrary elements of a ring. Furthermore we will show that if these elements are suitable roots of unity, we indeed get functors of $N$-complexes and not just $\infty$-complexes.

\begin{definition}\label{lifting bifunctor}
Let $R$ be a commutative ring, $\C,\C'$ and $\D$ be categories enriched over $\Mod_R$ such that $\D$ admits countable coproducts, $\xi \in R$, and let $$F:\C \otimes_{\Mod_R} \C' \to \D$$ be an additive $R$-linear functor. Then define for objects $X \in \CC_\infty(\C)$ and $ Y \in \CC_\infty(\C')$ the $\infty$-complex $F_\xi(X,Y)$ in $\D$ by
$$F_\xi(X,Y)^n := \bigoplus_{i+j=n}F(X^i, Y^j)$$
with the differential defined on summands by
$$d_{F_\xi(X,Y)}\circ \iota_{i,j}:=\iota_{i+1,j}\circ F(d_X,\id)+ \xi^i \iota_{i,j+1} \circ F(\id,d_Y).$$
Here $\iota_{k,l}$ denotes the canonical inclusion $$F(X^k,Y^l) \to \bigoplus_{i+j=n}F(X^i,Y^j)$$ for $k+l=n$.

Furthermore for morphisms of $\infty$-complexes $f:X \to X'$ and $g:Y \to Y'$ define the family $\left( F_\xi(f,g)^n \right)_{n\in \Z}$ of morphisms in $\D$ by
$$F_\xi(f,g)^n :=\bigoplus_{i+j=n} F(f^i,g^j).$$
\end{definition}

\begin{proposition}\label{infty bifunctor}
Let $R$ be a commutative ring, $\C,\C'$ and $\D$ be categories enriched over $\Mod_R$ such that $\D$ admits countable coproducts, $\xi \in R$, and let $$F:\C \otimes_{\Mod_R} \C' \to \D$$ be an additive $R$-linear functor. Then the assignments from Definition~\ref{lifting bifunctor} define a functor 
$$F_\xi(\_,\_): \CC_\infty(\C) \otimes \CC_\infty (\C') \to \CC_\infty (\D).$$
\begin{proof}
It suffices to check for $f:X \to Y$ and $g:X' \to Y'$ that  $$F_\xi(f,g)\circ d_{F_\xi(X,X')}=d_{F_\xi(Y,Y')} \circ F_\xi(f,g).$$ Since by definition $F_\xi(\_,\_)$ is functorial as well as it has the right codomain. Hence we check the statement on direct summands of the form $F(X^{i}, X'^{j})$:
\begin{align*}
F_\xi(f,g)^n\circ d_{F_\xi(X,X')}\circ \iota_{i,j}&=\bigoplus_{i'+j'=n}F(f^{i'} , g^{j'})\circ \left(\iota_{i+1,j}\circ F(d_X, \id)  + \xi^{i} \iota_{i,j+1} \circ F(\id , d_{X'})\right)\\
&= \iota_{i+1,j}\circ F(f^{i+1}\circ d_X , g^{j})  + \xi^{i} \iota_{i,j+1} \circ F(f^{i} , g^{j+1} \circ d_{X'})\\
&=\iota_{i+1,j}\circ F(d_Y \circ f^{i} , g^{j})+ \xi^{i} \iota_{i,j+1} \circ  F(f^{i} , d_{Y'} \circ g^{j} )\\
&=\left(\iota_{i+1,j}\circ F(d_Y, \id) + \xi^{i} \iota_{i,j+1} \circ F(\id , d_{Y'})\right)\circ F(f^{i} , g^{j})\\
&=d_{F_\xi(X,X')} \circ\iota_{i,j} \circ  F(f^{i} , g^{j})\\
&=d_{F_\xi(Y,Y')} \circ \left( \bigoplus_{i'+j'=n}F(f^{i'} , g^{j'})\right) \circ \iota_{i,j}\\
&=d_{F_\xi(Y,Y')} \circ F_\xi(f,g) \circ \iota_{i,j}. \qedhere
\end{align*}
\end{proof}
\end{proposition}

Having this general definition we can use our work on $q$-binomial coefficients to prove that $F_\xi(\_,\_)$ indeed defines a bifunctor on the category of $N$-complexes for $\xi$ a suitable root of unity. To make precise what suitable means we define the following:

\begin{definition}\label{dagger}
Let $(\C,N,\xi)$ be a triple consisting of an additive category $\C$ enriched over $\Mod_R$ for $R$ a ring, $N \in \N$ and $\xi \in R$. Then we say that $(\C,N,\xi)$ \textbf{allows the twisted Leibniz rule} if one of the following holds:
\begin{enumerate}[label=(\roman*)]
\item Every morphism in $\C$ is $p$-torsion for $p$ a prime number i.\,e.\ $pf=0$ for all morphisms $f$, $\xi=1$ and $N=p^n$ for some natural number $n\in \N$.
\item The element $\xi \in R$ is a primitve $N$-th root of unity.
\item  Every morphism in $\C$ is $p$-torsion for $p$ a prime number, the element $\xi \in R$ is a primitive $k$-th root of unity and $N=p^nk$ for some natural number $n\in \N$.
\end{enumerate}
\end{definition}

\begin{remark}
In the above definition, case (iii) specializes to case (i) by setting $\xi=1$ and $k=1$.
Furthermore one can think of case (iii) as an enhancement of case (ii), by applying the structure of case (i) to case (ii).
\end{remark}

\begin{theorem}\label{lifting bifunctor Ncomplex}
Let $R$ be an integral domain, $\C,\C'$ and $\D$ be additive categories enriched over $\Mod_R$ such that $\D$ admits countable coproducts, $N \in \N$, $\xi \in R$ and let $F:\C \otimes_{\Mod_R} \C' \to \D$ be an additive bifunctor. Then the functor defined in Definition~\ref{lifting bifunctor} restricts to a functor
$$F_\xi(\_,\_): \CC_N(\C) \otimes \CC_N (\C) \to \CC_N (\C)$$ if $(\D,N,\xi)$ allows the twisted Leibniz rule.
\begin{proof}
Since we already know by Proposition~\ref{infty bifunctor} that $F_\xi(\_,\_):\CC_\infty(\C) \times \CC_\infty(\C') \to \CC_\infty(\D)$ is a functor we only need to show in all cases that for $N$-complexes $X,Y$ we get that $F_\xi(X,Y)$ is an $N$-complex as well. For this we consider the three cases in Definition~\ref{dagger} separately and calculate omitting inclusions:
\begin{enumerate}[label=(\roman*)]
\item Since every morphism in $\D$ is a $p$-torsion for a prime number $p$ and $N=p^n$, we may use that for all morphisms $f$ in $\D$ we have $p f=f+f+...+f=0 \in \D$, hence it suffices to consider $\D$ as enriched over $R/pR$. Now we may use that $\binom{p}{k}=0$ in $R/pR$ for all $0<k<p$. First we generalize this to $\binom{N}{k}=0$ if $0<k < N$.
For this consider the following calculation in $R/pR[x,y]$
\begin{align*}
\sum_{k=1}^{N} \dbinom{N}{k} x^{N-k}y^k&=(x+y)^N=(((x+y)^p)^p...)^p=(((x^p+y^p)^p)^p...)^p\\
&=x^{p^n} + y^{p^n}=x^N+y^N.
\end{align*}
So the claim holds by comparing coefficients. Now using this we get
\begin{align*}
d_{F_1(X,Y)}^N=& d^{N-1}(F(d_X, \id_Y) + F(\id_X, d_Y))\\
=&d^{N-2}(F(d_X^2, \id_Y) + F(d_X, d_Y) + F(d_X, d_Y) + F(\id_X, d^2_Y))\\
=& d^{N-2}(F(d_X^2 , \id_Y) + 2 F( d_X , d_Y)  + F(\id_X , d^2_Y))\\
=&d^{N-3}(F(d_X^3, \id_Y) + 3 F(d_X^2, d_Y) + 3 F(d_X, d_Y^2) + F(\id_X, d^3_Y)\\
&\vdots\\
=&\sum_{k=0}^{N} \dbinom{N}{k} F(d^{N-k}_X, d^k_Y)\\
=&F(d_X^N, \id_Y) + F(\id_X , d_Y^N)\\
=&0
\end{align*}
and so $F_\xi(X, Y)$ is an $N$-complex. Here we used in the second to last equality $\binom{N}{k}=0$ if $0<k < N$.
\item
In this case we will use the two formulas from Proposition~\ref{q binomial formula} and Corollary~\ref{q binomial vanishing at root of unity}.
Here we need for $\binom{N}{k}_\xi=0$ that $\xi$ is a primitive $N$-th root of unity. 
Analogously to Lemma~\ref{q binom} we calculate for each summand (cf. Remark~\ref{calculation conceptual}):
\begin{align*} &\left. d_{F_\xi (X ,Y)}^N \right|_{F(X^i, Y^j)}\\
=& d^{N-1} \left(F(d_X, \id_Y) + \xi^i F(\id_X, d_Y)\middle) \right|_{F(X^i, Y^j)}\\
=&d^{N-2}\left( F(d_X^2, \id_Y) + \xi^{i+1} F(d_X, d_Y) + \xi^i F(d_X , d_Y) + \xi^{2i}F(\id_X , d^2_Y) \middle)\right|_{F(X^i, Y^j)}\\
=&d^{N-2}\left( \binom{2}{2}_{\xi}F(d_X^2, \id_Y)  + \xi^{i} \binom{2}{1}_{\xi}F(d_X , d_Y) +\xi^{2i}\binom{2}{0}_{\xi}F(\id_X, d^2_Y) \middle) \right|_{X^i\otimes Y^j}\\
=&d^{N-3} \left( \binom{2}{2}_{\xi}F(d_X^3 , \id_Y)  +\xi^{i+2}\binom{2}{2}_{\xi}F(d_X^2 , d_Y) +\xi^{i} \binom{2}{1}_{\xi}F(d_X^2 , d_Y) \right.\\
+& \left. \xi^{2i+1} \binom{2}{1}_{\xi}F(d_X , d_Y^2)+\xi^{2i}\binom{2}{0}_{\xi}F(d_X d^2_Y)+\xi^{3i}\binom{2}{0}_{\xi}F(\id_X , d^3_Y) \middle)\right|_{F(X^i, Y^j)}\\
=& d^{N-3} \left( \binom{3}{3}_{\xi}F(d_X^3 , \id_Y)  +\xi^{i}\binom{3}{2}_{\xi}F(d_X^2 , d_Y)  + \xi^{2i} \binom{3}{1}_{\xi}F(d_X ,d_Y^2)+\xi^{3i}\binom{3}{0}_{\xi}F(d_X , d^3_Y )\middle)\right|_{F(X^i, Y^j)}\\
&\vdots \\
=&\left. \sum_{k=0}^N \xi^{\frac{k(k-1)}{2}} \binom{N}{N-k}_{\xi} F(d_X^{N-k}, d_Y^{k})\right|_{F(X^i, Y^j)}\\
=& F(d_X^{N} , \id_Y) + F(\id_X , d_Y^N)\\
=&0.
\intertext{Hence} 
d_{F_\xi(X, Y)}^N=&F(d_X^{N}, \id_Y) +  F(\id_X , d_Y^N)=0
\end{align*}
on all of $X\otimes_\xi Y$ and so it is an $N$-complex as desired.
\item Since this case is a combination of the above cases we will use both of the above calculations. First we will use that for $\xi\in R$ a primitive $k$-th root of unity and $l \in \N$ we have
\begin{align*}
 d_{F_\xi(X, Y)}^k=& F(d_X^{k}, \id_Y) + F(\id_X , d_Y^k).
\end{align*}
Using this we observe that
\begin{align*}
 d_{F_\xi(X, Y)}^N&=d_{F_\xi(X, Y)}^{p^nk}= d_{F_\xi(X, Y)}^{k}\circ    d_{F_\xi(X, Y)}^{k} \circ ...\circ  d_{F_\xi(X, Y)}^{k},
\end{align*}
with $p^n$ factors in the factorization. Now consider $d_{F_\xi(X, Y)}^k$ instead of $d_{F_\xi(X, Y)}$. The morphism $d_{F_\xi(X, Y)}^k$ now behaves the same way as $d_{F_\xi(X, Y)}$, with $\xi$ replaced by $1$ and $N$ replaced by $p^n$, yielding an analogous calculation to case (1):
\begin{align*}
 d_{F_\xi(X, Y)}^N=&d_{F_\xi(X, Y)}^{N-k}(F(d_X^{k} , \id_Y) +  F(\id_X, d_Y^k))\\
 =&d_{F_\xi(X, Y)}^{N-2k}(F(d_X^{2k} , \id_Y) + 2F(d_X^{k} , d_Y^{k})   + F(\id_X, d_Y^{2k}))\\
 =&F_\xi(d_{X, Y}^{N-3k}(F(d_X^{3k} , \id_Y) + 3 F(d_X^{2k} , d_Y^{k}) + 3F(d_X^{k} d_Y^{2k})   + F(\id_X ,d_Y^{3k}))\\
 &\vdots \\
=&\sum_{l=0}^{p^n} \dbinom{p^n}{l} F(d^{(p^n-l)k}_X, d^{lk}_Y)\\
=&F(d_X^{p^nk},\id_Y) + F(\id_X , d_Y^{p^nk})\\
=&F(d_X^N, \id_Y) + F(\id_X , d_Y^N)\\
=&0,
\end{align*}
so $F_\xi(X,Y)$ is again an $N$-complex.
\end{enumerate}
\end{proof}
\end{theorem}

\begin{remark}\label{calculation conceptual}
The above proof can be also seen as an application of the formula 
$$\prod_{k=0}^{N-1}(x+q^ky)=\sum_{k=0}^{N}q^{\frac{k(k-1)}{2}}\binom{n}{k}_q x^{N-k}y^k $$ 
from Lemma~\ref{q binom} by setting $x:=F(d_X,\id_Y), y:= F(\id_X, d_Y)$ and $q:= \xi$. Observing that $\xi$ is a $k$-th primitive root of unity with $k \mid N$ yields that every distinct power of $\xi$ appears precisely $\frac{N}{k}$ times and so we get the desired result. However, we proved this $q$-binomial formula in the ring $\Z[x,y,q]$, so we would need a general argument to apply it here, which we do not think is necessary, since we only applied it in full strength in case (2).
\end{remark}

Another important type of functors we need to lift are functors which are covariant in one variable and contravariant in the other since these arise as internal homs. In particular they are of interest if we want to construct a closed monoidal structure on the category of $N$-complexes.

\begin{definition}\label{lifting contravariant/covariant}
Let $R$ be a commutative ring, $\C,\C'$ and $\D$ be categories enriched over $\Mod_R$ such that $\D$ admits countable products, $\xi \in R$ and $F(\_,\_):\C^{op} \otimes_{\Mod_R} \C' \to \D$ an $R$-linear bifunctor. Then define for two $\infty$-complexes $X,Y$ the $\infty$-complex $F_\xi(X,Y)$ as 
\begin{align*}
F_\xi(X,Y)^i&:=\prod_{j\in \Z} F(X^j,Y^{i+j})
\end{align*}
with the differential defined on factors by 
\begin{align*}
 \pi_{j,i+j} \circ d_{F_\xi(X,Y)}:= F(X^j,d_Y)\circ \pi_{j,i+j-1} - \xi^{i-1} F(d_X,Y^{i+j})\circ \pi_{j+1,i+j}.
\end{align*}
Here $\pi_{k,l}$ denotes the canonical morphism $$\pi_{k,l}: \prod_{j\in \Z} F(X^j,Y^{i+j})\to F(X^k,Y^{l})$$ for $i=l-k$.

Furthermore define for a morphism of $\infty$-complexes $f:X \to X'$ the family of morphisms in $\C$:
\begin{align*}
F_\xi(f,Y)&:=\left( \prod_{j\in \Z}F(f^j,Y^{i+j}):F_\xi(X',Y)^i \to F_\xi(X,Y)^i\right)_{i \in \Z},
\end{align*}
and for $f:Y \to Y'$ the family of morphisms in $\C'$:
\begin{align*}
F_\xi(X,f)^i&:=\left( \prod_{j\in \Z}F(X^j,f^{i+j}): F_\xi(X,Y)^i \to F_\xi(X,Y')^i \right)_{i \in \Z} .
\end{align*}
\end{definition}

\begin{proposition}\label{functoriality contravariant/covariant}
Let $R$ be a commutative ring, $\C,\C'$ and $\D$ be additive categories enriched over $\Mod_R$ such that $\D$ admits countable products, $\xi \in R$ and $F(\_,\_):\C^{op} \otimes_{\Mod_R} \C' \to \D$ an $R$-linear bifunctor. Then 
the assignments from Definition~\ref{lifting contravariant/covariant} define a functor 
\begin{align*}
F_\xi(\_,\_):\CC_\infty(\C)^{op} \times \CC_\infty(\C') \to \CC_\infty(\D).
\end{align*}
\begin{proof}
One just needs to check that for $f:X \to X'$ in $\C$ and $Y \in \C'$, respectively $Y\in \C$ and $f:X \to X'$ in $\C'$, the families $F_\xi(f,Y)$ and $F_\xi(Y,f)$ indeed define morphisms of $\infty$-complexes since the assignments are levelwise products of functors. Consider first $$F_\xi(f,Y): F_\xi(X',Y) \to F_\xi(X,Y)$$ for $f:X\to X' \in \C(X,X')$ and $Y \in \C'$. So we need to show that 
\begin{equation*}
\tikz[heighttwo,xscale=4,yscale=2,baseline]{
\node (1) at (0,1){$F_\xi(X',Y)^{i-1}$};
\node (2) at (1,1){$F_\xi(X,Y)^{i-1}$};
\node (3) at (0,0){$F_\xi(X',Y)^{i}$};
\node (4) at (1,0){$F_\xi(X,Y)^{i}$};
\draw[->]
(1) edge node[left] {$d_{F_\xi(X',Y)}$}(3)
(2) edge node[right] {$d_{F_\xi(X,Y)}$}(4)
(1) edge node[above] {$F_\xi(f,Y)^{i-1}$} (2)
(3) edge node[below] {$F_\xi(f,Y)^{i}$} (4)
;
}
\end{equation*} 
commutes. This we check on factors:
\begin{align*}
&\pi_{j,i+j}\circ \left( F_\xi(f,Y)^{i} \circ d_{F_\xi(X',Y)}\right)  \\
&=  \pi_{j,i+j}\circ \prod_{j'\in \Z} F(f^{j'},Y^{i+j'}) \circ d_{F_\xi(X',Y)}\\
&=F(f^j,Y) \circ \pi_{j,i+j} \circ  d_{F_\xi(X',Y)}\\
&= F(f^j,Y) \circ (F({X'}^j,d_Y)\circ \pi_{j,i+j-1} - \xi^i F(d_{X'},Y^{i+j})\circ \pi_{j+1,i+j})  \\
&= F(f^{j},d_Y)\circ \pi_{j,i+j-1} -\xi^i F(f^{j} \circ d_{X'} , Y^{j+i})\circ \pi_{j+1,i+j}  \\
&=  F(f^{j},d_Y)\circ \pi_{j,i+j-1} - \xi^i F((d_{X'} \circ f^{j})^{op}   , Y^{j+i})\circ \pi_{j+1,i+j}  \\
&= F(f^{j},d_Y)\circ \pi_{j,i+j-1}- \xi^i F((f^{j+1} \circ d_{X} )^{op} , Y^{j+i})\circ \pi_{j+1,i+j}  \\
&= F(f^{j},d_Y)\circ \pi_{j,i+j-1} - \xi^i F(d_{X}  \circ f^{j+1} , Y^{j+i})\circ \pi_{j+1,i+j}   \\
&=  F(X,d_Y)\circ F(f^{j},Y^{i+j-1}) \circ\pi_{j,i+j-1 }- \xi^i F(d_{X} , Y^{j+i})\circ F(f^{j+1} , Y^{j+i})\circ \pi_{j+1,i+j} \\
&= \left(F(X,d_Y) \circ \pi_{j,i+j-1}- \xi^i F(d_{X} , Y^{j+i})\circ  \pi_{j+1,i+j}  \right) \circ \prod_{j'\in\Z}F(f^{j'},Y^{i+j'-1}) \\
&= \pi_{j,i+j} \circ (d_{F_\xi(X,Y)} \circ F_\xi(f,Y)^{i-1}).
\end{align*}
Now we show for $F_\xi(Y,f): F_\xi(Y,X) \to F_\xi(Y,X')$ for $Y \in \C$ and $f:X \to X'$ a morphism in $\C'$ that the diagram
\begin{equation*}
\tikz[heighttwo,xscale=4,yscale=2,baseline]{
\node (1) at (0,1){$F_\xi(Y,X)^{i-1}$};
\node (2) at (1,1){$F_\xi(Y,X')^{i-1}$};
\node (3) at (0,0){$F_\xi(Y',X)^{i}$};
\node (4) at (1,0){$F_\xi(Y,X')^{i}$};
\draw[->]
(1) edge node[left] {$d_{F_\xi(Y,X)}$}(3)
(2) edge node[right] {$d_{F_\xi(Y,X')}$}(4)
(1) edge node[above] {$F_\xi(Y,f)^{i-1}$} (2)
(3) edge node[below] {$F_\xi(Y,f)^{i}$} (4)
;
}
\end{equation*} 
indeed commutes. Again we check on factors:
\begin{align*}
&\pi_{j,i+j}\circ \left( F_\xi(Y,f)^{i} \circ d_{F_\xi(Y,X)}\right) \\
 &= \pi_{j,i+j}\circ \left(\prod_{j'\in \Z}F(Y^{j'},f^{i+j'})\right) \circ d_{F_\xi(Y,X)} \\
&=F(Y^j,f^{i+j})\circ \pi_{j,i+j}\circ d_{F_\xi(Y,X)}\\
&=F(Y^j,f^{i+j})\circ  \left( F(Y^j,d_X)\circ \pi_{j,i+j-1}- \xi^i F(d_Y,X^{i+j})\circ \pi_{j+1,i+j} \right)\\
&=  F(Y^j,f^{i+j})\circ F(Y^j,d_X)\circ \pi_{j,i+j-1} - \xi^i F(Y^j,f^{i+j})\circ F(d_Y,X^{i+j})\circ \pi_{j+1,i+j} \\
&= F(Y^j,f^{i+j}\circ d_X)\circ \pi_{j,i+j-1}- \xi^i F(d_Y,f^{i+j})\circ \pi_{j+1,i+j}\\
&=   F(Y^j,d_X \circ  f^{i+j-1})\circ \pi_{j,i+j-1}- \xi^i F(d_Y,X^{i+j}) \circ  F(Y^{j+1},f^{i+j})\circ \pi_{j+1,i+j} \\
&=  F(Y^j,d_X) \circ  F(Y^j,f^{i+j-1})\circ \pi_{j,i+j-1}- \xi^i F(d_Y,X^{i+j}) \circ  F(Y^{j+1},f^{i+j})\circ \pi_{j+1,i+j}  \\
&= \left( F(Y^j,d_X) \circ \pi_{j,i+j-1}- \xi^i F(d_Y,X^{i+j}) \circ  \pi_{j+1,i+j}\right) \circ \left(\prod_{j'\in \Z}F(Y^{j'},f^{i+j'-1})\right)\\
&=\pi_{j,i+j} \circ d_{F_\xi(X,Y)} \circ F_\xi(Y,f)^{i-1}.
\end{align*}
Here the shift in $Y$ between the fourth and sixth line occurs since $F(\_,X)$ is contravariant.
\end{proof}
\end{proposition}

Similarly to the case of additive bifunctors, we want that $F_\xi(\_,\_)$ again defines a functor on the category of $N$-complexes for certain $\xi$ and $F(\_,\_)$ a contravariant/covariant bifunctor.

\begin{theorem}\label{induced contravariantcovariant}
Let $R$ be an integral domain, $\C,\C'$ and $\D$ be additive categories enriched over $\Mod_R$, $\xi \in R$ such that $\D$ admits countable products, $N \in \N$ such that $(\D,N,\xi)$ allows the twisted Leibniz rule and let $F(\_,\_):\C^{op} \otimes_{\Mod_R} \C' \to \D$ be an $R$-linear bifunctor. Then the functor  
\begin{align*}
F_\xi(\_,\_):\CC_\infty(\C)^{op} \times \CC_\infty(\C') &\to \CC_\infty(\D)
\intertext{induces a functor}
F_\xi(\_,\_):\CC_N(\C)^{op} \times \CC_N(\C') &\to \CC_N(\D).
\end{align*}
\begin{proof}
By Proposition~\ref{functoriality contravariant/covariant} we have that $F_\xi(\_,\_)$ induces a functor $\CC_\infty(\C)^{op} \times \CC_\infty(\C') \to \CC_\infty(\D)$, so it suffices to check in all cases, that for $X\in \CC_N(\C)$ and $Y \in \C'$ the complex $F_\xi(X,Y)$ is again an $N$-complex. This is analogous to the proof of Theorem~\ref{3.12}, still we give for completeness reasons the proof. Similarly to Theorem~\ref{3.12} we will check all three cases in Definition~\ref{dagger} separately. Furthermore we omit for legibility the projections in the calculations.
\begin{enumerate}[label=(\roman*)]
\item Since every morphism in $\D$ is a $p$-torsion for a prime number $p$ and $N=p^n$, we may use that for all morphisms $f$ in $\D$ we have $p f=0 \in \D$, hence it suffices to consider $\D$ as enriched over $R/pR$. Now we may use in $R/pR$ that $\binom{N}{k}=0$ for all $0<k<N$ just like in Theorem~\ref{lifting bifunctor Ncomplex}.
Using this one gets
\begin{align*}
\pi_{j,i+j} \circ  d_{F_1(X,Y)}^N =& (F(d_X,Y^{i+j}) -  F(X^j,d_Y))\circ d_{F_\xi (X,Y)}^{N-1}\\
=&(F(d^2_X,Y^{i+j}) - F(d_X,d_Y) - F(d_X,d_Y) + F(X^j,d^2_Y))\circ d_{F_\xi (X,Y)}^{N-2}\\
=& (F(d^2_X,Y^{i+j}) - 2 F(d_X,d_Y) + F(X^j,d^2_Y))\circ d_{F_\xi (X,Y)}^{N-2}\\
=&(F(d^3_X,Y^{i+j}) - 3F(d^2_X,d_Y) + 3 F(d_X,d^2_Y) -F(X^j,d^3_Y))\circ d_{F_\xi (X,Y)}^{N-3}\\
&\vdots\\
=&\sum_{k=1}^{N} (-1)^k \dbinom{N}{k} F(d^{N-k}_X,d^k_Y)\\
=&F(d^N_X,Y^{i+j}) + (-1)^N F(X^j,d^N_Y))\\
=&0
\end{align*}
and so $F_\xi(X,Y)$ is an $N$-complex.
\item
In this case we will use the following two formulas from Proposition~\ref{q binomial formula} and Corollary~\ref{q binomial vanishing at root of unity}:
Where we need for $\binom{N}{k}_\xi=0$ that $\xi$ is an $N$-th root of unity. 
Now we calculate analogously to the formula of Lemma~\ref{q binom} $d^N_{F_\xi(X,Y)}$ for each summand, keeping track of the respective degrees and setting $x=F(\id_X,d_Y),\; y=-F(d_X,\id_Y)$ and $q=\xi$:
\begin{align*} &\pi_{j,i+j} \circ  d_{F_\xi (X,Y)}^N \\
=&\left(  F(X^j,d_Y) - \xi^i F(d_X,Y^{i+j})\right)\circ d_{F_\xi (X,Y)}^{N-1}\\
=&\left(F(X^j,d^2_Y) - \xi^i F(d_X,d_Y) -\xi^{i+1} F(d_X,d_Y)+ \xi^{2i} F(d^2_X,Y^{i+j})    \right)\circ d_{F_\xi (X,Y)}^{N-2}\\
=&\left( \binom{2}{0}_{\xi}F(X^j,d^2_Y)  - \xi^{i} \binom{2}{1}_{\xi} F(d_X,d_Y) +\xi^{2i} \binom{2}{2}_{\xi}F(d^2_X,Y^{i+j}) \right)\circ d_{F_\xi (X,Y)}^{N-2}\\
=& \left( \binom{2}{0}_{\xi}F(X^j,d^3_Y)   -\xi^{i+2}\binom{2}{0}_{\xi}F(d_X,d^2_Y)  -\xi^{i} \binom{2}{1}_{\xi}F(d_X,d^2_Y) \right.\\
& + \left.\xi^{2i+1} \binom{2}{1}_{\xi}F(d^2_X,d_Y)+\xi^{2i}\binom{2}{2}_{\xi}F(d^2_X,d_Y)- \xi^{3i} \binom{2}{2}_{\xi}F(d^3_X,Y^{i+j}) \right)\circ d_{F_\xi (X,Y)}^{N-3}
\\
=&\left(\binom{3}{0}_{\xi}F(X^j,d^3_Y)-  \xi^{i} \binom{3}{1}_{\xi}F(d_X,d^2_Y)+\xi^{2i}\binom{3}{2}_{\xi}F(d^2_X,d_Y)- \xi^{3i} \binom{3}{3}_{\xi}F(d^3_X,Y^{i+j}) \right)\circ d_{F_\xi (X,Y)}^{N-3}\\
&\vdots \\
=& \sum_{k=0}^N \xi^{\frac{k(k-1)}{2}} \binom{N}{N-k}_{\xi}(-1)^k F(d^k_X,d^{N-k}_Y)\\
=& F(d^N_X,Y^{i+j}) + (-1)^k F(X^j,d^N_Y)\\
=&0.
\end{align*}
Hence 
\begin{equation*}
d_{F_\xi (X,Y)}^N=F(d^N_X,Y^{i+j}) + F(X^j,d^N_Y)=0
\end{equation*}
on all factors of $F_\xi (X,Y)$ and so it is an $N$-complex as desired.
\item Since this case is a combination of the above two we will use both of the above calculations. First we will use that for $\xi\in R$ a primitive $k$-th root of unity and $l \in \N$ we have
\begin{align*}
 d_{F_\xi(X, Y)}^k=& F(d_X^{k}, \id_Y) + (-1)^k F(\id_X , d_Y^k).
\end{align*}
Using this we observe that the following holds:
\begin{align*}
 d_{X\otimes_\xi Y}^N&=d_{X\otimes_\xi Y}^{p^nk}= d_{X\otimes_\xi Y}^{k}\circ    d_{X\otimes_\xi Y}^{k} \circ ...\circ  d_{X\otimes_\xi Y}^{k},
\end{align*}
with $p^n$ factors in the composition. Now consider $d_{F_\xi(X, Y)}^k$ instead of $d_{F_\xi(X, Y)}$. The morphism $d_{F_\xi(X, Y)}^k$ now behaves the same way as $d_{F_\xi(X, Y)}$, with $\xi$ replaced by $1$ and $N$ replaced by $p^n$, yielding an analogous calculation to case (1):
\begin{align*}
 \pi_{j,i+j} \circ d_{F_\xi (X,Y)}^N=&\left(F(d^k_X,Y^{i+j}) +  (-1)^k F(X^j,d^k_Y)\right) \circ d_{F_\xi (X,Y)}^{N-k}\\
 =&\left(F(d^{2k}_X,Y^{i+j}) + (-1)^k 2F(d^k_X,d^k_Y)   + (-1)^{2k} F(X^j,d^{2k}_Y)\right) \circ d_{F_\xi (X,Y)}^{N-2k}\\
 =&\left((F(d^{3k}_X,Y^{i+j}) + (-1)^k 3F(d^{2k}_X,d^k_Y) + (-1)^{2k} 3F(d^k_X,d^{2k}_Y)   + (-1)^{3k} F(X^j,d^{3k}_Y)\right) \circ d_{F_\xi (X,Y)}^{N-3k}\\
 &\vdots\\
=&\sum_{l=1}^{p^n} \dbinom{p^n}{l}(-1)^{lk}F(d^{(p^n-l)k}_X,d^{lk}_Y)\\
=&F(d^{p^nk}_X,Y^{i+j}) + (-1)^{nk} F(X^j,d^{p^nk}_Y)\\
=&F(d^{N}_X,Y^{i+j}) + (-1)^{nk} F(X^j,d^{N}_Y)\\
=&0.
\end{align*}
\end{enumerate}
So $F_\xi (X,Y)$ is again an $N$-complex.
\end{proof}
\end{theorem}
 Having established these tools for the lifting of functors, we can finally construct the monoidal structure on the category of $N$-complexes and show that it is closed if the underlying category is the category of modules over an integral domain.

\section{Monoidal structure}\label{section 4}

In this section we will generalize the tensor product on the category of chain complexes to arbitrary $N$-complexes. Since we want to do this in a very general context we need to introduce first the notions of countable distributive monoidal categories and monoidal categories enriched over a ring $R$. This will allow us to use liftings of bifunctors from the previous section to define a monoidal structure.

\subsection{Background on distributive monoidal structures}
First we will discuss general distributivity properties for monoidal categories in order to use the liftings from the previous section.

\begin{definition}
Let $(\C,\otimes,\unit,\alpha,\lambda,\rho)$ be a monoidal category that admits finite coproducts. Then define the following morphisms:
\begin{enumerate}

\item For $X,Y,Z\in \C$ the \textbf{left distributivity morphism} 
\begin{align*}
{}_{X,Y,Z}\delta : (X \otimes Y) \oplus (X \otimes Z) \to X \otimes (Y \oplus Z),
\end{align*}
 is the morphism induced by the universal property of the coproduct in the following diagram:
\begin{equation*}
\tikz[heighttwo,xscale=2,yscale=2,baseline]{
\node (1) at (0,1){$ X\otimes Y$};
\node (2) at (2,2){$ X\otimes Z$};
\node (3) at (2,1){$(X \otimes Y)\oplus (X \otimes Z)$};
\node (4) at (4,0){$ X \otimes (Y \oplus Z)$};
\node[overlay] (Point) at (4.7,-0.05){.};
\draw[->]
(1) edge node[above]{$  \iota_{X\otimes Y}$}(3)
(2) edge node[left]{$\iota_{X \otimes Z} $}(3)
(1) edge node[below left]{$ \id_X \otimes \iota_Y$}(4)
(2) edge node[above right]{$\id_X \otimes \iota_Z$}(4)
;
\draw[dashed,->]
(3) edge (4);
}
\end{equation*}

\item For $X,Y,Z\in \C$ the \textbf{right distributivity morphism} 
\begin{align*}
\delta_{Y,Z,X} : (Y \otimes X) \oplus (Z \otimes X) \to  (Y \oplus Z) \otimes X,
\end{align*}
is the morphism induced by the universal property of the coproduct and the morphisms 
\begin{align*}
 \iota_Y \otimes \id_X &: (Y \otimes X)\to(Y \oplus Z) \otimes X \\
 \iota_Z \otimes \id_X &: (Z \otimes X)\to(Y \oplus Z) \otimes X .
\end{align*} 
\end{enumerate}
\end{definition}

\begin{remark}
By construction the left distributivity morphism defines a natural transformations between the two functors
\begin{align*}
\C \times \C \times \C &\to \C \\
(X,Y,Z) &\mapsto (X\otimes Y)\oplus (X \otimes Z) \\ 
(X,Y,Z) &\mapsto X\otimes (Y\oplus Z).
\end{align*}
The same holds for the right distributivity morphism and 
\begin{align*}
\C \times \C \times \C &\to \C \\
(X,Y,Z) &\mapsto (Y\otimes X)\oplus (Z \otimes X) \\ 
(X,Y,Z) &\mapsto  (Y\oplus Z)\otimes X .
\end{align*}
\end{remark}

\begin{definition}\label{distributive}
A monoidal category $(\C,\otimes,\unit,\alpha,\lambda,\rho)$, is called
\begin{enumerate} 
\item \textbf{distributive} if it is additive, the morphisms ${}_{X,Y,Z}\delta$ and $\delta_{Y,Z,X}$ are isomorphisms for all $X,Y,Z \in \C$ and $X \otimes 0 \cong 0 \cong 0 \otimes X$ for all $X \in \C$.
\item \textbf{countable distributive} if $\C$ is distributive, admits countable coproducts and for all $X,Y_i \in \C$, for $i \in \Z$ the canonical morphisms
 \begin{align*}
 \tau_{X,Y_i} &: \bigoplus_{i \in \Z} (X \otimes Y_i)\to  X \otimes \bigoplus_{i \in \Z} Y_i \\ 
 \sigma_{Y_i,X} &: \bigoplus_{i \in \Z} ( Y_i \otimes X)\to  \left(\bigoplus_{i \in \Z} Y_i \right) \otimes X
 \end{align*} are isomorphisms.
\end{enumerate}
\end{definition}

\begin{definition}\label{3.5}
Let $(\C,\otimes_\C,\unit_\C,\alpha_\C,\lambda_\C,\rho_\C)$ be a monoidal category enriched over a symmetric monoidal category $(\D,\otimes_\D,\unit_\D,\alpha_\D,\lambda_\D,\rho_\D)$. 
We call $\C$ a \textbf{$\D$-monoidal} category if the bifunctor $\otimes_\C : \C \times \C  \to \C$ and structure morphisms factor through the canonical functor $\C \times \C \to\C \otimes_\D \C$, which is the identity on objects and $\otimes_\D$ on morphism spaces. 
\end{definition}

\begin{remark}
Definition~\ref{3.5} implies that $\_ \otimes_\C Y$ , $X\otimes_\C \_$ and the natural transformations involved in the monoidal structure of $\C$ are in fact $\D$-enriched functors and natural transformations of $\D$-enriched functors.
\end{remark}

Now we have a good setup for applying the constructions from Section 3 and we can define the monoidal structure on $\CC_N(\C)$. However it is not yet clear that this indeed defines a monoidal structure since we are in a very general setup at the moment.

\begin{definition}\label{3.7}
Let $(\C,\otimes, \unit,\alpha,\lambda,\rho)$ be a countable distributive $\Mod_R$-monoidal category, for $R$ a commutative ring, and $\xi \in R$. Then define the following  on $\CC_\infty(\C)$: 
\begin{enumerate}
\item A functor $\_ \otimes_\xi \_ : \CC_\infty(\C) \times \CC_\infty(\C) \to \CC_\infty(\C)$ just as in Definition~\ref{lifting bifunctor} using the bifunctor $\_ \otimes \_ :\C \otimes_\D \C \to \C$.
\item The $\infty$-complex $\unit_\xi:= \mu^0_1(\unit)$.
\item The morphism $\alpha_{\xi,X,Y,Z} : (X\otimes Y) \otimes Z \to X\otimes (Y \otimes Z)$ levelwise defined as (cf. Remark~\ref{3.8}) 
$$(\alpha_{\xi,X,Y,Z})^j:=\left(\bigoplus_{l+m=n}\sigma_{X^i\otimes Y^k,Z^l}\right)\circ \gamma_{l,i,k}^{-1} \circ \left(\bigoplus_{i+k+l=n}\alpha_{X^i,Y^k,Z^l}\right) \circ \gamma_{i,k,l} \circ \left( \bigoplus_{i+j=n}\tau^{-1}_{X^i,Y^k \otimes Z^l}\right).$$
 Here $\gamma_{i,k,l}$ is the canonical natural isomorphism $$\bigoplus_{i+j=n} \left(\bigoplus_{k+l=j}M_{i,k,l} \right) \to \bigoplus_{i+k+l=n}M_{i,k,l} .$$
\item The morphism $\lambda_{\xi,X} :\unit_\xi \otimes_\xi X  \to X$ levelwise by $(\lambda_\xi(X))^j:=\lambda(X^j)$.
\item And the morphism $\rho_{\xi,X} : X\otimes_\xi \unit_\xi \to X$  levelwise by $(\rho_\xi(X))^j:=\rho(X^j)$.
\end{enumerate}
\end{definition}

\begin{remark}\label{3.8} Regarding Definition~\ref{3.7} the following remarks are to be made:
\begin{enumerate}
\item The construction of $\alpha_{\xi,X,Y,Z}$ can be seen more intuitively in the following diagram in $\C$, which also explains the numerous assumptions on the monoidal structure of $\C$:
\begin{equation*}
\tikz[heighttwo,xscale=2,yscale=3,baseline]{
\node (1) at (0,2){$ \bigoplus\limits_{i+j=n}\left( X^i \otimes \left(\bigoplus\limits_{k+l=j}Y^k\otimes Z^l\right)\right)$};
\node (2) at (0,1){$ \bigoplus\limits_{i+j=n}\left(\bigoplus\limits_{k+l=j}X^i \otimes \left( Y^k\otimes Z^l \right)\right)$};
\node (3) at (0,0){$ \bigoplus\limits_{i+k+l=n}X^i \otimes \left( Y^k\otimes Z^l \right)$};
\node (4) at (4,0){$  \bigoplus\limits_{i+k+l=n} \left(X^i \otimes Y^k \right) \otimes Z^l$};
\node (5) at (4,1){$  \bigoplus\limits_{m+l=n}\left(\bigoplus\limits_{i+k=m}\left(X^i \otimes  Y^k\right)  \otimes Z^l \right) $};
\node (6) at (4,2){$ \bigoplus\limits_{m+l=n}\left( \left(\bigoplus\limits_{i+k=m}X^i \otimes Y^k\right)\otimes Z^l\right)$};

\node[overlay] (Point) at (5.2,-0.05){.};

\draw[->]
(1) edge node[left]{$  \bigoplus\limits_{i+j=n}\tau^{-1}_{X^i,Y^k \otimes Z^l}$}(2)
(2) edge node[left]{$\gamma_{i,k,l} $}(3)
(3) edge node[below]{$ \bigoplus\limits_{i+k+l=n}\alpha_{X^i,Y^k,Z^l}$}(4)
(4) edge node[right]{$\gamma_{l,i,k}^{-1} $}(5)
(5) edge node[right]{$ \bigoplus\limits_{m+l=n}\sigma_{X^i\otimes Y^k,Z^l}$}(6)
;
\draw[dashed,->]
(1) edge (6);
}
\end{equation*} 
\item By Definition~\ref{distributive} the morphisms $\lambda_{\xi,X}$ and $\rho_{\xi,X}$ are well defined.
\item The restriction on $\C$ being enriched over some ring can be easily dropped, since if $\C$ is just an additive category, it is already enriched over $\Z$. In this case however, one restricts oneself to $\xi \in \Z$.
\end{enumerate}
\end{remark}

Using these remarks we can prove that we obtain a monoidal structure on the category of $\infty$-complexes. This will already give us a monoidal structure on the category of $N$-complexes if $(\C,N,\xi)$ allows the twisted Leibniz rule.

\begin{proposition}\label{3.9}
Let $(\C,\otimes, \unit,\alpha,\lambda,\rho)$ be a countable distributive $\Mod_R$-monoidal category for $R$ a commutative ring and $\xi \in R$. Then the morphisms from Definition~\ref{3.7} define natural isomorphisms
\begin{align*}
\alpha_\xi :(\_ \otimes_\xi \_ ) \otimes_\xi \_ &\xrightarrow[]{\sim} \_ \otimes_\xi (\_ \otimes_\xi \_ )\\
\lambda_\xi :(\unit_\xi \otimes_\xi \_ ) &\xrightarrow[]{\sim} \id_\C \\
\rho_\xi :  (\_ \otimes_\xi \unit_\xi) &\xrightarrow[]{\sim} \id_\C .
\end{align*}
\begin{proof}
Every morphism involved in Definition~\ref{3.7} is a natural isomorphism, so all of the above morphisms define levelwise natural isomorphisms. Since $\C$ is a $\Mod_R$-monoidal category, these natural isomorphisms commute with $\xi$ and hence also with the respective differentials, so they define natural isomorphisms.
\end{proof}
\end{proposition}

The proof of the following theorem is a bit sketchy. However, it is really straight forward, but since we have to deal with countable coproducts combined with the pentagon identity and the triangle identity, a rigid proof would take some pages of diagrams which do not give additional insight.

\begin{theorem}\label{3.10}
Let $(\C,\otimes, \unit,\alpha,\lambda,\rho)$ be a countable distributive $\Mod_R$-monoidal category for $R$ a commutative ring and $\xi \in R$. Then $(\CC_\infty(\C),\otimes_\xi,\unit_\xi ,\alpha_\xi, \lambda_\xi,\rho_\xi)$ defines a monoidal structure on $\CC_\infty(\C)$ which is countable distributive.
\begin{proof}
All of the above structure morphisms are natural isomorphisms between the corresponding functors by Proposition~\ref{3.9}. Hence the only thing left to check are the triangle and the pentagon identities. By Proposition~\ref{3.9} and Remark~\ref{3.8} we have natural isomorphisms 
\begin{align*}
\bigoplus\limits_{i+j=n}\left( X^i \otimes \left(\bigoplus\limits_{k+l=j}Y^k\otimes Z^l\right)\right) &\cong \bigoplus\limits_{\mathclap{i+k+l=n}}X^i \otimes \left( Y^k\otimes Z^l \right)\\ \bigoplus\limits_{\mathclap{i+k+l=n}} \left(X^i \otimes Y^k \right) \otimes Z^l &\cong  \bigoplus\limits_{m+l=n}\left( \left(\bigoplus\limits_{i+k=m}X^i \otimes Y^k\right)\otimes Z^l\right).
\end{align*}
So it suffices if the pentagon and triangle identities are satisfied on summands of the form $X^i\otimes \left( Y^k\otimes Z^l \right)$. There they are clearly satisfied, since $(\C,\otimes, \unit,\alpha,\lambda,\rho)$ is a monoidal category and $$\unit_\xi^j=\begin{cases} \unit \text{ if } j=0, \\ 0 \text{ else.}\end{cases}$$

Now for countable distributivity observe first of all that by Proposition~\ref{complete/cocomplete} the category $\CC_\infty(\C)$ admits countable coproducts and is additive by Corollary~\ref{CN additive}. Furthermore, since $\_\otimes_\xi \_$ is defined levelwise as a coproduct in a countable distributive monoidal category, we get that $\_ \otimes_\xi \_$ is levelwise countable distributive and hence countable distributive since coproducts are computed levelwise.
\end{proof}
\end{theorem}

\begin{corollary}\label{3.12}
Let $(\C,\otimes, \unit,\alpha,\lambda,\rho)$ be a countable distributive $\Mod_R$-monoidal category, for $R$ an integral domain, $\xi \in R$ and $N\in \N$. Then $(\CC_N(\C),\otimes_\xi,\unit_\xi, \alpha_\xi, \lambda_\xi,\rho_\xi)$ is a countable distributive monoidal category if the triple $(\C,N,\xi)$ allows the twisted Leibniz rule (cf. Definition~\ref{dagger}).
\begin{proof}
Since by Proposition~\ref{3.10} $(\CC_\infty(\C),\otimes_\xi,\unit_\xi , \alpha_\xi, \lambda_\xi,\rho_\xi)$ is a countable distributive monoidal category it suffices to check in all cases, that for $X,Y \in \CC_N(\C)$ the complex $X\otimes_\xi Y$ is again an $N$-complex. But this is just the statement of Proposition~\ref{lifting bifunctor Ncomplex}, since $\_\otimes \_ $ is in all cases a bifunctor and $(\C,N,\xi)$ allows the twisted Leibniz rule.
\end{proof}
\end{corollary}

\subsection{Closed monoidal structure}
$\;$

Since most interesting categories are not just monoidal but actually closed monoidal, we also want to lift the adjunctions needed for the structure of a closed monoidal category and an analogous structure on $\CC_N(\C)$ for $N \in \N \cup \{\infty\}$.
For this we will use again the liftings we constructed in Section 3 . However, before we can define the internal hom we need to make the following general observation.

\begin{proposition}\label{adjoint contravariant covariant}
Let $F:\C \times \C' \to \D$ be a bifunctor, such that for all $X \in \C'$ the functor $F(\_,X)$ admits a right adjoint $G(X,\_): \D \to \C$. Then $G(X,\_)$ induces a bifunctor $$G(\_,\_):\C'^{op} \times \D \to \C.$$
\begin{proof}
The functoriality in the right variable is given by assumption, hence it just remains to be shown that a morphism $f:X \to X'$ in $\C'$ induces naturally a morphism $$G(f,Y):G(X',Y)\to G(X,Y).$$ For this consider the following composition:
\begin{align*}
F(G(X',Y),X) \xrightarrow{F(\id,f)} F(G(X',Y),X')  \xrightarrow{\epsilon} Y,
\end{align*}
where $\epsilon$ is the counit of the adjunction. This composition is adjoint to a morphism $$G(f,Y):G(X',Y) \to G(X,Y),$$ which is by construction compatible with the structure of $G(X,\_):\C \to \D$.
\end{proof} 
\end{proposition}

Now using this observation we can lift both functors involved in the structure of a closed monoidal structure, and show that these are still compatible in the desired way if we consider modules over a ring.

We conjecture that the following statement holds in general for a closed countable distributive $\Mod_R$-monoidal category. However the constructions involving $\id_X\otimes d_Y$ and $[d_Y,\id_Z]$ blow up in the general case. Hence we will only prove the case of modules over a ring $R$, since there the proof is less confusing and one has good control of the adjoints.

\begin{theorem}\label{closed monoidal infty}
Let $R$ be  a commutative ring and $\xi \in R$. Then we have an adjunction $$\_ \otimes_\xi X : \CC_\infty(\Mod_R) \leftrightarrow \CC_\infty(\Mod_R) : [X,\_]_\xi.$$ 

Here $$[\_,\_]_\xi:\CC_\infty(\Mod_R)^{op} \times \CC_\infty(\Mod_R) \to \CC_\infty(\Mod_R)$$ is the functor induced by $[\_,\_]:\Mod_R^{op} \otimes_{\Mod_R} \Mod_R \to \Mod_R$ and Definition~\ref{induced contravariantcovariant}.

\begin{proof}
Throughout this proof we denote by $\Gr(\C)$ the category of graded objects in a category $\C$.

First observe that we already have in $\Gr(\Mod_R)$ an isomorphism
$$\Gr(\Mod_R)(X\otimes_\xi Y,Z) \xrightarrow{\sim} \Gr(\Mod_R)(X,[Y,Z]_\xi)$$
by Definition~\ref{lifting contravariant/covariant} and the tensor-hom adjunction in $\Gr(\Mod_R)$.
So it suffices to check that this isomorphism preserves and reflects that a morphism commutes with the differentials.
For this consider a morphism $f \in \CC_\infty(\C)(X\otimes_\xi Y,Z)$. In particular $f\in \Gr(\Mod_R)(X\otimes_\xi Y,Z)$ and $f$ has the property $$ f \circ d_{X \otimes Y} -d_Z \circ f =0.  $$
Evaluation of this map at an element $x \otimes y \in X^i \otimes Y^{n-i} \subset X \otimes Y$ yields now 
\begin{align*}
 (f\circ d_{X \otimes Y} -d_Y \circ f) (x\otimes y) &= f(d_X(x)\otimes y + \xi^i x \otimes d_Y(y))- d_Z(f(x\otimes y))\\
 &=f(d_X(x)\otimes y) + \xi^i f(x \otimes d_Y(y))- d_Z(f(x\otimes y)).
\end{align*}
This is adjoint to the map defined levelwise as
\begin{align*}
X^i \to& [Y,Z]_\xi^{i+1}\\
x \mapsto& \begin{pmatrix} Y^n \to Z^{n+i}  \hfill\\
  \; \; \; y \mapsto f(d_X(x)\otimes y) + \xi^i f(x \otimes d_Y(y))- d_Z(f(x\otimes y))
  \end{pmatrix}
\end{align*} 
since the adjunction in $\Gr(\Mod_R)$ is an adjunction of additive functors.
Now consider the adjoint map $\widehat{f}:X \to [Y,Z]_\xi$ of $f$. Levelwise this is the map
\begin{align*}
X^i &\to [Y,Z]_\xi^{i+1}\\
x &\mapsto  \begin{pmatrix}
Y \to Z  \hfill\\
  \; y \mapsto  f(x \otimes y)
\end{pmatrix}.
\end{align*}
Precomposing with $d_X$ yields that the map $\widehat{f}\circ d_X$ is levelwise of the form
\begin{align*}
X^i \to& [Y,Z]_\xi^{i+1}\\
x \mapsto& \begin{pmatrix} Y^n \to Z^{n+i}  \hfill\\
  \; \; \; y \mapsto f(d_X(x) \otimes y)\end{pmatrix},
\end{align*}
and postcomposing with $d_{[Y,Z]}$ gives that the levelwise structure of $d_{[Y,Z]}\circ \widehat{f}$ is
\begin{align*}
X^i \to& [Y,Z]_\xi^{i+1}\\
x \mapsto& \begin{pmatrix} Y^n \to Z^{n+i}  \hfill\\
  \; \; \; y \mapsto  d_Z(f(x\otimes y)) - \xi^i f(x \otimes d_Y(y))\end{pmatrix}.
\end{align*}
Now just as above, $\widehat{f}$ is actually a morphism of $\infty$-complexes if and only if 
$$0=\widehat{f}\circ d_X -d_{[Y,Z]} \circ \widehat{f}.$$
After evaluating at an $x \in X^i$, this is just the map
\begin{align*}
X^i \to& [Y,Z]_\xi^{i+1}\\
x \mapsto& \begin{pmatrix} Y^n \to Z^{n+i} \hfill\\
  \; \; \;y \mapsto f(d_X(x)\otimes y) + \xi^i f(x \otimes d_Y(y))- d_Z(f(x\otimes y))\end{pmatrix}.
\end{align*}
Hence the adjoint map of $d_Z \circ f - f \circ d_{X \otimes_\xi Y}$ and $\widehat{f}\circ d_Z - d_{[Y,Z]_\xi} \circ \widehat{f}$ have to coincide. 
So we may conclude, using that the adjunction is an adjunction between pointed functors, that $f$ is a morphism of $\infty$-complexes  if and only if $\widehat{f}$ is a morphism of $\infty$-complexes.
\end{proof}

\end{theorem}

\begin{corollary}
Let $N \in \N$, $R$ be an integral domain and $\xi \in R$, then we have an adjunction $$\_ \otimes_\xi X : \CC_N(\Mod_R) \leftrightarrow \CC_N(\Mod_R) : [X,\_]_\xi$$ if $(\Mod_R,N,\xi)$ allows the twisted Leibniz rule.

In particular we get that $(\CC_N(\Mod_R),\otimes_\xi,\unit_\xi, \alpha_\xi, \lambda_\xi,\rho_\xi)$ is a closed monoidal category if $(\Mod_R,N,\xi)$ allows the twisted Leibniz rule.
\begin{proof}
By Theorem~\ref{closed monoidal infty} it just remains to show that $\_\otimes_\xi\_$ and $[\_,\_]_\xi$ induce functors of $N$-complexes in $\Mod_R$, which is just Corollary~\ref{3.12} and Theorem~\ref{induced contravariantcovariant}.
\end{proof}
\end{corollary}

%% file: Frobenius_exact_structure.tex
\section{Frobenius exact structure on $\CC_N(\A)$}\label{section 5}

Similar to the the chain complex case one can construct a Frobenius exact structure on the category of $N$-complexes in an additive category $\A$. We will use this to define the homotopy category of $N$-complexes over $\A$ as the stable category with respect to this exact structure. This yields immediately that the homotopy category of $N$-complexes is triangulated by Happel's Theorem. In particular the work done here repeats arguments by Iyama, Kato and Miyashi from their paper \cite{N-complex} thoroughly. For the notions of exact and Frobenius exact  categories we refer to \cite{Buhler20101} and Appendix A in the paper by B. Keller \cite{keller}.

\begin{definition}
Let $\A$ be an additive category and $N  \in \N$. Then define $\Se_\oplus$ to be the collection of all levelwise split short exact sequences in $\CC_N(\A)$, i.\,e.\ sequences which are levelwise isomorphic to a sequence of the form
\begin{equation*}
\tikz[heighttwo,xscale=2,yscale=2,baseline]{

\node (01) at (0.3,0.3){$0$};
\node (X1) at (1,0.3){$X$};
\node (X2) at (2,0.3){$X\oplus Y$};
\node (X3) at (3,0.3){$Y$};
\node (02) at (3.7,0.3){$0$};
    
\draw[->]
(01) edge (X1)
(X3) edge (02)
(X1) edge node[above] {$\begin{bmatrix}
\id \\ 0
\end{bmatrix}$} (X2)
(X2) edge node[above] {$\begin{bmatrix}
0 & \id 
\end{bmatrix}$} (X3);

}.
\end{equation*}
\end{definition}

We will show that $S_\oplus$ defines a Frobenius exact structure, hence we need to show the following three properties of $\Se_\oplus$:
\begin{enumerate}
\item $\Se_\oplus$ is an exact structure,
\item $\Se_\oplus$ has enough $\Se_\oplus$-projective objects and $\Se_\oplus$-injective objects,
\item the $\Se_\oplus$-projective objects are precisely the $\Se_\oplus$-injective objects.
\end{enumerate}

First we show that it indeed is an exact structure:

\begin{proposition}\label{4.2}
Let $\A$ be an additive category and $N  \in \N$. Then $\Se_\oplus$ defines an exact structure on $\CC_N(\A)$.
\begin{proof}
Observe that $\Se_\oplus$ defines levelwise an exact structure. Now direct sums, pushouts and pullbacks in $\CC_\infty(\A)$ are computed pointwise. Hence $\Se_\oplus$ defines an exact structure on $\CC_\infty(\A)$ and also on $\CC_N(\A)$, since it is a full subcategory.
\end{proof}
\end{proposition}

To see that $\Se_\oplus$ has enough $\Se_\oplus$-projective objects and $\Se_\oplus$-injective objects we will construct for an $N$-complex $X$ an $\Se_\oplus$-injective object $\II_N(X)$ and an $\Se_\oplus$-projective object $\PP_N(X)$ such that we have an admissible monomorphisms $X \to \II_N(X)$ and an admissible epimorphism $\PP_N(X) \to X$. The construction of these is motivated by the case of chain complex, which can be recovered from the given construction by setting $N=2$.

\begin{definition}\label{4.3}
Let $\A$ be an additive category, $N \in \N$ and $X \in \CC_N(\A)$. Then define
\begin{enumerate}
\item the $N$-complex
\begin{equation*}
\II_N(X):=\bigoplus_{i\in\Z}\mu^{i}_N(X^i),
\end{equation*}
\item the $N$-complex
\begin{equation*}
\PP_N(X):=\bigoplus_{i\in\Z}\mu^{i+N-1}_N(X^i),
\end{equation*}
\item
and the morphisms $\lambda_X^{\II}:X \to \II_N(X) $ and $\lambda_X^{\PP}:\PP_N(X)  \to X $ by setting levelwise
\begin{equation*}
(\lambda^\II_X)^i:= \begin{bmatrix} \id_{X^i} & d_X & d_X^2 & \cdots & d_X^{N-2}& d_X^{N-1} \end{bmatrix}^T
\end{equation*} and 
\begin{equation*}(\lambda_X^{\PP})^i:=\begin{bmatrix} d_X^{N-1} & d_X^{N-2} &\cdots & d_X^2 & d_X & \id_{X^{i-N+1}} \end{bmatrix}.
\end{equation*}

\end{enumerate}
\end{definition}

Since the data of an $\Se_\oplus$-injective hull and an $\Se_\oplus$-projective cover contain the morphisms as well, we check that these are morphisms of $N$-complexes.

\begin{proposition}
Let $\A$ be an additive category and $N  \in \N$. Then the maps $\lambda_X^\II : X \to \II_N(X)$ and $\lambda_X^\PP :\PP_N(X) \to X$ define morphisms of $N$-complexes.

\begin{proof} We check separately:\\
\begin{enumerate}
\item[$\lambda_X^\II$:] One needs to check that $d_{\II_N(X)}\circ \lambda_X^\II = \lambda_X^\II \circ d_X$. we calculate 
\begin{align*}
(d_{\II_N(X)} \circ \lambda_X^\II)^i &= d_{\II_N(X)}\circ \begin{bmatrix} \id_{X^i} & d_X & d_X^2 & \cdots & d_X^{N-2}& d_X^{N-1} \end{bmatrix}^T\\
&= \begin{bmatrix} d_X & d_X^2 &d_X^{3}& \cdots & d_X^{N-1} & 0 \end{bmatrix}^T
\end{align*}
and
\begin{align*}
(\lambda_X^\II \circ d_X)^i =& \begin{bmatrix} \id_{X^{i+1}} & d_X & d_X^2 & \cdots & d_X^{N-2}& d_X^{N-1} \end{bmatrix}^T \circ d_X\\
=&  \begin{bmatrix}  d_X & d_X^2 & d_X^3 \cdots & d_X^{N-1}& d_X^N \end{bmatrix}^T.
\end{align*}
So by $d_X^N=0$ the claim holds.
\item[$\lambda_X^\PP$:]
One needs to check again that $d_X\circ \lambda_X^\PP = \lambda_X^\PP \circ d_{\PP(X)}.$ Now calculation yields
\begin{align*}
(d_X \circ \lambda_X^\PP)^i &= d_X\circ \begin{bmatrix} d_X^{N-1} & d_X^{N-2} &\cdots & d_X^2 & d_X & \id_{X^{i-N+1}} \end{bmatrix}\\
 &= \begin{bmatrix} d_X^N & d_X^{N-1} &\cdots & d_X^3 & d_X^2 & d_X  \end{bmatrix}
\end{align*}
and
\begin{align*}
(\lambda_X^\PP \circ d_{\PP_N(X)})^i &=\begin{bmatrix} d_X^{N-1} & d_X^{N-2} &\cdots & d_X^2 & d_X & \id_{X^{i-N+1}} \end{bmatrix} \circ d_{\PP_N(X)} \\
&=\begin{bmatrix} 0 &d_X^{N-1}  &\cdots & d_X^3 &d_X^2 & d_X \end{bmatrix}.
\end{align*}
This again concludes the proof since $d_X^N=0$.\qedhere
\end{enumerate}
\end{proof}
\end{proposition}


Now we also check that $\PP_N(X)$ and $\II_N(X)$ define an $\Se_\oplus$-cover and $\Se_\oplus$-hull respectively.

\begin{proposition} \label{4.4}
There are sequences in $\Se_\oplus$ of the form 
\begin{equation*}
\tikz[heighttwo,xscale=2,yscale=2,baseline]{

\node (01) at (0.3,0.3){$0$};
\node (X1) at (1,0.3){$X$};
\node (X2) at (2,0.3){$\II_N(X)$};
\node (X3) at (3,0.3){$Y$};
\node (02) at (3.7,0.3){$0$};
    
\draw[->]
(01) edge  (X1)
(X3) edge (02)
(X1) edge node[above] {$\lambda^\II_X$}(X2)
(X2) edge  (X3);
},
\end{equation*}
and 
\begin{equation*}
\tikz[heighttwo,xscale=2,yscale=2,baseline]{

\node (01) at (0.3,0.3){$0$};
\node (X1) at (1,0.3){$Y'$};
\node (X2) at (2,0.3){$\PP_N(X)$};
\node (X3) at (3,0.3){$X$};
\node (02) at (3.7,0.3){$0$};
    
\draw[->]
(01) edge  (X1)
(X3) edge (02)
(X1) edge node[above] {$\lambda^\PP_Y$}(X2)
(X2) edge  (X3);
}.
\end{equation*}
\begin{proof}
This becomes obvious if we consider $Y=\bigoplus_{i\in\Z}\mu^{i-1}_{N-1}(X^i)$ in the first sequence and $Y'=\bigoplus_{i\in\Z}\mu^{i+N}_{N-1}(X^i)$ in the second.
\end{proof}
\end{proposition}

To check that $\II_N(X)$ and $\PP_N(X)$ really are $\Se_\oplus$-injective (resp. $\Se_\oplus$-projective) objects, we need to show by the definition of $\Se_\oplus$-injective (resp. $\Se_\oplus$-projective) objects that every sequence in $\Se_\oplus$ such that the rightmost term is $\Se_\oplus$-projective, or the leftmost term is $\Se_\oplus$-injective admits a splitting in the category of $N$-complexes. For this we will explicitly construct the splittings and prove that they indeed are morphisms of $N$-complexes.

\begin{proposition}\label{4.5}
Let $\A$ be an additive category and $N\in \N$. Then we have
\begin{enumerate}
\item $\II_N(X)$ is $\Se_\oplus$-injective,
\item $\PP_N(X)$ is $\Se_\oplus$-projective.
\end{enumerate}
\begin{proof} We will prove in both cases  that if we have a levelwise splitting morphism $\II_N(X) \to Y$ respectively $Y \to \PP_N(X)$ that this morphism already admits a splitting in the $\CC_N(\A)$.
\begin{enumerate}
\item Let $f: \II_N(X) \to Y$ be a morphism of $N$-complexes such that there exists for every $i \in \Z$ a levelwise split $h^i: Y^i \to \II_N(X)^i$, satisfying $h^i \circ f^i =\id_{\II_N(X)}$. We need to provide a morphism of $N$-complexes $g:Y \to \II_N(X)$ such that $g \circ f = \id_{\II_N(X)}$.\\
The $f^i$ and $h^i$ morphisms are by definition of the form
\begin{equation*}
f^i= \begin{bmatrix}f_1^i & f_2^i& \cdots & f_N^i\end{bmatrix} \; , \;h^i=\begin{bmatrix}
h_1^i & h_2^i& \cdots &h_N^i\end{bmatrix}^T,
\end{equation*}
 now the fact that $h^i \circ f^i = \id_{\II_N(X)^i}$ just means $h_j^i\circ f_j^i = \id_{X^{i+j}}$ and $h_j^i\circ f_k^i = 0$ if $k \neq j$. Since $f$ is a morphism of $N$-complexes, we know that $d_Y \circ f^i=f^{i+1} \circ d_{\II_N(X)}$, but this just implies, that $d_Y \circ f_j^i= f_{j-1}^{i+1}$ if $1 < j \le N$ and zero otherwise. Furthermore, observe that by the definition of $\II_N(X)$ the following equation holds for $2\le j \le N $:
\begin{equation*}
 h_j^i\circ f_j^i = \id_{X^{i+j}}=h_{j-1}^{i+1} \circ f_{j-1}^{i+1}.
\end{equation*}
Now define the morphism $g: Y \to \II_N(X)$ levelwise as 
$$g^i:=\begin{bmatrix}h_1^i&h_1^{i+1}\circ d_Y& \cdots &h_1^{i+N-1}\circ d_Y^{N-1}\end{bmatrix}^T.$$
 Note that this $g$ is well defined by the definition of $\II_N(X)$. Furthermore this indeed defines a morphism of $N$-complexes since 
\begin{align*}
(g\circ d_Y)^i& = \begin{bmatrix} h_1^{i+1}\circ d_Y & h_1^{i+2} \circ d_Y^2 &\cdots & h_1^{i+N-1} \circ d_Y^{N-1}& h_1^{i+N} \circ d_Y^N\end{bmatrix}^T \\
&=\begin{bmatrix}h_1^{i+1}\circ d_Y & h_1^{i+2} \circ d_Y^2&\cdots& h_1^{i+N-1} \circ d_Y^{N-1}& 0\end{bmatrix}^T \\
&= (d_{\II_N(X)} \circ g)^i.
\end{align*}
The morphism $g$ is a section in $\CC_N(\A)$, since
\begin{align*}
g \circ f &=
\begin{bmatrix}
h_1^i \\
h_1^{i+1}\circ d_Y\\
h_1^{i+2}\circ d_Y^2\\
\vdots \\
h_1^{i+N-1}\circ d_Y^{N-1}]^T
\end{bmatrix}
\circ
\begin{bmatrix}
f_1^i & f_2^i & f_3^i & ... & f_N^i
\end{bmatrix} 
\\
&=
\begin{bmatrix}
h_1^i\circ f_1^i & h_1^i \circ f_2^i & h_1^i \circ f_3^i & \cdots & h^i_1\circ f^i_N \\
h_1^{i+1}\circ d_Y \circ f_1^i & h_1^{i+1}\circ d_Y \circ f_2^i & h_1^{i+1} \circ d_Y \circ f_3^i & \cdots & h_1^{i+1} \circ d_Y \circ f_N^i \\
h_1^{i+2}\circ d_Y^2 \circ f_1^i & h_1^{i+2}\circ d_Y^2 \circ f_2^i & h_1^{i+2} \circ d_Y^2 \circ f_3^i & \cdots & h_1^{i+2} \circ d_Y^2 \circ f_N^i
\\
\vdots & \vdots & \vdots & \ddots & \vdots \\
h_1^{i+N-1}\circ d_Y^{N-1} \circ f_1^i & h_1^{i+N-1}\circ d_Y^{N-1} \circ f_2^i & h_1^{i+N-1} \circ d_Y^{N-1} \circ f_3^i & \cdots & h_1^{i+N-1} \circ d_Y^{N-1} \circ f_N^i
\end{bmatrix}\\
\displaybreak[3]
&=
\begin{bmatrix}
h_1^i\circ f_1^i & h_1^i \circ f_2^i & h_1^i \circ f_3^i & \cdots & h^i_1\circ f^i_N \\
h_1^{i+1}\circ 0 & h_1^{i+1}\circ f_1^{i+1} & h_1^{i+1} \circ f_2^{i+1} & \cdots & h_1^{i+1} \circ f_{N-1}^{i+1} \\
h_1^{i+2}\circ 0 & h_1^{i+2} \circ 0 & h_1^{i+2} \circ f_3^{i+2} & \cdots & h_1^{i+2} \circ f_{N-2}^{i+2}
\\
\vdots & \vdots & \vdots & \ddots & \vdots \\
h_1^{i+N-1}\circ 0 & h_1^{i+N-1}\circ 0 & h_1^{i+N-1} \circ 0 & \cdots & h_1^{i+N-1} \circ f_1^{i+N-1}
\end{bmatrix}\\
&=
\begin{bmatrix}
h_1^i\circ f_1^i & 0 & 0 & \cdots & 0 \\
0 & h_1^{i+1}\circ f_1^{i+1} &0 & \cdots & 0\\
0 &0 & h_1^{i+2} \circ  f_1^{i+2} & \cdots &0\\
\vdots & \vdots & \vdots & \ddots & \vdots \\
0 & 0 & 0 & \cdots & h_1^{i+N-1} \circ f_1^{i+N-1}
\end{bmatrix}\\
&=
\begin{bmatrix}
\id_{X^i} & 0 & 0 & \cdots & 0 \\
0 & \id_{X^{i+1}} &0 & \cdots & 0\\
0 &0 &\id_{X^{i+2}} & \cdots &0\\
\vdots & \vdots & \vdots & \ddots & \vdots \\
0 & 0 & 0 & \cdots & \id_X{^{i+N-1}}
\end{bmatrix}\\
&= \id_{\II_N(X)}^i,
\end{align*}
where we use in the third to last row that $d_Y \circ f_j^i= f_{j-1}^{i+1}$ and $d_Y \circ f_1^i = 0$. In the second to last row we used $f_j^i\circ h_k^i = 0$ for $k \neq j$.
\item This case is just the dual of (1), but for completeness we give the proof as well.

Let $f: Y \to \PP_N(X)$ be a morphism of $N$-complexes such that there exists for all $i \in \Z$ a levelwise split $h^i: \PP_N(X)^i \to Y^i$, satisfying $f^i \circ h^i =\id_{\PP_N(X)}^i$, we now need to provide a morphism of $N$-complexes $g:\PP_N(X) \to Y$ such that $f \circ g = \id_{\PP(X)}$.\\
Again the morphisms $f^i$ and $h^i$ are of the form
\begin{equation*}
f^i=\begin{bmatrix}f^i_1&f_2&f^i_3&\cdots &f^i_N\end{bmatrix}^T \; , \; h^i=\begin{bmatrix}h^i_1&h^i_2&h^i_3&\cdots &h^i_N\end{bmatrix}.
\end{equation*}
 The fact that $f^i \circ h^i = \id_{\PP_N(X)^i}$ yields $f_j^i\circ h_j^i = \id_{X^{i-N+j+1}}$ and $f_j^i\circ h_k^i = 0$ if $k \neq j$. Since $f$ is a morphism of $N$-complexes, we know that $f^i \circ d_Y=d_{\PP_N(X)} \circ f^{i+1} $, but this implies that $ f_j^i \circ d_Y= f_{j-1}^{i+1}$ if $1 < j \le N$ and zero otherwise. Furthermore, observe that by the definition of $\PP_N(X)$ the following equation holds for $2\le j \le N $:
 \begin{equation*}
 f_j^i  \circ h_j^i= \id_{X^{i-N+j}}=f_{j-1}^{i+1} \circ h_{j-1}^{i+1}.
\end{equation*}
Hence define $g: \PP_N(X) \to Y$ dually to (1) as 
$$g^i:=\begin{bmatrix} h_1^i&d_Y \circ h_1^{i+1} & \cdots &d_Y^{N-1} \circ h_1^{i+N-1}\end{bmatrix}.$$ 
Note that this $g$ is again well defined by the definition of $\PP_N(X)$, and defines a morphism of $N$-complexes since 
\begin{align*}
(d_Y \circ g)^i& = \begin{bmatrix} d_Y \circ h_1^{i-1} &d_Y^2 \circ h_1^i &\cdots & d_Y^{N-1} \circ h_1^{i+N-3}&d_Y^N \circ h_1^{i+N-2}\end{bmatrix} \\
&=\begin{bmatrix} d_Y \circ h_1^{i-1} &d_Y^2 \circ h_1^i &\cdots & d_Y^{N-1} \circ h_1^{i+N-3}&0\end{bmatrix} \\
&= (g \circ d_{\PP_N(X)})^i.
\end{align*}

the morphism $g$ is a section in $\CC_N(\A)$, since
\begin{align*}
f \circ g &=
\begin{bmatrix}
f_1^i \\
f_2^i\\
f_3^i\\
\vdots \\
f_N^i
\end{bmatrix}
\circ
\begin{bmatrix}
h_1^i & d_Y \circ h_1^{i+1} & d_Y^2 \circ h_1^{i+2} & ... & d_Y^{N-1} \circ h_1^{i+N-1}
\end{bmatrix} 
\\
&=
\begin{bmatrix}
f_1^i\circ h_1^i & f_1^i \circ d_Y \circ h_1^{i+1} & f_1^i \circ  d_Y^2 \circ h_1^{i+2} & \cdots & f^i_1\circ d_Y^{N-1}\circ h_1^{i+N-1}\\
f_2^i \circ h_1^i & f_2^i \circ d_Y \circ h_1^{i+1}& f_2^i \circ  d_Y^2 \circ h_1^{i+2} & \cdots & f_2^i \circ  d_Y^{N-1} \circ h_1^{i+N-1} \\
f_3^i\circ h_1^i  & f_3^i \circ d_Y \circ h_1^{i+1} &f_3^i \circ  d_Y^2 \circ h_1^{i+2} & \cdots & f_3^i \circ  d_Y^{N-1} \circ h_1^{i+N-1}
\\
\vdots & \vdots & \vdots & \ddots & \vdots \\
f_N^i\circ h_1^i  & f_N^i \circ d_Y \circ h_1^{i+1} & f_N^i \circ  d_Y^2 \circ h_1^{i+2} & \cdots & f_N^i \circ  d_Y^{N-1} \circ h_1^{i+N-1}
\end{bmatrix}\\
&=
\begin{bmatrix}
f_1^i\circ h_1^i & 0 \circ h_1^{i+1} & 0 \circ h_1^{i+2} & \cdots & 0 \circ h_1^{i+N-1}\\
f_2^i \circ h_1^i & f_1^{i+1}  \circ h_1^{i+1}& 0 \circ h_1^{i+2} & \cdots & 0 \circ h_1^{i+N-1} \\
f_3^i\circ h_1^i  & f_2^{i+1}  \circ h_1^{i+1} &f_1^{i+2} \circ h_1^{i+2} & \cdots & 0 \circ h_1^{i+N-1}
\\
\vdots & \vdots & \vdots & \ddots & \vdots \\
f_N^i\circ h_1^i  & f_N^i \circ d_Y \circ h_1^{i+1} & f_N^i \circ  d_Y^2 \circ h_1^{i+2} & ... & f_N^i \circ  d_Y^{N-1} \circ h_1^{i+N-1}
\end{bmatrix}
\displaybreak[3]
\\
&=
\begin{bmatrix}
 f_1^i \circ h_1^i & 0 & 0 & \cdots & 0 \\
0 &  f_1^{i+1} \circ h_1^{i+1} &0 & \cdots & 0\\
0 &0 & f_1^{i+2} \circ h_1^{i+2}& \cdots &0\\
\vdots & \vdots & \vdots & \ddots & \vdots \\
0 & 0 & 0 & ... & f_1^{i+N-1} \circ h_1^{i+N-1}
\end{bmatrix}\\
&=
\begin{bmatrix}
\id_{X^i-N+j} & 0 & 0 & \cdots & 0 \\
0 & \id_{X^{i-N+1}} &0 & \cdots & 0\\
0 &0 &\id_{X^{i-N+2}} & \cdots &0\\
\vdots & \vdots & \vdots & \ddots & \vdots \\
0 & 0 & 0 & \cdots & \id_X{^{i}}
\end{bmatrix} 
\\&= \id_{\PP_N(X)}^i,
\end{align*}
where we use in the third to last row $ f_j^i \circ d_Y= f_{j-1}^{i+1}$ and $ f_1^i \circ d_Y= 0$. In the second to last row we use $f_j^i\circ h_k^i = 0$ for $k \neq j$.
\end{enumerate}
\end{proof}
\end{proposition}

Now the last thing we need for $\Se_\oplus$ to be a Frobenius exact structure on $\CC_N(\A)$ is that $\Se_\oplus$-injectives and $\Se_\oplus$-projectives coincide.

\begin{proposition}\label{4.6}
Let $\A$ be an additive category and $N  \in \N$. Then an object in $\CC_N(\A)$ is $S_\oplus$-projective if and only if it is $S_\oplus$-injective.
\begin{proof}
Let $P$ be $\Se_\oplus$-projective, then we have a sequence in $\Se_\oplus$ of the form:
\begin{equation*}
\tikz[heighttwo,xscale=2,yscale=2,baseline]{

\node (01) at (0.3,0.3){$0$};
\node (X1) at (1,0.3){$X$};
\node (X2) at (2,0.3){$\PP_N(P)$};
\node (X3) at (3,0.3){$P$};
\node (02) at (3.7,0.3){$0$};
    
\draw[->]
(01) edge (X1)
(X3) edge (02)
(X1) edge node[above] {} (X2)
(X2) edge node[above] {$\lambda^\PP_P$} (X3)

}.
\end{equation*}
Now since $P$ is $\Se_\oplus$-projective, this short exact sequence splits, hence $\PP_N(P)\cong P \oplus X$. But since $\PP_N(P)=\II_N(P[-N+1])$ we have that $P$ is a direct summand of an $\Se_\oplus$-injective object, and hence $\Se_\oplus$-injective itself. Here $P[-N+1]$ denotes the $N$-complex of $P$ shifted by $-N+1$.

Conversely let $I$ be an $\Se_\oplus$-injective object. Then we have the following short exact sequence in $\Se_\oplus$:
\begin{equation*}
\tikz[heighttwo,xscale=2,yscale=2,baseline]{

\node (01) at (0.3,0.3){$0$};
\node (X1) at (1,0.3){$I$};
\node (X2) at (2,0.3){$\II_N(I)$};
\node (X3) at (3,0.3){$X$};
\node (02) at (3.7,0.3){$0$};
    
\draw[->]

(01) edge (X1)
(X3) edge (02)
(X1) edge node[above] {$\lambda^\II_I$} (X2)
(X2) edge node[above] {} (X3)

}.
\end{equation*}
Now again the short exact sequence splits since $I$ is $\Se_\oplus$-injective, and we get $\II_N(I)\cong I \oplus X$. Now $\II_N(I)=\PP_N(I[N-1])$ and hence $I$ is a summand of an $\Se_\oplus$-projective object and hence itself $\Se_\oplus$-projective.
\end{proof}
\end{proposition}

Now we can put together the above properties of $\Se_\oplus$ to get:

\begin{theorem}
For an additive category $\A$ and $N  \in \N$ the exact structure $\Se_\oplus$ defines a Frobenius exact structure on the category $\CC_N(\A)$.
\begin{proof}
By Proposition~\ref{4.2}, The collection $\Se_\oplus$ defines an exact structure. It has enough $\Se_\oplus$-projectives and $\Se_\oplus$-injectives, as any object admits an admissible monomorphism $\lambda^\II_X:X \to \II_N(X)$ and an admissible epimorphism $\lambda_X^\PP :\PP_N(X) \to X$ by Proposition~\ref{4.4} and Proposition~\ref{4.5}. Finally $\Se_\oplus$-injective objects are precisely the $\Se_\oplus$-projective objects by Proposition~\ref{4.6}. So all defining criteria of Frobenius exactness are fulfilled.
\end{proof}
\end{theorem}

Having established this Frobenius exact structure we define the homotopy category of $N$-complexes and use Happel's Theorem~I.2.6 \cite{happel} to conclude that the homotopy category is well defined and triangulated in a canonical way.

\begin{remark}
We would like to remind that by 13.4 in \cite{Buhler20101} the stable category $\underline{\C}_\Se$ of a Frobenius exact category $\C_\Se$ is defined to have the same objects as $\C$ and for $X,Y \in \C$ the morphism space $\underline{\C}_\Se(X,Y):=\C(X,Y)/I$, where $I$ is the two-sided ideal of morphisms factoring through an $\Se$-injective-projective object. 
\end{remark}

\begin{definition}
Let $\A$ be an additive category and $N \in \N$. Then denote by $\KK_N(\A)$ the stable category of $\CC_N(\A)_{S_\oplus}$.
\end{definition}

\begin{corollary}\label{4.9}
The category $\KK_N(A)$ is triangulated for $\A$ an additive category and $N \in \N$. Furthermore every $\Se_\oplus$-short exact sequence 
\begin{equation*}
\tikz[heighttwo,xscale=2,yscale=2,baseline]{

\node (01) at (0.3,0){$0$};
\node (X1) at (1,0){$X$};
\node (X2) at (2,0){$Y$};
\node (X3) at (3,0){$Z$};
\node (02) at (3.7,0){$0$};
    
\draw[->]

(01) edge (X1)
(X3) edge (02)
(X1) edge node[above] {$f$} (X2)
(X2) edge node[above] {$g$} (X3)

}
\end{equation*}
induces a triangle
\begin{equation*}
\tikz[heighttwo,xscale=2,yscale=2,baseline]{

\node (X1) at (1,0){$X$};
\node (X2) at (2,0){$Y$};
\node (X3) at (3,0){$Z$};
\node (X4) at (4,0){$\Sigma X$};

\draw[->]

(X1) edge node[above] {$f$} (X2)
(X2) edge node[above] {$g$} (X3)
(X3) edge node[above] {$h$} (X4)

}
\end{equation*} and every triangle is isomorphic to a triangle of this form.
\begin{proof}
This is a direct application of Happel's Theorem~I.2.6 \cite{happel}.
\end{proof}
\end{corollary}

From now on we will often need that every morphism can be completed into a triangle in $\KK_N(\A)$. We will refer to a completion of a morphism as its cone. Since we do not need any specific structure of these cones we just fix any choice of completion.

\begin{definition}
Let  $f:X \to Y$ be a morphism in $\CC_N(\A)$ for $\A$ an additive category and $N\in \N$. Then we refer to \textbf{the cone of $f$} in $\KK_N(\A)$ as $\Cone_{K_N(\A)}(f)$, or just $\Cone(f)$ if it is clear from context. This is just any choice of completion of the morphism into a triangle
\begin{equation*}
\tikz[heighttwo,xscale=2,yscale=2,baseline]{

\node (X1) at (1,0.3){$X$};
\node (X2) at (2,0.3){$Y$};
\node (X3) at (3,0.3){$Z$};
\node (X4) at (4,0.3){$\Sigma X$};

\draw[->]

(X1) edge node[above] {$f$} (X2)
(X2) edge node[above] {$g$} (X3)
(X3) edge node[above] {$h$} (X4)

}.
\end{equation*}
\end{definition}

For a precise construction of the cone objects together with the defining morphisms arising by Happel's Theorem we refer the reader to ``Derived categories of $N$-complexes'' \cite{N-complex}.

\begin{remark}
By the definition of a triangulated category these completions are unique up to non-unique isomorphism.
\end{remark}

%% file: derived_category.tex
\section{Properties of the homotopy category of $N$-complexes}\label{section 6}

Since a goal of this work is to imitate the definition of the derived category of an abelian category, we will need some general properties of the homotopy category $\KK_N(\A)$. This will give us more control over the objects, morphisms and coproducts in $\KK_N(\A)$.

First of all we will need to be able to factor through specific $\Se$-injective and $\Se$-projective objects. This is possible for any exact structure which has enough projectives and injectives respectively. Using this property we will be able to show that coproducts in $\KK_N(\A)$ can be calculated naively and prove a more calculation based approach towards the homotopy category.

\begin{lemma} \label{chosen projective}
Let $\C_\Se$ be an exact category. Then the following statements hold:
\begin{enumerate}
\item If $X\in \C$ admits an $\Se$-injective hull $X \xrightarrow{\lambda_X} I_X$ we have that a morphism $f:X \to Y$ factors as $f:X \to I \to Y$, where $I$ is an $\Se$-injective object, if and only if $f$ factors as $f:X  \xrightarrow{\lambda_X} I_X \to Y$.
\item If $Y\in \C$ admits an $\Se$-projective resolution $P_Y \xrightarrow{\lambda_Y} Y$ we have that a morphism $f:X \to Y$ factors as $f:X \to P \to Y$, where $P$ is an $\Se$-projective object, if and only if $f$ factors as $f:X \to P_Y \xrightarrow{\lambda_Y} Y$.
\end{enumerate}
\begin{proof}
\begin{enumerate}

\item The if part is clear, since if $f:X \to Y$ factors over a chosen $\Se$-injective, it especially factors over an $\Se$-injective object.

For the only if part, let $f:X \xrightarrow{\lambda} I \xrightarrow{g} Y$ be a factorisation of a morphism $f$ such that $I$ is $\Se$-injective. Now consider the following diagram in $\C$:
\begin{equation*}
\tikz[heighttwo,xscale=1.7,yscale=1.7,baseline]{

\node (Y) at (0,0){$Y$};
\node (I) at (0,1){$I$};
\node (X) at (0,2){$X$};
\node (IX) at (1,2){$I_X$};
\node (Z1) at (2,2){$Z$};
\node (I') at (1,1){$M$};
\node (Z2) at (2,1){$Z$};
\node (po) at (0.75,1.25){$\ulcorner$};

\node[overlay] (Point) at (2.2,0){,};

\draw[->]
(X) edge [bend right=30] node[left] {$f$} (Y)
(X) edge  node[right]{$\lambda$}(I)
(X) edge node[above]{$\lambda_X$}(IX)
(IX) edge (Z1)
(I) edge  node[right]{$g$} (Y)
(IX) edge node[right]{$h$}(I')
(I) edge node[above]{$j$}(I')
(I') edge (Z2);
\draw[double distance = 2pt] 
(Z1)--(Z2);
}
\end{equation*}
where the top sequence is induced by the $\Se$-injective hull $I_X$ for $X$. The left square is a pushout, which we can form since $(X \xrightarrow{\lambda_X} I_X \to Z)\in \Se$.
Since the top sequence is in $\Se$ we get by the stability of exact structures under pushouts that also the bottom one is in $\Se$. The bottom sequence now splits since $I$ is $\Se$-injective by assumption. Call this split $s:M \to I$. Hence we get the diagram in $\C$:
\begin{equation*}
\tikz[heighttwo,xscale=1.7,yscale=1.7,baseline]{

\node (Y) at (0,0){$Y$};
\node (I) at (0,1){$I$};
\node (X) at (0,2){$X$};
\node (IX) at (1,2){$I_X$};
\node (Z1) at (2,2){$Z$};
\node (I') at (1,1){$M$};
\node (Z2) at (2,1){$Z$};
\node (po) at (0.75,1.25){$\ulcorner$};

\node[overlay] (Point) at (2.2,0){.};

\draw[->]
(X) edge [bend right=30] node[left] {$f$} (Y)
(X) edge  node[right]{$\lambda$}(I)
(X) edge node[above]{$\lambda_X$}(IX)
(IX) edge (Z1)
(I) edge  node[right]{$g$} (Y)
(IX) edge node[right]{$h$}(I')
(I) edge node[above]{$j$}(I')
(I') edge (Z2);
\draw[double distance = 2pt] 
(Z1)--(Z2);
\draw[ ->]
(I') edge [bend left=30] node[below] {$s$} (I);
}
\end{equation*}

 Using this to factor $f$ we get the claimed factorization: $$f=g \circ \lambda =g \circ s \circ j \circ \lambda =g \circ s  \circ h \circ \lambda_X.$$
 
 \item Again the if part is clear since a specific choice of an $\Se$-projective object in particular gives us an $\Se$-projective object.
 For the only if part consider again a morphism $f:X \to Y$ which factors as $f:X \xrightarrow{g} P \xrightarrow{\lambda} Y$ for some $\Se$-projective object $P$. Now consider the $\Se$-projective resolution $P_Y \xrightarrow{\lambda_Y} Y$, and the following diagram in $\C$:
 
 \begin{equation*}
\tikz[heighttwo,xscale=1.7,yscale=1.7,baseline]{

\node (Y) at (2,0){$Y$};
\node (P) at (2,1){$P$};
\node (X) at (2,2){$X$};
\node (PY) at (1,0){$P_Y$};
\node (Z1) at (0,0){$Z$};
\node (P') at (1,1){$M$};
\node (Z2) at (0,1){$Z$};
\node (pb) at (1.25,0.75){$\lrcorner$};

\node[overlay] (Point) at (2.2,0){,};

\draw[->]
(X) edge [bend left=30] node[right] {$f$} (Y)
(P) edge  node[left]{$\lambda$}(Y)
(PY) edge node[below]{$\lambda_Y$}(Y)
(Z1) edge (PY)
(X) edge  node[left]{$g$} (P)
(P') edge node[left]{$h$}(PY)
(P') edge node[below]{$j$}(P)
(Z2) edge (P');
\draw[double distance = 2pt] 
(Z1)--(Z2);
}
\end{equation*}
where the right square is a pullback, which we can form since the lower sequence is in $\Se$. Now we get again using that $\Se$ is an exact structure, that the top row is in $\Se$. So the fact that $P$ is $\Se$-projective yields that $j$ splits. Call this split $s:P \to M$. This means that we get the diagram:
\begin{equation*}
\tikz[heighttwo,xscale=1.7,yscale=1.7,baseline]{

\node (Y) at (2,0){$Y$};
\node (P) at (2,1){$P$};
\node (X) at (2,2){$X$};
\node (PY) at (1,0){$P_Y$};
\node (Z1) at (0,0){$Z$};
\node (P') at (1,1){$M$};
\node (Z2) at (0,1){$Z$};
\node (pb) at (1.25,0.75){$\lrcorner$};

\node[overlay] (Point) at (2.2,0){.};

\draw[->]
(X) edge [bend left=30] node[right] {$f$} (Y)
(P) edge  node[left]{$\lambda$}(Y)
(PY) edge node[below]{$\lambda_Y$}(Y)
(Z1) edge (PY)
(X) edge  node[left]{$g$} (P)
(P') edge node[left]{$h$}(PY)
(P') edge node[below]{$j$}(P)
(P) edge [bend right=30]node[above]{$s$}(P')
(Z2) edge (P');
\draw[double distance = 2pt] 
(Z1)--(Z2);
}
\end{equation*}
Using this diagram we can calculate $$f=\lambda \circ g = \lambda \circ j \circ s \circ g= \lambda_Y \circ h \circ s \circ g,$$
which concludes the proof.\qedhere
\end{enumerate}
\end{proof}
\end{lemma}

As mentioned before we can use this argument to show that $\KK_N(\A)$ admits coproducts and give a precise description of them.

\begin{lemma}\label{naive colimits}
Let $\F_\Se$ be a Frobenius exact category such that $\F$ admits $I$-indexed coproducts for $I$ a set. Then the stable category $\underline{\F}_\Se$ also admits $I$-indexed coproducts and we have that the canonical functor $\F \to \underline{\F}_\Se$ preserves coproducts,  i.\,e.\ the morphisms exhibiting the coproduct in $\F$ are representatives of the morphisms exhibiting the coproduct in $\underline{\F}_\Se$.
\begin{proof}
We will show that for a coproduct $\left(\bigoplus\limits_{i \in I}X_i,\iota_i\right)$ in $\F$ the maps $\left[\iota_i\right]$ exhibit $\bigoplus\limits_{i\in I} X_i$ as a coproduct in $\underline{\F}_\Se$.

First consider the coproduct $\left(\bigoplus\limits_{i \in I}X_i,\iota_i\right)$ in $\F$. The $\iota_i$ induces morphisms $[\iota_i]:X_i \to \bigoplus\limits_{i\in I}X_i$ in $\underline{\F}_\Se$. Now consider for each $i \in I$ a morphism $\psi_i:X_i \to T$ in $\underline{\F}_\Se$ and pick for each of these a representative $\varphi_i$ such that $[\varphi_i]=\psi_i$, which are morphisms $\varphi_i:X_i \to T$ in $\F$. So we get by the universal property of the coproduct a unique morphism $\varphi: \bigoplus\limits_{i\in I} X_i \to T$, such that the following diagram commutes for each $i \in I$:
\begin{equation*}
\tikz[heighttwo,xscale=2,yscale=2,baseline]{

\node (Xi) at (0,0.8){$X_i$};
\node (X) at (1,0.8){$\bigoplus X_i$};
\node (T)  at (0.5,0) {$T.$};
    
\draw[->]
(Xi) edge node[above] {$\iota_i$} (X)
(Xi) edge node[left] {$\varphi_i$} (T)
(X) edge node[right] {$\varphi$} (T);
}
\end{equation*}
This implies by the definition of the stable category that the same diagram commutes in $\underline{\F}_\Se$.

Hence it just remains to show uniqueness of $\varphi$. Since $\F$ is additive it actually suffices to show that if we have a morphism $f:\bigoplus\limits_{i\in I}X_i \to T$ in $\underline{\F}_\Se$ such that all $f \circ [\iota_i]$ vanish, so does $f$. Hence assume $f \circ [\iota_i]=0$ and choose representatives $ [f_i]=f \circ [\iota_i]$. These $f_i$ now factor by definition of the stable category and Proposition~\ref{chosen projective} as $f_i = \lambda_T \circ \psi_i$ for some $\psi_i:X_i \to P_T$ and $\lambda_T: P_T \to T$ some $\Se$-projective resolution of $T$. So we get for all $i \in I$ a commutative diagram in $\F$
\begin{equation*}
\tikz[heighttwo,xscale=1.7,yscale=1.7,baseline]{

\node (Xi) at (0,1.6){$X_i$};
\node (X) at (2,1.6){$\bigoplus X_i$};
\node (T)  at (1,0) {$T.$};  
\node (P) at (1,1){$P_T$};
\draw[->]
(Xi) edge node[above right] {$\psi_i$}(P)
(P) edge node[above right] {$\lambda_T$} (T)
(Xi) edge node[above] {$\iota_i$} (X)
(Xi) edge node[left] {$f_i$} (T)
(X) edge node[right] {$f$} (T);
}
\end{equation*}
By the universal property of the coproduct this induces a map $\psi: \bigoplus\limits_{i\in I} X_i \to P_T$, such that the top, the left and the outer triangle commute in the following diagram in $\F$:
\begin{equation*}
\tikz[heighttwo,xscale=1.7,yscale=1.7,baseline]{

\node (Xi) at (0,1.6){$X_i$};
\node (X) at (2,1.6){$\bigoplus X_i$};
\node (T)  at (1,0) {$T.$};  
\node (P) at (1,1){$P_T$};
\draw[->]
(X) edge node[above left] {$\psi$}(P)
(Xi) edge node[above right] {$\psi_i$}(P)
(P) edge node[above right] {$\lambda_T$} (T)
(Xi) edge node[above] {$\iota_i$} (X)
(Xi) edge node[left] {$f_i$} (T)
(X) edge node[right] {$f$} (T);
}
\end{equation*}
This implies that the whole diagram commutes. The commutativity of this diagram now just means $[f] =0$ in $\underline{\F}_\Se$, which concludes the proof.
\end{proof}
\end{lemma}

\begin{remark}
One can show an analogous statement for products. There one just needs to use $\Se$-injective hulls instead of $\Se$-projective resolutions. However this is not important for the rest of our constructions, and so we omit it.
\end{remark}

We are also interested in having a more calculation based approach towards the structure of the homotopy category.

\begin{definition}
Let $\A$ be an additive category, $N \in \N$ and let $f:X\to Y$ be a morphism of $N$-complexes in $\A$. Then $f$ is said to be \textbf{nullhomotopic} if there exists for all $j \in \Z$ an $s^j \in \A(X^j,Y^{j-N+1})$ such that for all $i \in \Z$ we get:
\begin{equation*}
f^i=\sum_{j=0}^{N-1} d^{N-j-1}\circ s^{i+j} \circ d^{j}.
\end{equation*}
\end{definition}

This definition should of course be equivalent to the factoring of a morphism over an $\Se_\oplus$-injective-projective object.

\begin{lemma}\label{homotopy formula}

Let $\A$ be an additive category and $f:X\to Y$ be a morphism in $\CC_N(\A)$. Then the following statements are equivalent:
\begin{enumerate}
\item the morphism $f$ is nullhomotopic,
\item the morphism $f$ factors through an $\Se_\oplus$-injective-projective object.
\end{enumerate}
\begin{proof}
We show the implications separately:
\begin{description}
\item[\normalfont (1)$\implies$(2)]
By definition of a nullhomotopic morphism there exist $s^t \in \A (X^t,Y^{t-N+1})$ for all $t \in \Z$ such that 
\begin{equation*}
f^i=\sum_{j=0}^{N-1} d^{N-j-1}\circ s^{i+j} \circ d^{j}.
\end{equation*}
Consider the map $\lambda^I_N: X \to \II_N(X)$. The maps $s^i: X^i \to Y^{i-N+1}$ induce maps $\hat{s}^i: \mu^i_N(X^i)^i \to Y$ via the adjunction $\CC_N(\A)(\mu_N^i (X^i), Y) \cong \A (X^i,Y^{i-N+1})$ from Proposition~\ref{functoriality mu}. But this gives a map $\hat{s}: \II_N(X) \to Y$ since the source is just a direct sum of $\mu_N^i(X^i)$. Unravelling the definition of $\II_N (X)$, $\lambda_X^\II$ and the adjunction one gets:
\begin{align*}
(\lambda^\II_X)^i&=\begin{bmatrix}\id & d & d^2 & \cdots & d^{N-2}& d^{N-1} \end{bmatrix}^T :
 X^i \to X^i\oplus X^{i+1}\oplus \cdots \oplus X^{i+N-2}\oplus X^{i+N-1})  \\
\hat{s}^i&= \begin{bmatrix} d^{N-1} \circ s^i & d^{N-2} \circ s^{i+1} & d^{N-3} \circ s^{i+2} & \cdots  & s^{i+N-1} \end{bmatrix} :
X^i\oplus X^{i+1}\oplus  \cdots \oplus X^{i+N-1}  \to  Y^i ).
\end{align*}
Composing these morphisms yields now:
\begin{equation*}
(\hat{s} \circ \lambda_X^\II)^i = \sum_{j=0}^{N-1} d^{N-j-1}s^{i+j}d^{j} = f^i.
\end{equation*}
So since $\II_N(X)$ is $\Se_\oplus$-injective-projective this proves the claim.
\item[\normalfont (2)$\implies$(1)]
Let $f: X \to Y$ factor as $f=s\circ g$ where $g: X \to P$, $s:P \to Y$ and $P$ is $\Se_\oplus$-injective-projective. Then we may assume by Lemma~\ref{chosen projective} that $P \cong \II_N(X)$ and $g= \lambda^\II_X$. Now the claim follows cf. Remark~\ref{homotopy summand}, since we have that $f$ factors as $\lambda^\II_X \circ s$ where $\lambda^\II_X:X \to \II(X)$ and $s: \II_N(X) \to Y$. More precisely these $s^i$ now define levelwise maps $s^i: X^i\oplus X^{i+1}\oplus X^{i+2}...\oplus X^{i+N-1}  \to  Y^i$, the fact that these commute with the differential of $\II_N(X)$ yields
$$s^i\circ \iota_{X^j} = s^{i+1} \circ \iota_{X^{j}} \text{ for }0<j-i < N.$$
Here $\iota_{X^j}:X^j \hookrightarrow \bigoplus\limits_{l=i}^{\mathclap{i+N-1}}X^l$ is the canonicalinclusion for $i \le j \le i+N-1$.

 This means that $s: \II_N(X) \to Y$ is equivalently defined by $\widehat{s^i}: X^i \to Y^{i-N+1}$ and 
\begin{equation*}
s^i= \begin{bmatrix} d^{N-1} \circ \widehat{s}^i & d^{N-2} \circ \widehat{s}^{i+1} & d^{N-3} \circ \widehat{s}^{i+2} & \cdots  & \widehat{s}^{i+N-1} \end{bmatrix} :
X^i\oplus X^{i+1}\oplus X^{i+2} \cdots \oplus X^{i+N-1}  \to  Y^i .
\end{equation*} So we get 
\begin{equation*}
(s \circ \lambda_X^\II)^i = \sum_{j=0}^{N-1} d^{N-j-1}\circ \widehat{s}^{i+j}\circ d^{j} = f^i,
\end{equation*}
as claimed. \qedhere
\end{description}
\end{proof}
\end{lemma}

\begin{remark}\label{homotopy summand}
The construction in terms of adjunctions can be seen more intuitively using the following diagram in $\A$, in which we only consider one summand of $\mu_N^i (X)$:
\begin{equation*}
\tikz[heightone,xscale=2.5,yscale=2.5,baseline]{

\node (P00) at (-1,0){$\cdot \cdot \cdot$}; 
\node (P10) at (3,0){$\cdot \cdot \cdot$}; 
\node (P20) at (7,0){$\cdot \cdot \cdot$}; 
\node (P01) at (-1,1){$\cdot \cdot \cdot$}; 
\node (P11) at (3,1){$\cdot \cdot \cdot$}; 
\node (P21) at (7,1){$\cdot \cdot \cdot$}; 
\node (P02) at (-1,2){$\cdot \cdot \cdot$}; 
\node (P12) at (3,2){$\cdot \cdot \cdot$}; 
\node (P22) at (7,2){$\cdot \cdot \cdot$}; 

 \node (X_i-N-1) at (0,0){$X^{i-N}$};
 \node (X_i-N) at (1,0){$X^{i-N+1}$};
 \node (X_i-N+1) at (2,0){$X^{i-N+2}$};

  \node (X_i-1) at (4,0){$X^{i-1}$};
   \node (X_i) at (5,0){$X^{i}$};
    \node (X_i+1) at (6,0){$X^{i+1}$};

     \node (Xi_N-i-1) at (0,1){$0$};
 \node (Xi_N-i) at (1,1){$X^i$};
 \node (Xi_N-i+1) at (2,1){$X^i$};
 
  \node (Xi_i-1) at (4,1){$X^i$};
   \node (Xi_i) at (5,1){$X^i$};
    \node (Xi_i+1) at (6,1){$0$};

    \node (Y_i-N-1) at (0,2){$X^{i-N}$};
 \node (Y_i-N) at (1,2){$X^{i-N+1}$};
 \node (Y_i-N+1) at (2,2){$X^{i-N+2}$};
 
  \node (Y_i-1) at (4,2){$X^{i-1}$};
   \node (Y_i) at (5,2){$X^{i}$};
    \node (Y_i+1) at (6,2){$X^{i+1}$};
    
\draw[->]

(Y_i-N-1) edge node[right] {$0$} (Xi_N-i-1)
(Y_i-N) edge node[right] {$d^{N-1}$} (Xi_N-i)
(Y_i-N+1) edge node[right] {$d^{N-2}$} (Xi_N-i+1)
(Y_i-1) edge node[right] {$d$} (Xi_i-1)
(Y_i) edge node[right] {$\id_{X^i}$} (Xi_i)
(Y_i+1) edge node[right] {$0$} (Xi_i+1)

(Xi_N-i-1) edge node[right] {$0$} (X_i-N-1)
(Xi_N-i) edge node[right] {$s^i$} (X_i-N)
(Xi_N-i+1) edge node[right] {$d\circ s^i$} (X_i-N+1)
(Xi_i-1) edge node[right] {$d^{N-2}\circ s^i$} (X_i-1)
(Xi_i) edge node[right] {$d^{N-1}\circ s^i$} (X_i)
(Xi_i+1) edge node[right] {$0$} (X_i+1)

(P02) edge node[above] {$d$} (Y_i-N-1)
(Y_i-N-1) edge node[above] {$d$} (Y_i-N)
(Y_i-N) edge node[above] {$d$} (Y_i-N+1)
(Y_i-N+1) edge node [above] {$d$} (P12)

(P12) edge node[above] {$d$} (Y_i-1)
(Y_i-1) edge node[above] {$d$} (Y_i)
(Y_i) edge node[above] {$d$} (Y_i+1)
(Y_i+1) edge node[above] {$d$} (P22)

(P01) edge node[above] {$0$} (Xi_N-i-1)
(Xi_N-i-1) edge node[above] {$0$} (Xi_N-i)
(Xi_N-i) edge node[above] {$\id_{X^i}$} (Xi_N-i+1)
(Xi_N-i+1) edge node[above] {$\id_{X^i}$} (P11)

(P11) edge node[above] {$\id_{X^i}$} (Xi_i-1)
(Xi_i-1) edge node[above] {$\id_{X^i}$} (Xi_i)
(Xi_i) edge node[above] {$0$} (Xi_i+1)
(Xi_i+1) edge node[above] {$0$} (P21)

(P00) edge node[above] {$d$} (X_i-N-1)
(X_i-N-1) edge node[above] {$d$} (X_i-N)
(X_i-N) edge node[above] {$d$} (X_i-N+1)
(X_i-N+1) edge node[above] {$d$} (P10)

(P10) edge node[above] {$d$} (X_i-1)
(X_i-1) edge node[above] {$d$} (X_i)
(X_i) edge node[above] {$d$} (X_i+1)
(X_i+1) edge node[above] {$d$} (P20)
;

;

}.
\end{equation*}
By considering the composition of the two vertical morphisms one gets the desired summands for the homotopy equation.
\end{remark}

\section{Two approaches to the derived category of $N$-complexes}\label{section 7}

After this work on the homotopy category of $N$-complexes and its properties we are finally able to inspect the notion of cohomology of $N$-complexes. In particular this yields the notions of $N$-quasi-isomorphisms and the derived category of $N$-complexes. As in the case of chain complexes this category is defined as a Verdier quotient of the homotopy category of $N$-complexes. However, in chain complexes over a ring one can realize the derived category already as a subcategory of the homotopy category, namely as the full subcategory of h-projectives. 

The aim of this section is to understand the notion of $N$-cohomology and $N$-quasi-isomorphisms and use this in the definition of the derived category of $N$-complexes and the property of being h-projective. Finally we will show that for the category of modules over a ring the full subcategory of the homotopy category consisting of h-projective $N$-complexes already realizes the derived category. To prove this we will use a construction of semifree complexes which is inspired by an argument of Drinfeld \cite{drinfeld}.

Each of the approaches mentioned above has its advantages. For example the construction as a Verdier quotient is more straightforward, intuitive and comes with a canonical functor. On the other hand the construction using h-projectives allows better control over the morphisms in the derived category. This is due to the fact that all morphisms are already contained in the homotopy category and calculations using roof diagrams are not needed. 

In Section 8 we will construct another realization of the derived category of $N$-complexes using model structures. The advantage of this is that it shows that the derived category is indeed the Gabriel-Zismann localisation along certain morphisms.

In order to get a first notion of what the derived category of $N$-complexes should be we need to define the cohomology of an $N$-complex. Furthermore we will show that it still defines a functor just as in the usual case. A problem however is that the cohomology of $N$-complexes consists of more than one object in each degree, but we will see that we loose less than one might think due to that.

\begin{definition}\label{5.1}
Let $\A$ be an abelian category, $N \in \N$ and $X \in \CC_N(\A)$. Then the $r$\textbf{-th cohomology in degree} $i$, for $1 \le r \le N-1$ and $i \in \Z$, is the following object in $\A$:
\begin{equation*}
H_{(r)}^i(X):= \ker(d^r)^i/\im(d^{N-r})^i.
\end{equation*}
\end{definition}

\begin{proposition}\label{5.2}
Definition~\ref{5.1} can be extended to an additive functor $$H_{(r)}^i(\_):\CC_N(\A) \to \A.$$
\begin{proof}
Let $f:X \to Y$ be a morphism in $\CC_N(\A)$. By definition $f$ commutes with the differentials, hence $f\circ d_X^r=d_Y^r \circ f$. So $f$ induces a unique map $\ker(d^r)^i(f) :\ker(d_X^r)^i \to \ker(d^r_Y)^i $ by the universal property of $\ker(d^r_Y)$ such that the diagram 
\begin{equation*}
\tikz[heighttwo,xscale=2,yscale=2,baseline]{

\node (K) at (1,0){$\ker(d_Y^r)^i$};
\node (Yi) at (2,0){$Y^i$};
\node (Yi+r) at (3,0){$Y^{i+r}$};
\node (T)  at (2,-1) {$\ker(d_X^r)^i$};

\draw[->]

(T) edge node[right] {$f^i|_{\ker(d_Y^r)}$} (Yi)
(Yi) edge node[above] {$d^r$} (Yi+r);
\draw[right hook ->]
(K) to (Yi);
\draw[dashed, ->]
(T) edge node[below left] {$\exists ! \ker(d^r)^i(f)$}(K);
}
\end{equation*} commutes, since 
\begin{equation*}
d^r \circ f|_{\ker(d_X^r)} = f \circ d_X^r |_{\ker(d_X^r)}=0.
\end{equation*}
Hence to get a functor from $N$-complexes to the ground category it suffices to check that $f(\im(d_X^{N-r})) \subset \im(d_Y^{N-r})$, but this is satisfied, since by definition $f$ commutes with the differential.

Now it remains to check that this indeed is an additive functor. But for this we only need to check that the assignment $f \mapsto \ker(d^r)^i(f)$ is additive. Hence let $f=g +g'$. In this case both $\ker(d^r)^i(f)$ and $\ker(d^r)^i(g)+\ker(d^r)^i(g')$ would make the diagram commute, and hence they have to coincide by the universal property of the kernel.
\end{proof}
\end{proposition}

Having defined the notion of cohomology, it is natural to define the notion of an $N$-quasi-isomorphism:

\begin{definition}
Let $\A$ be an abelian category and $N \in \N$. We call a morphism $f: X\to Y$ in $\CC_N(\A)$ an \textbf{$N$-quasi-isomorphism} if $H_{(r)}^i(f)$ is an isomorphism for all $i \in \Z$ and $1 \le r \le N-1$.
\end{definition}

With these basic definitions we will now show that cohomology defines a functor on the category $\KK_N(\A)$ as well.

\begin{proposition}\label{5.4}
Let $\A$ be an abelian category, $N \in \N$ and $X\in \CC_N(\A)$ such that $X\cong 0$ in the category $\KK_N(\A)$. Then $H_{(r)}^i(X)=0$ for all $i \in \Z$ and $1 \le r \le N-1$.
\begin{proof}
The assumption $X \cong 0$ in $\KK_N(\A)$ implies that $\id_X$ factors through $\PP_N(X)$. Hence $X$ is a direct summand of $\PP_N(X)$ in $\CC_N(\A)$. So by Proposition~\ref{5.2} it is enough to show for all $X \in \CC_N(\A)$ that $H^i_{(r)}(\PP_N(X))=0$. Now by the definition of $\PP_N(X)$ we have
\begin{equation*}
\PP_N(X):=\bigoplus_{j \in \Z} \mu_N^j(X^j).
\end{equation*}
So it suffices to check the claim for $\mu_N^j(X^j)$. Now observe that $\ker(d_{\mu_N^j(X^j)}^r)^i=X^j$ if $j-r < i \le j$, and zero else. So we just need to consider the case $j-r < i \le j$, but then $$\im(d_{\mu_N^j(X^j)}^{N-r})^i=X^j=\ker(d_{\mu_N^j(X^j)}^r)^i$$ and hence $H_{(r)}^i(\mu_N^j(X^j))=0$. This finishes the proof since $H^i_{(r)}(\_)$ is additive.
\end{proof}
\end{proposition}

\begin{corollary}\label{5.5}
Let $\A$ be an abelian category, $N \in \N$, $i \in \Z$ and $1 \le r \le N-1$. Then $H_{(r)}^i(\_)$ induces an additive functor $H_{(r)}^i(\_):\KK_N(\A) \to \A$.
\begin{proof}
By Proposition~\ref{5.2} and Proposition~\ref{5.4} the functor $H_{(r)}^i(\_):\CC_N(\A)\to \A$ is additive and sends all morphisms to zero that factor through an $\Se_\oplus$-injective-projective object. But $\KK_N(\A)$ is the stable category of $\CC_N(\A)_{\Se_\oplus}$ and hence $H_{(r)}^i(\_)$ factors through $\KK_N(\A)$.
\end{proof}
\end{corollary}

\begin{corollary}\label{5.6}
Let $\A$ be an abelian category, $N \in \N$ and $f: X \to Y$ a morphism in $\CC_N(\A)$, such that $f$ induces an isomorphism in $\KK_N(\A)$. Then $f$ is an $N$-quasi-isomorphism.
\begin{proof}
Obvious by Corollary~\ref{5.5}.
\end{proof}
\end{corollary}

In the case of $N=2$ the long exact sequence in cohomology is a powerful tool to calculate cohomology.

\begin{lemma}\label{5.7}
Let $\A$ be an abelian category and 
\begin{equation*}
\tikz[heighttwo,xscale=2,yscale=2,baseline]{

\node (01) at (0,0.3){$0$};
\node (X) at (1,0.3){$X$};
\node (Y) at (2,0.3){$Y$};
\node (Z) at (3,0.3){$Z$};
\node (02) at (4,0.3){$0$};

\draw[->]
(01) edge (X)
(X) edge node[above]{$f$} (Y)
(Y) edge node[above]{$g$} (Z)
(Z) edge (02);
}
\end{equation*}
a short exact sequence in $\CC_N(\A)$. Then there exists for every $1 \le r \le N-1$ a long exact sequence of the form:
\begin{equation*}
\tikz[heighttwo,xscale=4,yscale=2,baseline]{

\node (P1) at (2,3){$\cdots$};
\node (HZi-N) at (3,3){$H^{i-N}_{(r)}(Z)$};
\node (HXi-N+r) at (1,2) {$H^{i-N+r}_{(N-r)}(X)$};
\node (HYi-N+r) at (2,2) {$H^{i-N+r}_{(N-r)}(Y)$};
\node (HZi-N+r) at (3,2) {$H^{i-N+r}_{(N-r)}(Z)$};
\node (HXi) at (1,1){$H^i_{(r)}(X)$};
\node (HYi) at (2,1){$H^i_{(r)}(Y)$};
\node (HZi) at (3,1){$H^i_{(r)}(Z)$};
\node (HXi+r) at (1,0){$H^{i+r}_{(r)}(X)$};
\node (P2) at (2,0){$\cdots$};

\node[overlay] (Point) at (3.2,0){.};

\draw[overlay,->]
(P1) edge (HZi-N)
(HXi-N+r) edge node[above]{$H^{i-N+r}_{(N-r)}(f)$} (HYi-N+r)
(HYi-N+r) edge node[above]{$H^{i-N+r}_{(N-r)}(g)$} (HZi-N+r)
(HXi) edge node[above]{$H^i_{(r)}(f)$} (HYi)
(HYi) edge node[above]{$H^i_{(r)}(g)$} (HZi)
(HXi+r) edge (P2)
(HZi-N) edge[out=355,in=175]  (HXi-N+r)
(HZi-N+r) edge[out=355,in=175] (HXi)
(HZi) edge[out=355,in=175]  (HXi+r);
}
\end{equation*}
\begin{proof}


Fix $1 \le r \le N-1$. Then we may build a short exact sequence
\begin{equation*}
\tikz[heighttwo,xscale=2,yscale=2,baseline]{

\node (01) at (0,0){$0$};
\node (X) at (1,0){$\widehat{X}$};
\node (Y) at (2,0){$\widehat{Y}$};
\node (Z) at (3,0){$\widehat{Z}$};
\node (02) at (4,0){$0$};

\draw[->]
(01) edge (X)
(X) edge node[above]{$\widehat{f}$} (Y)
(Y) edge node[above]{$\widehat{g}$} (Z)
(Z) edge (02);
}
\end{equation*}
in $\CC_2(\A)$ by considering the following diagram in $\A$:
\begin{equation*}
\tikz[heighttwo,xscale=3,yscale=1.9,baseline]{

\node(d-1)[overlay] at (0,2){\tiny(0)};
\node(d0)[overlay] at (0,3){\tiny(-1)};
\node(d1)[overlay] at (0,1){\tiny(1)};

\node (PX1) at (1,4){$\vdots$};
\node (PY1) at (2,4){$\vdots$};
\node (PZ1) at (3,4){$\vdots$};
\node (Xi-N+r) at (1,3){$X^{i-N+r}$};
\node (Yi-N+r) at (2,3){$Y^{i-N+r}$};
\node (Zi-N+r) at (3,3){$Z^{i-N+r}$};
\node (Xi) at (1,2){$X^i$};
\node (Yi) at (2,2){$Y^i$};
\node (Zi) at (3,2){$Z^i$};
\node (Xi+r) at (1,1){$X^{i+r}$};
\node (Yi+r) at (2,1){$Y^{i+r}$};
\node (Zi+r) at (3,1){$Z^{i+r}$};
\node (PX2) at (1,0){$\vdots$};
\node (PY2) at (2,0){$\vdots$};
\node (PZ2) at (3,0){$\vdots$};

\draw[overlay,->]
(PX1) edge node[left] {}(Xi-N+r)
(Xi-N+r) edge node[left] {$d_X^{N-r}$}(Xi)
(Xi) edge node[left] {$d_X^{r}$}(Xi+r)
(Xi+r) edge node[left] {}(PX2)
(PY1) edge node[left] {}(Yi-N+r)
(Yi-N+r) edge node[left] {$d_Y^{N-r}$}(Yi)
(Yi) edge node[left] {$d_Y^{r}$}(Yi+r)
(Yi+r) edge node[left] {}(PY2)
(PZ1) edge node[left] {}(Zi-N+r)
(Zi-N+r) edge node[left] {$d_Z^{N-r}$}(Zi)
(Zi) edge node[left] {$d_Z^{r}$}(Zi+r)
(Zi+r) edge node[left] {}(PZ2)
(Xi-N+r) edge node[above]{$f^{i-N+r}$} (Yi-N+r)
(Yi-N+r) edge node[above]{$g^{i-N+r}$} (Zi-N+r)
(Xi) edge node[above]{$f^i$} (Yi)
(Yi) edge node[above]{$g^i$} (Zi)
(Xi+r) edge node[above]{$f^{i+r}$} (Yi+r)
(Yi+r) edge node[above]{$g^{i+r}$} (Zi+r);
}.
\end{equation*}
Now we get a long exact sequence in ordinary cohomology of the form
\begin{equation*}
\tikz[heighttwo,xscale=4,yscale=2,baseline]{

\node (P1) at (2,3){$\cdots$};
\node (HZi-N) at (3,3){$H^{-2}_{}(\widehat{X})$};
\node (HXi-N+r) at (1,2) {$H^{-1}_{}(\widehat{X})$};
\node (HYi-N+r) at (2,2) {$H^{-1}_{}(\widehat{Y})$};
\node (HZi-N+r) at (3,2) {$H^{-1}_{}(\widehat{Z})$};
\node (HXi) at (1,1){$H^0_{}(\widehat{X})$};
\node (HYi) at (2,1){$H^0_{}(\widehat{Y})$};
\node (HZi) at (3,1){$H^0_{}(\widehat{Z})$};
\node (HXi+r) at (1,0){$H^{1}_{}(\widehat{X})$};
\node (P2) at (2,0){$\cdots$};

\node[overlay] (Point) at (3.2,0){.};
    
\draw[overlay,->]
(P1) edge (HZi-N)
(HXi-N+r) edge node[above]{$H^{-1}_{}(\widehat{f})$} (HYi-N+r)
(HYi-N+r) edge node[above]{$H^{-1}_{}(\widehat{g})$} (HZi-N+r)
(HXi) edge node[above]{$H^0_{}(\widehat{f})$} (HYi)
(HYi) edge node[above]{$H^0_{}(\widehat{g})$} (HZi)
(HXi+r) edge (P2)
(HZi-N) edge[out=355,in=175]  (HXi-N+r)
(HZi-N+r) edge[out=355,in=175]  (HXi)
(HZi) edge[out=355,in=175](HXi+r);
}
\end{equation*}
But this is just the desired sequence by the definition of $\widehat{X}$, $\widehat{Y}$ and $\widehat{Z}$.
\end{proof}
\end{lemma}

\begin{corollary}\label{5.8}
Let $\A$ be an abelian category, $N \in \N$  and 
\begin{equation*}
\tikz[heighttwo,xscale=2,yscale=1,baseline]{

\node (X) at (1,0){$X$};
\node (Y) at (2,0){$Y$};
\node (Z) at (3,0){$Z$};
\node (X1) at (4,0) {$\Sigma(X)$};

\draw[->]
(X) edge node[above]{$f$} (Y)
(Y) edge node[above]{$g$} (Z)
(Z) edge node[above]{$h$} (X1);
}
\end{equation*}
be a triangle in $\KK_N(\A)$. Then this triangle induces long exact sequences like in Lemma~\ref{5.7}.
\begin{proof}
The triangle comes by Corollary~\ref{4.9} from an $\Se_\oplus$-exact sequence, hence we get the long exact sequence for this short exact sequence by Lemma~\ref{5.7}. Since all objects in this short exact sequence are isomorphic in $\KK_N(\A)$ to the objects in our triangle we get the statement by Corollary~\ref{5.6}.
\end{proof}
\end{corollary}

One of the most important features of the long exact sequence in ordinary cohomology is a characterisation of a $2$-quasi-isomorphism by its cone. This is still possible in the general case.

\begin{corollary}\label{5.9}
Let $\A$ be an abelian category, $N \in \N$ and let $f: X\to Y$ be a morphism in $\CC_N(\A)$. Then $f$ is an $N$-quasi-isomorphism if and only if $H_{(r)}^i(\Cone(f))=0$ for all $i\in \Z$ and $1 \le r \le N-1$.
\begin{proof}
We have a triangle of the form 
\begin{equation*}
\tikz[heighttwo,xscale=2,yscale=2,baseline]{

\node (X) at (1,0.3){$X$};
\node (Y) at (2,0.3){$Y$};
\node (Z) at (3,0.3){$\Cone(f)$};
\node (X1) at (4,0.3) {$\Sigma(X)$};

\draw[->]
(X) edge node[above]{$f$} (Y)
(Y) edge node[above]{$g$} (Z)
(Z) edge node[above]{$h$} (X1);
},
\end{equation*}
hence we get long exact sequences 
\begin{equation*}
\tikz[heighttwo,xscale=1,yscale=2,baseline]{

\node (P1) at (4,3){$\cdots$};
\node (HZi-N) at (8,3){$H^{i-N}_{(r)}(\Cone(f))$};
\node (HXi-N+r) at (0,2) {$H^{i-N+r}_{(N-r)}(X)$};
\node (HYi-N+r) at (4,2) {$H^{i-N+r}_{(N-r)}(Y)$};
\node (HZi-N+r) at (8,2) {$H^{i-N+r}_{(N-r)}(\Cone(f))$};
\node (HXi) at (0,1){$H^i_{(r)}(X)$};
\node (HYi) at (4,1){$H^i_{(r)}(Y)$};
\node (HZi) at (8,1){$H^i_{(r)}(\Cone(f))$};
\node (HXi+r) at (0,0){$H^{i+r}_{(r)}(X)$};
\node (P2) at (4,0){$\cdots$};

\node[overlay] (Point) at (9.5,0){.};

\draw[overlay,->]
(P1) edge (HZi-N)
(HXi-N+r) edge node[above]{$H^{i-N+r}_{(N-r)}(f)$} (HYi-N+r)
(HYi-N+r) edge node[above]{$H^{i-N+r}_{(N-r)}(g)$} (HZi-N+r)
(HXi) edge node[above]{$H^i_{(r)}(f)$} (HYi)
(HYi) edge node[above]{$H^i_{(r)}(g)$} (HZi)
(HXi+r) edge (P2)
(HZi-N) edge[out=355,in=175](HXi-N+r)
(HZi-N+r) edge[out=355,in=175] (HXi)
(HZi) edge[out=355,in=175] (HXi+r);
}
\end{equation*}
But by the exactness of these sequences we get that $H^i_{(r)}(f)$ is an isomorphism for all $i \in \Z$ and $1 \le r \le N-1$, if and only if $H^i_{(r)}(\Cone(f))=0$ for all $i \in \Z$ and $1 \le r \le N-1$.
\end{proof}
\end{corollary}

\subsection{Verdier quotient}
$\;$

Having established these general tools for cohomology of $N$-complexes we can start to  construct the derived category in the way of Osamu Iyama, Kiriko Kato and Jun-Ichi Miyachi \cite{N-complex}. This construction is quite straightforward if one keeps in mind that the objects we want to divide out are the ones with vanishing cohomology. Now it just remains to check that these are stable enough to define a Verdier quotient.

\begin{definition}
Let $\A$ be an abelian category and $N \in \N$. We call an object $X \in \CC_N(\A)$ \textbf{acyclic}, if $H^i_{(r)}(X)=0$ for all $i \in \Z$ and $1 \le r \le N-1$. We refer to the full subcategory of $\KK_N(\A)$ consisting of acyclic objects as $\KK_N(\A)^\emptyset$.
\end{definition}

\begin{proposition}\label{acyclics thick triangulated}
Let $\A$ be an abelian category and $N \in \N$. Then $\KK_N(\A)^\emptyset$ is a thick triangulated subcategory of $\KK_N(\A)$.
\begin{proof} We check the two properties separately:
\begin{description}
\item[\normalfont thick] 
Let $X \in \KK_N(\A)$ be a direct summand of $Y \in \KK_N(\A)^\emptyset$. Then there exists a $Z$ such that $X \oplus Z \cong Y$, but this implies since $H^i_{(r)}$ is additive for all $1 \le r \le N$ and $i \in \Z$ that $$H^i_{(r)}(X)\oplus H^i_{(r)}(Z) \cong H^i_{(r)}(Y)\cong 0,$$ hence $H^i_{(r)}(X)=0$, so we get $X \in \KK_N(\A)^\emptyset$.

\item[\normalfont triangulated]
 We need to check that for all triangles
\begin{equation*}
\tikz[heighttwo,xscale=2,yscale=2,baseline]{

\node (X) at (1,0){$X$};
\node (Y) at (2,0){$Y$};
\node (Z) at (3,0){$Z$};
\node (X1) at (4,0) {$\Sigma(X)$};

\draw[->]
(X) edge node[above]{$f$} (Y)
(Y) edge node[above]{$g$} (Z)
(Z) edge node[above]{$h$} (X1);
}
\end{equation*}
 in $\KK_N(\A)$ such that $X$ and $Z$ are acyclic, also $Y$ is acyclic, and that if $X \in \KK_N(\A)^\emptyset$, also $\Sigma(\KK_N(\A))$ and $\Sigma^{-1}(\KK_N(\A))$ are objects in $\KK_N(\A)^\emptyset$.
 
 First consider the triangle above, then we get for all $1 \le r \le N-1$ and $i \in \Z$ by Corollary~\ref{5.8} a long exact sequence of the form 
\begin{equation*}
\tikz[heighttwo,xscale=4,yscale=2,baseline]{

\node (P1) at (2,2){$\cdots$};
\node (HZi-N) at (3,2){$0$};
\node (HXi-N+r) at (1,1) {$0$};
\node (HYi-N+r) at (2,1) {$H^{i-N+r}_{(N-r)}(Y)$};
\node (HZi-N+r) at (3,1) {$0$};
\node (HXi) at (1,0){$0$};
\node (HYi) at (2,0){$H^i_{(r)}(Y)$};
\node (HZi) at (3,0){$0$};
\node (HXi+r) at (1,-1){$0$};
\node (P2) at (2,-1){$\cdots$};

\node (point) at (3.2,-1){.};
    
\draw[overlay,->]
(P1) edge (HZi-N)
(HXi-N+r) edge node[above]{} (HYi-N+r)
(HYi-N+r) edge node[above]{} (HZi-N+r)
(HXi) edge node[above]{} (HYi)
(HYi) edge node[above]{} (HZi)
(HXi+r) edge (P2)
(HZi-N) edge[out=355,in=175] node[above] {} (HXi-N+r)
(HZi-N+r) edge[out=355,in=175] node[above] {} (HXi)
(HZi) edge[out=355,in=175] node[above] {} (HXi+r);
}
\end{equation*}
So $H^i_{(r)}(Y)=0=H^{i-N+r}_{(N-r)}(Y)$, but this just means $Y \in \KK_N(\A)^\emptyset$.
Now consider $X \in \KK_N(\A)^\emptyset$ for the stability under the suspension functor. By definition of the stable category and its triangulated structure we have a triangle of the form 
\begin{equation*}
\tikz[heighttwo,xscale=2,yscale=2,baseline]{

\node (X) at (1,0){$X$};
\node (Y) at (2,0){$\PP_N(X)$};
\node (Z) at (3,0){$\Sigma(X)$};
\node (X1) at (4,0) {$\Sigma(X).$};

\draw[->]
(X) edge node[above]{$\lambda_X^\PP$} (Y)
(Y) edge  (Z)
(Z) edge  (X1);
}
\end{equation*}
We get again by Corollary~\ref{5.8} long exact sequences 
\begin{equation*}
\tikz[heighttwo,xscale=1,yscale=2,baseline]{

\node (P1) at (4,3){$\cdots$};
\node (HZi-N) at (8,3){$H^{i-N}_{(r)}(\Sigma(X))$};
\node (HXi-N+r) at (0,2) {$H^{i-N+r}_{(N-r)}(X)$};
\node (HYi-N+r) at (4,2) {$0$};
\node (HZi-N+r) at (8,2) {$H^{i-N+r}_{(N-r)}(\Sigma(X))$};
\node (HXi) at (0,1){$H^i_{(r)}(X)$};
\node (HYi) at (4,1){$0$};
\node (HZi) at (8,1){$H^i_{(r)}(\Sigma(X))$};
\node (HXi+r) at (0,0){$H^{i+r}_{(r)}(X)$};
\node (P2) at (4,0){$\cdots$};

\node[overlay] (Point) at (9.6,0){.};

\draw[overlay,->]
(P1) edge (HZi-N)
(HXi-N+r) edge  (HYi-N+r)
(HYi-N+r) edge (HZi-N+r)
(HXi) edge  (HYi)
(HYi) edge  (HZi)
(HXi+r) edge (P2)
(HZi-N) edge[out=355,in=175]  (HXi-N+r)
(HZi-N+r) edge[out=355,in=175] (HXi)
(HZi) edge[out=355,in=175] (HXi+r);
}
\end{equation*}
So we have for all $1 \le r \le N-1$ and $i \in \Z$ 
$$H^i_{(r)}(X)\cong H^{i-N+r}_{(N-r)}(\Sigma(X)).$$
 This proves the claim for $\Sigma(X)$. For $\Sigma^{-1}(X)$ consider the defining triangle of $\Sigma^{-1}(X)$:
\begin{equation*}
\tikz[heighttwo,xscale=2,yscale=2,baseline]{

\node (X) at (1,0){$\Sigma^{-1}(X)$};
\node (Y) at (2,0){$\II_N(X)$};
\node (Z) at (3,0){$X$};
\node (X1) at (4,0) {$X.$};

\draw[->]
(X) edge  (Y)
(Y) edge  node[above]{$\lambda_X^\II$}(Z)
(Z) edge  (X1);
}
\end{equation*}
This gives again long exact sequences:
\begin{equation*}
\tikz[heighttwo,xscale=1,yscale=2,baseline]{

\node (P1) at (4,3){$\cdots$};
\node (HZi-N) at (8,3){$H^{i-N}_{(r)}(X)$};
\node (HXi-N+r) at (0,2) {$H^{i-N+r}_{(N-r)}(\Sigma^{-1}(X))$};
\node (HYi-N+r) at (4,2) {$0$};
\node (HZi-N+r) at (8,2) {$H^{i-N+r}_{(N-r)}(X)$};
\node (HXi) at (0,1){$H^i_{(r)}(\Sigma^{-1}(X))$};
\node (HYi) at (4,1){$0$};
\node (HZi) at (8,1){$H^i_{(r)}(X)$};
\node (HXi+r) at (0,0){$H^{i+r}_{(r)}(\Sigma^{-1}(X))$};
\node (P2) at (4,0){$\cdots$};

\node[overlay] (Point) at (9,0){.};

\draw[overlay,->]
(P1) edge (HZi-N)
(HXi-N+r) edge  (HYi-N+r)
(HYi-N+r) edge (HZi-N+r)
(HXi) edge  (HYi)
(HYi) edge  (HZi)
(HXi+r) edge (P2)
(HZi-N) edge[out=355,in=175]  (HXi-N+r)
(HZi-N+r) edge[out=355,in=175] (HXi)
(HZi) edge[out=355,in=175] (HXi+r);
}
\end{equation*}
So we get for all $i \in \Z$ and $1\le r \le N-1$ that $$H^i_{(r)}(\Sigma^{-1}(X))\cong H^{i-N+r}_{(N-r)}(X).$$
which finished the claim. \qedhere
\end{description} 
\end{proof}
\end{proposition}

Now we will finally define the derived category in the way of \cite{N-complex}.

\begin{definition}
Let $\A$ be an abelian category and $N \in \N$. Then the \textbf{derived category of $N$-complexes in $\A$} is defined as the Verdier-quotient
\begin{align*}
\DD_N(\A):= \KK_N(\A)/\KK_N(\A)^\emptyset.
\end{align*}
\end{definition}

The following Lemma for $N$-complexes in a ring can be found in \cite{kapranov} as Proposition~1.5, and allows to calculate acyclicity of an $N$-complex using just one fixed $r$ with $1 \le r \le N-1$. Although this is a quite powerful statement the proof does not give a lot insight, so we omit it.

\begin{lemma}\label{kapranov formula}
Let $R$ be a ring and $N \in \N$. Then we have for $X\in \CC_N(\Mod_R)$ that $H_{(r)}^i=0$ for all $i \in \Z$ and $1 \le r \le N-1$ if and only if there is one $1 \le r \le N-1$ such that $H_{(r)}^i=0$ for all $i \in \Z$.
\begin{proof}
Proposition~1.5 in \cite{kapranov}.
\end{proof}
\end{lemma}

We will later need the following statement in order to control acyclicity on the level of morphisms.

\begin{proposition}\label{R compact generator}
Let $R$ be a ring and $N \in \N$. Then we have that $X \in K_N^\emptyset(\Mod_R)$ if and only if for all $i\in \Z$ we have $\KK_N(\Mod_R)(\mu_1^i(R),X)=0$.
\begin{proof}
By Lemma~\ref{kapranov formula} it suffices to show that $H_{(1)}^j(X)=0$ for all $j \in \Z$ if and only if $\KK_N(\Mod_R)(\mu_1^i(R),X)=0$ for all $i \in \Z$. 

Hence consider First $x \in \ker(d_X)^i\subset X^i$. We want to show that $x \in \im(d^{N-1})$ . 

This $x$ defines a map 
\begin{align*}
 \mu_1^i(R) &\to X^i \\ 1  &\mapsto x
 \end{align*}
 which by assumption has to factor as $\mu_1^i(R) \to \II_N(\mu_{1}^i(R)) \xrightarrow{\varphi} X$ in $\CC_N(\Mod_R)$. But we have $\II_N(\mu_{1}^i(R)) = \mu_{N}^{i}(R)$. So define $y \in X^{i-N+1}$ by 
 $$y:=\varphi^{i-N+1}(1).$$
  This now yields 
 $$d_X^{N-1}(y)=d_X^{N-1}\circ \varphi^{i-N+1} (1)=\varphi^i \circ d_{\mu_N^i(R)}^{N-1}(1)=\varphi^i(1)=x,$$
 so $d^{N-1}(y)=x$. Since this now implies that $x \in \im(d^{N-1})$ we have $$H_{(1)}^i(X)= \ker(d) / \im(d^{N-1})=0$$as desired.
 
 For the other implication let $\varphi:\mu_1^i(R) \to X$. This morphism is uniquely defined by $x:=\varphi^i(1)$. This is an element of $\ker(d_X)$, since $d \circ \varphi = \varphi \circ d =0$. The acyclicity of $X$ now means that there exists a $y \in X^{i-N+1}$ such that $d^{N-1}(y)=x$. Now consider the morphism $\psi: \mu_N^{i}(R) \to X$ defined by $\psi^{i-N+1}(1)=y$. This $\psi$ has the property $\psi^i(1)=x$ and hence $\varphi$ factors through $\mu_N^i(R)$. But this just means that $\varphi=0$ in $\KK_N(\Mod_R)$.   
\end{proof}
\end{proposition}

\subsection{The category of h-projectives}
$\;$

Similar to the case of ordinary chain complexes we will establish the notion of h-projective $N$-complexes. This allows to realize the derived category of $N$-complexes of a ring already as subcategory of the homotopy category. Furthermore it advantageously avoids roof diagrams for calculations of morphisms. However, due to our general setting, these two approaches might not coincide. Nonetheless we will show in the next subsection that this works out for $N$-complexes of modules over ring.

\begin{definition}
Let $\A$ be an abelian category and $N \in \N$. An object $X \in \CC_N(\A)$ is called \textbf{h-projective} if for all acyclic objects $Y$ in $\KK_N(\A)$ the abelian group $\KK_N(\A)(X,Y)$ is trivial. We refer to the full subcategory of h-projective objects in $\KK_N(\A)$ as $\hproj_{N}(\A)$.
\end{definition}

Since the derived category has naturally the structure of a triangulated category, this should hold as well for the category of h-projectives. 

\begin{proposition}\label{5.14}
Let $\A$ be an abelian category and $N \in \N$. Then the category $\hproj_{N}(\A)$ is a thick triangulated subcategory of $\KK_N(\A)$ which is closed under coproducts if $\A$ is closed under coproducts.
\begin{proof}
We again check the properties separately:
\begin{description}
\item[\normalfont triangulated]
First consider a triangle
\begin{equation*}
\tikz[heighttwo,xscale=2,yscale=2,baseline]{

\node (X) at (1,0){$X$};
\node (Y) at (2,0){$Y$};
\node (Z) at (3,0){$Z$};
\node (X1) at (4,0) {$\Sigma(X)$};

\draw[->]
(X) edge node[above]{$f$} (Y)
(Y) edge node[above]{$g$} (Z)
(Z) edge node[above]{$h$} (X1);
}
\end{equation*}
 in $\KK_N(\A)$ such that $X$ and $Z$ are h-projective. Then we get for all $T \in \KK_N(\A)$ a long exact sequence since $\KK_N(\A)$ is a triangulated category
\begin{equation*}
\tikz[heighttwo,xscale=4,yscale=2,baseline]{

\node (P1) at (2,2){$\cdots$};
\node (HZi-N+r) at (3,2) {$\KK_N(\A)(\Sigma X,T)$};
\node (HXi) at (1,1){$\KK_N(\A)(Z,T)$};
\node (HYi) at (2,1){$\KK_N(\A)(Y,T)$};
\node (HZi) at (3,1){$\KK_N(\A)(X,T)$};
\node (HXi+r) at (1,0){$\KK_N(\A)(\Sigma^{-1} Z,T)$};
\node (P2) at (2,0){$\cdots$};

\node[overlay] (Point) at (3.2,0){.};
    
\draw[overlay,->]
(P1) edge (HZi-N+r)
(HXi) edge node[above]{$\circ g$} (HYi)
(HYi) edge node[above]{$\circ f$} (HZi)
(HXi+r) edge (P2)
(HZi-N+r) edge[out=355,in=175] (HXi)
(HZi) edge[out=355,in=175]   (HXi+r);
}
\end{equation*}
Hence we get for all acyclic $T$ an exact sequence of the form
\begin{equation*}
\tikz[heighttwo,xscale=2,yscale=2,baseline]{

\node (X) at (1,0){$0$};
\node (Y) at (2,0){$\KK_N(\A)(Y,T)$};
\node (Z) at (3,0){$0$};

\draw[->]
(X) edge node[above]{} (Y)
(Y) edge node[above]{} (Z);
}
\end{equation*}
and so $Y$ has to be h-projective as well.

To show the stability under suspension, just observe that we have
  $$ \KK_N(\A)(\Sigma(X),T)\cong \KK_N(\A)(X,\Sigma^{-1} (T)) \; \text{ and }\KK_N(\A)(\Sigma^{-1}(X),T)\cong \KK_N(\A)(X,\Sigma (T))$$
   by the definition of a triangulated category for all $T \in \KK_N(\A)^\emptyset$.
   But this now shows the claim since $\KK_N(\A)^\emptyset$ is triangulated.

\item[\normalfont thick]Let $X$ be a direct summand of $Y$, where $Y$ is an h-projective $N$-complex. Then we have $X\oplus Z \cong Y$ for some $Z \in \KK_N(\A)$. Now since  $\KK_N(\A)(\_,T)$ is an additive functor, we get: 
\begin{align*}
\KK_N(\A)(X\oplus Z,T)=\KK_N(\A)(X,T)\oplus \KK_N(\A)(Z,T).
\end{align*}
But again specialization to $T$ acyclic yields 
\begin{align*}
\KK_N(\A)(X,T)\oplus \KK_N(\A)(Z,T)=\KK_N(\A)(X\oplus Z,T)=0
\end{align*}
and so $\KK_N(\A)(X,T)=0$.
\item[\normalfont coproducts] Assume that $\A$ admits coproducts. Then by Proposition~\ref{complete/cocomplete} and Lemma~\ref{naive colimits} the category $\KK_N(\A)$ admits coproducts. Now the claim follows since for $T \in \KK_N(\A)^\emptyset$:
 $$\KK_N(\A)(\bigoplus_{i \in I} X_i, T)= \prod_{i \in I} \KK_N(X_i,T)= \prod_{i \in I} 0 =0.$$
\end{description}
So all claimed properties hold.
\end{proof}
\end{proposition}

Furthermore the derived category has the property that any acyclic object in it is zero, which should hold as well for h-projectives.

\begin{proposition}\label{5.15}
Let $\A$ be an abelian category, $N \in \N$ and $X \in \hproj_{N}(\A)$. Then $X$ is acyclic if and only if it is isomorphic to zero in $\KK_N(\A)$.
\begin{proof}
If $X\cong 0$ in $\KK_N(\A)$, then $H_{(r)}^i(X)=H_{(r)}^i(0)=0$ and hence it is acyclic. For the converse observe that, since $X$ is both acyclic and h-projective, we have $\KK_N(\A)(X,X)=0$, now consider $\id_X \in \KK_N(\A)(X,X)=0$. So we get $\id_X=0$ and so $X \cong 0$ in $\KK_N(\A)$.
\end{proof}
\end{proposition}

Now finally the last property that makes the h-projectives resemble the derived category we can prove in general is that a morphism is an isomorphism in $\hproj_{N}(\A)$ if and only if it is an $N$-quasi-isomorphism.

\begin{proposition}\label{quasi iso between hproj}
Let $\A$ be an abelian category, $N \in \N$ and $f:X \to Y$ a morphism in $\hproj_N(\A)$. Then $f$ is an $N$-quasi-isomorphism if and only if $f$ is an isomorphism in $\KK_N(\A)$.
\begin{proof}$\;$\\
\begin{itemize}
\item[``$\Leftarrow$'':]This implication is just Corollary~\ref{5.6}.
\item[``$\Rightarrow$'':]Let $f:X \to Y$ be an $N$-quasi-isomorphism between h-projective objects. Then we get by Corollary~\ref{5.9} a triangle
\begin{equation*}
\tikz[heighttwo,xscale=2,yscale=2,baseline]{

\node (X) at (1,0){$X$};
\node (Y) at (2,0){$Y$};
\node (Z) at (3,0){$\Cone(f)$};
\node (X1) at (4,0) {$\Sigma(X)$};

\draw[->]
(X) edge node[above]{$f$} (Y)
(Y) edge node[above]{$g$} (Z)
(Z) edge node[above]{$h$} (X1);
}
\end{equation*}
in $\KK_N(\A)$ such that $H_{(r)}^i(\Cone(f))=0$ for all $i \in \Z$ and $1 \le r \le N-1$. Hence $\Cone(f)$ is acyclic. On the other hand $\Cone(f)$ has to be h-projective by Proposition~\ref{5.14}. But this means by Proposition~\ref{5.15} that $\Cone(f)\cong 0$ in $\KK_N(\A)$, and so $f$ is an isomorphism in $\KK_N(\A)$. \qedhere
\end{itemize}
\end{proof}
\end{proposition}

The following proposition will allow us to use certain h-projectives as building blocks for our semifree resolutions.

\begin{proposition}\label{levelwise free h-projective}
Let $R$ be a ring and $N \in \N$. Then every $N$-complex in $\Mod_R$ which is levelwise free with vanishing differential is $h$-projective.
\begin{proof}
Let $X$ be such an $N$-complex and $T$ an acyclic $N$-complex. Then we have 
\begin{align*}
\KK_N(\Mod_R)(X,T)&=\prod_{i\in \Z}\KK_N(\Mod_R)(\mu_1^i(R^{\oplus I_i}),T))\\
&=\prod_{i\in \Z}\prod_{I_i} \KK_N(\Mod_R)(\mu_1^i(R),T)=\prod_{i\in \Z}\prod_{I_i}0=0.
\end{align*}
Here the second to last equation holds by Proposition~\ref{R compact generator}.
\end{proof}
\end{proposition}

Although we do not have a natural equivalence $\hproj_N(\A)\cong \DD_N(\A)$ in general, we still get a fully faithful functor in general. To show this we will need the following general property of triangulated categories.

\begin{proposition}\label{splitting triangle}
Let $\T$ be a triangulated category and 
\begin{equation*}
\tikz[heighttwo,xscale=2,yscale=2,baseline]{

\node (X) at (1,0){$X$};
\node (Y) at (2,0){$Y$};
\node (Z) at (3,0){$Z$};
\node (X1) at (4,0) {$\Sigma(X)$};

\draw[->]
(X) edge node[above]{$f$} (Y)
(Y) edge node[above]{$g$} (Z)
(Z) edge node[above]{$h$} (X1);
}
\end{equation*}
a triangle in $\T$. Then $h=0$ if and only if $g$ admits a splitting $s:Z\to Y$.
\begin{proof}
First let $h=0$. Then the following diagram commutes and hence the dashed arrow exists:
\begin{equation*}
\tikz[heighttwo,xscale=2,yscale=2,baseline]{

\node (X1) at (1,0){$X$};
\node (Y1) at (2,0){$Y$};
\node (Z1) at (3,0){$Z$};
\node (SX1) at (4,0) {$\Sigma(X)$.};
\node (X2) at (1,1){$X$};
\node (Y2) at (2,1){$X \oplus Z$};
\node (Z2) at (3,1){$Z$};
\node (SX2) at (4,1) {$\Sigma(X)$};

\draw[->]
(X1) edge node[above]{$f$} (Y1)
(Y1) edge node[above]{$g$} (Z1)
(Z1) edge node[above]{$0$} (SX1)
(X2) edge node[above]{$\begin{bmatrix}
\id \\
0
\end{bmatrix}$} (Y2)
(Y2) edge node[above]{$\begin{bmatrix}
0 & \id
\end{bmatrix}$} (Z2)
(Z2) edge node[above]{$0$} (SX2)
(X2) edge node[right]{$\id$} (X1)
(Z2) edge node[right]{$\id$} (Z1)
(SX2) edge node[right]{$\id$} (SX1);
\draw[dashed, ->]
(Y2) edge node[right]{$\widehat{s}$}(Y1);
}
\end{equation*}
The commutativity of the middle square now yields 
$$g \circ \widehat{s} \circ \begin{bmatrix}
0\\
\id_Z
\end{bmatrix}=\id_Z.$$
So we get the split $s:=\widehat{s} \circ \begin{bmatrix}
0\\
\id_Z
\end{bmatrix} : Z \to Y $.

For the other direction just observe that if he have a split $s:Z \to Y$ the diagram
\begin{equation*}
\tikz[heighttwo,xscale=2,yscale=2,baseline]{

\node (X1) at (1,0){$X$};
\node (Y1) at (2,0){$Y$};
\node (Z1) at (3,0){$Z$};
\node (SX1) at (4,0) {$\Sigma(X)$.};
\node (X2) at (1,1){$X$};
\node (Y2) at (2,1){$X \oplus Z$};
\node (Z2) at (3,1){$Z$};
\node (SX2) at (4,1) {$\Sigma(X)$,};

\draw[->]
(X1) edge node[above]{$f$} (Y1)
(Y1) edge node[above]{$g$} (Z1)
(Z1) edge node[above]{$h$} (SX1)
(X2) edge node[above]{$\begin{bmatrix}
\id \\
0
\end{bmatrix}$} (Y2)
(Y2) edge node[above]{$\begin{bmatrix}
0 & \id
\end{bmatrix}$} (Z2)
(Z2) edge node[above]{$0$} (SX2)
(X2) edge node[right]{$\id$} (X1)
(Z2) edge node[right]{$\id$} (Z1)
(SX2) edge node[right]{$\id$} (SX1)
(Y2) edge node[right]{$\begin{bmatrix}
f & s
\end{bmatrix}$}(Y1)
}
\end{equation*}
commutes. So $$h=h \circ \id = \id \circ 0 = 0$$ 
as claimed. \qedhere
\end{proof}
\end{proposition}

\begin{remark}
We will use in the following proof that morphisms in the derived category, or any Verdier quotient can be characterized as roof diagrams or coroof diagrams. We refer the reader for the notions of these to the book by Neeman \cite{neeman}. In particular we refer to Definition~2.1.11, Lemma~2.1.14 and Remark~2.1.25.
\end{remark}

\begin{theorem}\label{restriction fully faithful}
The functor $\pi:\KK_N(\A) \to \DD_N(\A)$  induces a bijection
$$\KK_N(\A)(X,Y) \xrightarrow{\pi} \DD_N({\A})(X,Y) \text{ for } X \in \hproj_N(\A).$$
 In particular the restriction of $\pi$ to $\hproj_N(\A)$ is fully faithful.
\begin{proof}
Let $f \in \DD_N(\A)(X,Y)$, where $X \in \hproj_N(\A)$. Then we have that $f$ is described by a roof diagram of the form $(s,g )$ for $s:W\to X$ an $N$-quasi-isomorphism and $g:W \to Y$ by the construction of $\DD_N(\A)$ as Verdier quotient. Now consider the diagram
\begin{equation*}
\tikz[heighttwo,xscale=1,yscale=1,baseline]{

\node (X) at (0,1){$X$};
\node (Y) at (2,1){$Y$};
\node (C) at (-1,0){$C(s)$};
\node (SC) at (2,3) {$\Sigma^{-1}(C(s))$};
\node (W) at (1,2){$W$};

\draw[->]
(W) edge node[above left]{$s$} (X)
(X) edge node[above left]{$\varphi$} (C)
(SC) edge node[right] {$\psi$}(W)
(W) edge node[above right]{$g$} (Y);
}
\end{equation*}
in $\KK_N(\A)$, where $C(s)$ is the cone of $s$ in $\KK_N(\A)$. Since $s$ is an $N$-quasi-isomorphism we get by Corollary~\ref{5.9} that $C(s)$ is acyclic, but this yields $\varphi = 0$ since $X \in \hproj_N(\A)$. Using that $\KK_N(\A)$ is triangulated and Proposition~\ref{splitting triangle} we get a splitting $t: X \to W$ of $s$ and the following diagram commutes
\begin{equation*}
\tikz[heighttwo,xscale=1,yscale=1,baseline]{

\node (X) at (0,0.5){$X$};
\node (X2) at (0,2.5){$X$};
\node (Y) at (2,0.5){$Y$};
\node (W) at (1,1.5){$W$};

\draw[->]
(W) edge node[below right]{$s$} (X)
(W) edge node[above right]{$g$} (Y)
(X2) edge node[above right]{$t$}(W)
(X2) edge node[left]{$\id$}(X);

}.
\end{equation*}
Hence $f$ can be represented by the equivalent roof $(\id ,g \circ t ) $, so $f=\pi(g \circ t)$ which proves that  $\pi$ induces a surjection. 

To show that it is also injective consider a morphism $g:X \to Y$ such that $\pi(g)=0$. The morphism $\pi (g)$ is now described by the coroof diagram $(g,\id_Y)$  which is by assumption equivalent to $(0,\id_Y)$ hence we have an $N$-complex $Z$ and an $N$-quasi-isomorphisms $\varphi:Y \to Z$ such that the following diagram commutes in $\KK_N(\A)$:
\begin{equation*}
\tikz[heighttwo,xscale=2,yscale=2,baseline]{

\node (X) at (0,0){$X$};
\node (Y1) at (1,1){$Y$};
\node (Y2) at (2,0){$Y$};
\node (Y3) at (1,-1){$Y.$};
\node (Z) at (1,0){$Z$};
  
\draw[->]
(X) edge node[above left]{$g$} (Y1)
(Y2) edge node[above right]{$\id_Y$} (Y1)
(X) edge node[above right]{$0$}(Y3)
(Y2) edge node[below right]{$\id_Y$}(Y3)
(X) edge node[above]{$0$} (Z)
(Y1)edge node[right]{$\varphi$}(Z)
(Y2)edge node[above]{$\varphi$}(Z)
(Y3)edge node[right]{$\varphi$}(Z);
}
\end{equation*}
In particular we get $\varphi \circ g =0 $ in $\KK_N(\A)$. Now since $\varphi$ is an $N$-quasi-isomorphism we get the following triangle in $\KK_N(\A)$
\begin{equation*}
\tikz[heighttwo,xscale=2,yscale=2,baseline]{

\node (X1) at (1,0){$Y$};
\node (X2) at (2,0){$Z$};
\node (X3) at (3.5,0){$\Cone_{\KK_N(\A)}(\varphi)$};
\node (X4) at (5,0){$\Sigma Y$};

\draw[->]

(X1) edge node[above] {$\varphi$}(X2)
(X2) edge (X3)
(X3) edge  (X4);

}
\end{equation*}
where $\Cone_{\KK_N(\A)}(\varphi)$ is acyclic. Since $\KK_N(\A)$ is triangulated this induces a long exact sequence of the form
\begin{equation*}
\tikz[heighttwo,xscale=3,yscale=2,baseline]{
\node (P1) at (-1,0){$\cdots$};
\node (X0) at (0,0){$\KK_N(X,\Sigma^{-1}C)$};
\node (X1) at (1,0){$\KK_N(X,Y)$};
\node (X2) at (2,0){$\KK_N(X,Z)$};
\node (X3) at (3,0){$\KK_N(X,C)$};
\node (P2) at (4,0){$\cdots$};
    
\draw[->]
(P1) edge (X0)
(X0) edge (X1)
(X1) edge node[above] {$\varphi\circ $}(X2)
(X2) edge (X3)
(X3) edge  (P2);

}
\end{equation*}
where we denote by $C$ the cone of $\varphi$ in $\KK_N(\A)$. Since $X$ is h-projective we get $$\KK_N(X,C)\cong 0 \cong \KK_N(X,\Sigma^{-1}C),$$ as $C$ is acyclic by  Proposition~\ref{5.1} and $\KK_N(\A)^{\emptyset}$ is triangulated by Proposition~\ref{acyclics thick triangulated}. But this means that $\varphi\circ :\KK_N(X,Y) \to \KK_N(X,Z)$ is an isomorphism. So since $\varphi \circ g =0$ we get that $g$ already has to vanish in $\KK_N(X,Y)$. This shows that $\pi$ is injective on $\KK_N(\A)(X,Y)$.

To conclude the in particular part just observe that $\pi$ induces by the above proof a bijection $\KK_N(\A)(X,Y)\to \DD_N(\A)(X,Y)$ for all $X,Y \in \hproj_N(\A)$.
\end{proof}
\end{theorem}

\subsection{Semifree resolutions}
$\;$

To conclude the section we will show that the canonical functor $\KK_N(\A) \to \DD_N(\A)$ restricts to an equivalence $\hproj_N(\A)\to \DD_N(\A)$ if we have $\A= \Mod_R$ for a ring $R$. By Theorem~\ref{restriction fully faithful} we only have to show that this functor is essentially surjective. To do this we will give here a construction following the idea presented in the remark after Lemma~13.3 in ``Dg-quotients of Dg-categories'' by Drinfeld \cite{drinfeld}. 

The recipe for this construction is to consider a certain type of objects, namely semifree modules, and show that every complex admits a semifree replacement in $\DD_N(\A)$. Since these are h-projective this will give essential surjectivity of the projection functor.

So let us first define what it means for an $N$-complex to be semifree and why these $N$-complexes are h-projective.

\begin{definition}
Let $R$ be a ring and $N \in \N$. We call an $N$-complex $F$ in $\Mod_R$ \textbf{semifree} if there exists a filtration $0=F_0\subset F_1 \subset F_2 \subset F_3 \subset...$, such that:
\begin{enumerate}
\item $F=\bigcup_{n\in \N} F_n,$
\item $F_i / F_{i-1}$ is levelwise free with vanishing differential $\forall i \geq 1.$
\end{enumerate}
\end{definition}

\begin{proposition}\label{semifree => h-projective}
Let $R$ be a ring, $N \in \N$ and $M$ a semifree $N$-complex in $R$. Then we have $M\in \hproj_N(\Mod_R)$.
\begin{proof}
Since $M$ is semifree we have a filtration $0=F_0\subset F_1 \subset F_2 \subset F_3 ...$ such that $M= \colim_i F_i$. Now we have by Proposition~\ref{levelwise free h-projective} that $F_0$ and $F_i/F_{i-1}$ are h-projective and levelwise free. Assume that $F_{i-1}$ is h-projective, then we have a short exact sequence of the form
\begin{equation*}
\tikz[heighttwo,xscale=2,yscale=2,baseline]{

\node (01) at (0,0){$0$};
\node (X1) at (1,0){$F_{i-1}$};
\node (X2) at (2,0){$F_{i}$};
\node (X3) at (3,0){$F_{i}/F_{i-1}$};
\node (02) at (4,0){$0$};
    
\draw[->]

(01) edge (X1)
(X3) edge (02)
(X1) edge (X2)
(X2) edge (X3);

}
\end{equation*}
which levelwise splits, since $F_i/F_{i-1}$ is levelwise free. So the sequence induces a triangle 
\begin{equation*}
\tikz[heighttwo,xscale=2,yscale=2,baseline]{

\node (X1) at (1,0){$F_{i-1}$};
\node (X2) at (2,0){$F_{i}$};
\node (X3) at (3,0){$F_{i}/F_{i-1}$};
\node (X4) at (4,0){$\Sigma F_{i-1}$};

\draw[->]

(X1) edge (X2)
(X2) edge (X3)
(X3) edge  (X4);

}
\end{equation*}
in $\KK_N(\A)$ by Corollary~\ref{4.9}. Now since $\hproj_N(\Mod_R)$ is a triangulated subcategory and both $F_{i-1}$ and $F_{i}/F_{i-1}$ are h-projective, we get that $F_i$ is h-projective. So inductively all $F_j$ are h-projective and all inclusions $\varphi_i: F_i \hookrightarrow F_{i+1}$ levelwise split.
To show now that finally $M$ is h-projective consider the short exact sequence 
\begin{equation*}
\tikz[heighttwo,xscale=2,yscale=2,baseline]{

\node (01) at (0,0){$0$};
\node (X1) at (1,0){$\bigoplus\limits_{i=0}^{\infty} F_i$};
\node (X2) at (2,0){$\bigoplus\limits_{i=0}^{\infty} F_i$};
\node (X3) at (3,0){$M$};
\node (02) at (4,0){$0$};
    
\draw[->]

(01) edge (X1)
(X3) edge (02)
(X1) edge (X2);
\draw[->]
(X1) edge node[above] {$\Psi$}(X2);
\draw[->>]
(X2) edge (X3);
}
\end{equation*}
exhibiting $M$ as $\colim_i F_i$. Here $\Psi$ is the morphism induced by the composition $$F_i \xrightarrow{\begin{bmatrix}
\id_{F_i}\\ -\varphi_i
\end{bmatrix}} F_i \oplus F_{i+1} \xhookrightarrow{\begin{bmatrix}
\iota_i & \iota_{i+1}
\end{bmatrix}} \bigoplus\limits_{i=0}^{\infty} F_i , $$ where $\varphi_i$ is the inclusion $F_i \hookrightarrow \F_{i+1}$ and $\iota_j$ is the inclusion $F_j \hookrightarrow \bigoplus\limits_{i=0}^{\infty} F_i$. Now observe that $\Psi$ is of the form 
$$\begin{bmatrix}
\id_{F_0} & 0 &0&0 & \cdots \\
-\varphi_0 & \id_{F_1} &0&0&\cdots \\
0 & -\varphi_1 & \id_{F_2} & 0& \cdots \\
0&0&-\varphi_2 & \id_{F_3} & \cdots \\
\vdots & \vdots & \vdots & \vdots & \ddots 
\end{bmatrix}.$$
Consider the following morphism in $\Gr(\Mod_R)$, the category of graded modules over $R$:
$$ \begin{bmatrix}
0 & -\psi_0 & -\psi_0 \circ \psi_1 & -\psi_0 \circ \psi_1 \circ \psi_2 & \cdots \\
0 & 0 & \psi_1 & \psi_1 \circ \psi_2 & \cdots \\
0 &0 & 0 & \psi_2 & \cdots \\
0&0&0 & 0 & \cdots \\
\vdots & \vdots & \vdots & \vdots & \ddots 
\end{bmatrix},$$
where $\psi_i$ denotes a split of $\varphi_i$ in $\Gr(\Mod_R)$, which exists by the argument above. This is a levelwise split of $\Psi$, since  one can calculate
\begin{align*}
\begin{bmatrix}
0 & -\psi_0 & -\psi_0 \circ \psi_1 & -\psi_0 \circ \psi_1 \circ \psi_2 & \cdots \\
0 & 0 & -\psi_1 & -\psi_1 \circ \psi_2 & \cdots \\
0 &0 & 0 & -\psi_2 & \cdots \\
0&0&0 & 0 & \cdots \\
\vdots & \vdots & \vdots & \vdots & \ddots 
\end{bmatrix}
\begin{bmatrix}
\id_{F_0} & 0 &0&0 & \cdots \\
-\varphi_0 & \id_{F_1} &0&0&\cdots \\
0 & -\varphi_1 & \id_{F_2} & 0& \cdots \\
0&0&-\varphi_2 & \id_{F_3} & \cdots \\
\vdots & \vdots & \vdots & \vdots & \ddots 
\end{bmatrix}= \begin{bmatrix}
\id_{F_0} & 0 & 0 &0  & \cdots \\
0& \id_{F_1} &0&0&\cdots \\
0 & 0& \id_{F_2} & 0& \cdots \\
0&0&0 & \id_{F_3} & \cdots \\
\vdots & \vdots & \vdots & \vdots & \ddots 
\end{bmatrix}
\end{align*}
 where the zeroes above the diagonal appear since the entries there are $$-\psi_{i}\circ ... \circ \psi_{j} + \psi_{i}\circ ... \circ \psi_{j}\circ \psi_{j+1}\circ \varphi_{j+1}=-\psi_{i}\circ ... \circ \psi_{j} + \psi_{i}\circ ... \circ \psi_{j}=0.$$ Hence the sequence levelwise splits. This gives again that the sequence induces a triangle
\begin{equation*}
\tikz[heighttwo,xscale=2,yscale=2,baseline]{

\node (X1) at (1,0.3){$\bigoplus\limits_{i=0}^{\infty} F_i$};
\node (X2) at (2,0.3){$\bigoplus\limits_{i=0}^{\infty} F_i$};
\node (X3) at (3,0.3){$M$};
\node (X4) at (4,0.3){$\Sigma \bigoplus\limits_{i=0}^{\infty} F_i$};

\draw[->]

(X1) edge (X2)
(X2) edge (X3)
(X3) edge  (X4);

}
\end{equation*}
in $\KK_N(\A)$. Now since $\hproj_N(\Mod_R)$ is by Proposition~\ref{5.14} closed under coproducts we have that  $\bigoplus\limits_{i=0}^{\infty} F_i $ is h-projective and again since $\hproj_N (\Mod_R)$ is triangulated we get that $M$ is h-projective as claimed. 
\end{proof}
\end{proposition}

Having established this, let us construct our semifree resolutions:

\begin{theorem}\label{semifree resolution}
Let $R$ be a ring, $N \in \N$ and $M$ an $N$-complex in $\Mod_R$. Then there exists a semifree $N$-complex $F$ and a morphism $f: F \to M$ such that $f$ is an $N$-quasi-isomorphism. Furthermore, $f$ can be chosen surjective and such that the induced maps $\ker(d^r_F) \to \ker(d^r_M)$ are surjective for all $1 \le r \le N-1$.
\begin{proof}
Let the $N$-complex $F_1$ be defined by $F_1^j := \langle \ker(d_M\mid_{M^j})\rangle_R$ and by setting $d_{F_1}=0$, let $f_1:F_1 \to M$ be the morphism defined levelwise by the inclusions $\iota: \ker(d_M\mid_{M^j}) \hookrightarrow M^j$. This defines a morphism of $N$-complexes since $f_1 \circ d_{F_1} = 0 = d_M \circ f_1$.

Now assume $F_{i-1}$ and $f_{i-1}$ are already constructed and let  $B_{i-1}$ be the graded set defined levelwise by
\begin{equation*}
B_{i-1}^{j}:= \{(b,y_b)\in \ker(d^{N-1}_{F_{i-1}}\mid_{F_{i-1}^{j+1}})\times M^j \mid f_{i-1}(b)=d(y_b)\}.
\end{equation*}
 Then define  
\begin{align*}
F_i^j &:= F_{i-1}^j \oplus \langle B_{i-1}^j\rangle_R,\\
d_{F_i} &:= \begin{bmatrix} d_{F_{i-1}} & d_{B_{i-1}} \\
0 & 0 \end{bmatrix}
\end{align*}
where $d_{B_{i-1}}((b,y_b)):=  b$.
This makes $F_i$ into an $N$-complex, since $d_{F_{i-1}}^N=0$ by assumption and
\begin{align*}
d^N_{F_i}((b,y_b)) &=d_{F_{i-1}}^{N-1}(b)=0,
\end{align*}
since by definition $b \in \ker(d^{N-1}_{F_{i-1}}\mid_{F_{i-1}^{j+1}})$.
Furthermore define 
\begin{align*}
f_i:F_i & \to M ,\\
(x , r(b,y_b)) & \mapsto f_{i-1}(x)+ r y_b.
\end{align*}
This defines a morphism of $N$-complexes since:
\begin{align*}
d\circ f_i(x , r(b,y_b)) =&d(f_{i-1}(x)+ r y_b) \\
=&f_{i-1}(d(x))+r d(y_b) \\
=&f_{i-1}(d(x)) +r f_{i-1}(b) \\
=&f_{i-1}(d(x)+r b) \\
=&f_{i-1}\circ d_{F_i} (x,r(b,y_b) ). 
\end{align*}

Now define 
\begin{equation*}
F:=\lim_{\longrightarrow}(0=F_0 \hookrightarrow F_1 \hookrightarrow F_2 \hookrightarrow F_3 \hookrightarrow ...).
\end{equation*} 
We get a natural map by the definition of the colimit which will be the desired $f:F \to M$.

\begin{equation*}
\tikz[heightone,xscale=3,yscale=3,baseline]{

\node (F0) at (0,1){$F_0=0$};
\node (F1) at (1,1){$F_1$};
\node (F2) at (2,1){$F_2$};
\node (F3) at (3,1){$F_3$};
\node (F4) at (4,1){$F_4$};

\node (P1) at (4.5,1){$\cdot\cdot\cdot$};

\node (P2) at (5,0.5) {\Large{$\rightsquigarrow$}};

\node (M1) at (0,0){$M$};
\node (F) at (6,1){$F$};
\node (M2) at (6,0){$M$};

\draw[->]
(F0) edge node[left] {\tiny{$f_0=0$}} (M1)
(F1) edge node[above left] {\tiny{$f_1$}} (M1)
(F2) edge node[above ] {\tiny{$f_2$}} (M1)
(F3) edge node[above ] {\tiny{$f_3$}} (M1)
(F4) edge node[below right] {\tiny{$f_4$}} (M1);

\draw[right hook ->]
(F0) edge node[left] {} (F1)
(F1) edge node[left] {} (F2)
(F2) edge node[left] {} (F3)
(F3) edge node[left] {} (F4);

\draw[, ->]
(F) edge node[left] {\tiny{$f$}} (M2);
}
\end{equation*}

The just defined $F$ is semifree since it is defined as the colimit of inclusions,  so we have a filtration $0=F_0\subset F_1 \subset F_2\subset F_3 \subset F_4...$ and $F=\bigcup_{n\in \N} F_i$. Furthermore we have that 
\begin{equation*}
(F_i/F_{i-1})^j= \langle B_{i-1}^j\rangle_R,
\end{equation*}
and hence $F_i/F_{i-1}$ is levelwise free, and the differential vanishes since by definition $d_{F_i}(x)\in F_{i-1}$. So we have that $F$ is semifree.

Now it just remains to show that $f:F\to M$ is a surjective $N$-quasi-isomorphism which is surjective on $\ker(d^r_M)$ for all $1 \le r \le N-1$:

Consider $f|_{\ker(d^r_F)}:\ker(d^r_F) \to \ker(d^r_M)$ for $1\le r \le N$:\\
Let $y \in \ker(d_M^{r})$ and define $y_i:=d^{r-i}(y)\in \ker(d^i)$, in particular $y_1\in \ker(d)$, hence there exists a $x_1:=1 y_1 \in \ker(d_{F_1})$ such that $f(x_1)=y_1$. So we get $(x_1,y_2) \in B_1$. By construction this implies:
\begin{equation*}
f_2(0, 1 (x_1,y_2))=y_2.
\end{equation*}
So let $x_2:= (0,  1 (x_1,y_2))\in \ker(d^2_{F_2})$.

Assume now that $x_i\in \ker(d^i_{F_i})$ for some $1\le i \le r-1$ has been constructed such that $f_i(x_i)=y_i$, then $(x_i,y_{i+1})\in B_i$ and one gets

\begin{equation*}
f_{i+1}(0,  1 (x_i,y_{i+1}))=y_{i+1}.
\end{equation*}
So define $x_{i+1}:=(0, 1 (x_i,y_{i+1}))\in \ker(d^{i+1}_{F_i+1})$. 

Now one gets inductively $x_{r} \in \ker(d^r_{F_{r}})$ with $f_{r}(x_{r})=y_{r}=y$ and hence the induced map $f_{r}:\ker(d^r_{F_r})\to \ker(d^r_M)$ is surjective. This shows that $f|_{\ker(d^r_F)}:\ker(d^r_F) \to \ker(d^r_M)$ is surjective. In particular $f$ is surjective since $\ker(d_M^N)=M$.

To prove that $f$ is an $N$-quasi-isomorphism it remains to show that $f$ is injective on cohomology.

Let $x\in \ker(d^r_F)$ with $y:=f(x)\in \im(d^{N-r}_M)\subset \ker(d_M^r)$. Then there exists an $F_j$ such that $x \in F_j$.
Now, since by assumption $f(x)\in \im(d^{N-r}_M)$, there exists a $y_{N-r}\in M$ such that $d^{N-r}(y_{N-r})=y$ and so $(x,d^{N-r-1}y_{N-r}) \in B_j$ and 
\begin{equation*}
f_{j+1}(0\oplus  1(x,d^{N-r-1}y_{N-r}))=d^{N-r-1}y_{N-r}.
\end{equation*}
Define $x_1:=0\oplus  1(x,d^{N-r-1}y_{N-r})$.

Assume now that $x_i$ has been constructed for $1 \le i \le N-r-1$ such that $f(x_i)=d^{N-r-i}(y_{N-r})$. Then $(x_i,d^{N-r-i-1}y_{N-r})\in B_{j-i}$ and 
\begin{equation*}
f_{j+i+1}(0\oplus  1(x_i,d^{N-r-(i+1)}y_{N-r}))=d^{N-r-(i+1)}y_{N-r}.
\end{equation*}
Hence define $x_{i+1}:=0\oplus 1(x_i,d^{N-r-(i+1)}y_{N-r})$.

So one gets inductively $x_{N-r}$. Now we compute by the definition of $F$: 
\begin{align*}
d_F^{N-r}(x_{N-r})=&d_F^{N-r}(x_{N-r-1}, y_{N-r})\\
=&d_F^{N-r-1}(x_{N-r-1})=d_F^{N-r-1}(x_{N-r-2},d_M y_{N-r})\\
=&d_F^{N-r-2}(x_{N-r-2})=d_F^{N-r-2}(x_{N-r-3},d_M^2 y_{N-r})\\
&\vdots\\
=&d_F(x_{1})=d_F(x,d_M^{N-r} y_{N-r})=x.
\end{align*}
Hence $x \in \im(d^{N-r}_F)$. This now shows that $f$ is injective in cohomology which concludes the proof.
\end{proof}
\end{theorem}

Finally using this we can show that $\hproj_N(\Mod_R)$ indeed is a model for the derived category of $N$-complexes over a ring $R$.

\begin{corollary}\label{equiv derived hproj}
Let $R$ be a ring and $N \in \N$. Then we have that the canonical functor $\pi : \KK_N(\Mod_R) \to \DD_N(\Mod_R)$ restricts to an equivalence of categories:
\begin{equation*}
\hproj_N(\Mod_R) \xrightarrow{\sim} \DD_N(\Mod_R).
\end{equation*}
\begin{proof}
We know by Theorem~\ref{restriction fully faithful} that $\pi\mid_{\hproj_N(\Mod_R)}$ is fully faithful, hence it just remains to show that it is essentially surjective. But by Theorem~\ref{semifree resolution} we have that every $X \in \KK_N(\Mod_R)$ is $N$-quasi-isomorphic to a semifree $N$-complex $F_X$, which is by Proposition~\ref{semifree => h-projective} h-projective. So $\pi(X)\cong\pi(F_X)$, which implies essential surjectivity.
\end{proof}
\end{corollary}

%% file: model_structure.tex
\section{A model structure on $\CC_N(\Mod_R)$}\label{section 8}

We will establish a model structure on the category of $N$-complexes over a ring $R$, which will give as model theoretic homotopy category the derived category of $N$-complexes over $R$. In particular this shows that the derived category of $N$-complexes over a ring is indeed the Gabriel-Zismann localisation of $\CC_N(\Mod_R)$ along the $N$-quasi-isomorphisms. To show this, we will always consider a ring and the category of $N$-complexes over it. So we fix throughout this section a ring $R$ and  a natural number $N \in \N$. 

The same model structure was constructed by J. Gillespie and M. Hovey in \cite{Gillespie2}. However our construction is independent and more hands on, since J. Gillespie and M. Hovey used the theory of graded n-Gorenstein K-algebras whereas our proof is direct.  Furthermore it is no surprise that both constructions yield this model structure, since it is the obvious generalization of the projective model structure on the category of chain complexes.

The first thing we will need are the following collections of morphisms, since we will use the small object argument to cofibrantly generate a model structure. To do so we will use Theorem~2.1.19 in \cite{hovey}.
\begin{definition}
Define the following classes of morphisms in $\CC_N(\Mod_R)$:
\begin{enumerate}
\item $J:=\left\{\iota_{i}:0 \to \mu_N^i(R) \middle| i \in \Z\right\}$.
\item $I:=J \cup\left\{\iota_{i,r}:\mu_{r}^{i}(R) \to \mu_N^i(R) \middle| i \in \Z\ 0< r < N\right\} $, where $\iota_{i,r}^k=\id_R$ whenever possible and zero otherwise.
\item $W:=\{\text{N-quasi-isomorphisms}\}\subset \Mor(\CC_N(\Mod_R))$.
\end{enumerate}
\end{definition}

The set $I$ will become the set of generating cofibrations and $J$ the set of generating trivial cofibrations. However in order to classify the induced fibrations, trivial fibrations and trivial cofibrations we will need the notion of left and right lifting property respectively.

\begin{definition}
Let $\C$ be a category and $\Psi \subset \Mor(\C)$ a class of morphisms. Then we say that:
\begin{enumerate}
\item A morphism $\varphi: X \to Y$ has the \textbf{left lifting property with respect to $\Psi$} if for all $(\psi:M \to N) \in \Psi$ and all commutative squares of the form
\begin{equation*}
\tikz[heighttwo,xscale=3,yscale=2,baseline]{
\node (N) at (1,0) {$N$};
\node (M) at (1,1) {$M$};
\node (Y) at (0,0) {$Y$};
\node (X) at (0,1) {$X$};
    
\draw[->]
(X) edge node[left] {$\varphi$}(Y)
(X)edge(M)
(Y)edge(N)
(M)edge node[right] {$\psi$} (N);
}
\end{equation*}
there exists a lift $h:Y \to M$ such that the following diagram commutes:
\begin{equation*}
\tikz[heighttwo,xscale=3,yscale=2,baseline]{
\node (N) at (1,0) {$N.$};
\node (M) at (1,1) {$M$};
\node (Y) at (0,0) {$Y$};
\node (X) at (0,1) {$X$};
    
\draw[->]
(X) edge node[left] {$\varphi$}(Y)
(Y) edge node[above left] {$h$}(M)
(X)edge(M)
(Y)edge(N)
(M)edge node[right] {$\psi$} (N);
}
\end{equation*}
\item A morphism $\varphi: X \to Y$ has the \textbf{right lifting property with respect to $\Psi$} if for all $(\psi:M \to N) \in \Psi$ and all commutative squares of the form
\begin{equation*}
\tikz[heighttwo,xscale=3,yscale=2,baseline]{
\node (N) at (0,0) {$N$};
\node (M) at (0,1) {$M$};
\node (Y) at (1,0) {$Y$};
\node (X) at (1,1) {$X$};
    
\draw[->]
(X) edge node[left] {$\varphi$}(Y)
(M)edge(X)
(N)edge(Y)
(M)edge node[right] {$\psi$} (N);
}
\end{equation*}
there exists a lift $h:N \to X$ such that the following diagram commutes:
\begin{equation*}
\tikz[heighttwo,xscale=3,yscale=2,baseline]{
\node (N) at (0,0) {$N$};
\node (M) at (0,1) {$M$};
\node (Y) at (1,0) {$Y.$};
\node (X) at (1,1) {$X$};
    
\draw[->]
(X) edge node[left] {$\varphi$}(Y)
(N) edge node[above left] {$h$}(X)
(M)edge(X)
(N)edge(Y)
(M)edge node[right] {$\psi$} (N);
}
\end{equation*}
\item We denote the collections of morphisms with the left and right lifting property with respect to $\Psi$ by $\boxslash \Psi$ and $\Psi \boxslash$ respectively.
\end{enumerate}
\end{definition}

\begin{remark}
The morphisms $\iota_{i,N-1}$ are of the form:
\begin{equation*}
\tikz[heighttwo,xscale=2,yscale=2,baseline]{
\node (LX) at (-1,1){$\cdots$};
\node (X1) at (0,1){$0$};
\node(X2) at (1,1) {$0$};
\node (X3) at (2,1){$R$};
\node (P1) at (3,1){$\cdots$};
\node (XN-1) at (4,1) {$R$};
\node(XN) at (5,1) {$R$};
\node(XN+1) at (6,1) {$0$};
\node (RX) at (7,1){$\cdots$};
 
\node (LY) at (-1,0){$\cdots$};
\node (Y1) at (0,0){$0$};
\node (Y2) at (1,0){$R$};
\node (Y3) at (2,0){$R$};
\node (P2) at (3,0){$\cdots$};
\node (YN-1) at (4,0) {$R$};
\node (YN) at (5,0) {$R$};
\node (YN+1) at (6,0) {$0$};
\node (RY) at (7,0){$\cdots$};

\node[overlay] (Point) at (7.3,0){.};

\node (P3) at (3,0.5) {$\cdots$};   
    
\draw[->]
(LX) edge node[above] {} (X1)
(X1) edge node[above] {} (X2)
(X2) edge node[above] {} (X3)
(X2) edge node[above] {} (Y2)
(XN) edge node[above] {} (XN+1)
(XN+1) edge node[above] {} (RX)
(LY) edge node[above] {} (Y1)
(Y1) edge node[above] {} (Y2)
(YN) edge node[above] {} (YN+1)
(X1) edge node[right] {} (Y1)
(XN+1) edge node[right] {} (YN+1)
(YN+1) edge node[above] {} (RY)
(X3) edge node[right] {} (P1)
(P1) edge node[above] {$\id$} (XN-1)
(XN-1)edge node[above] {$\id$}  (XN)
(Y2) edge node[below] {$\id$} (Y3)
(Y3)edge node[below] {$\id$} (P2)
(P2) edge node[below] {$\id$} (YN-1)
(X3)edge node[right] {$\id$} (Y3)
(XN-1)edge node[left] {$\id$}(YN-1)
(XN)edge node[right] {$\id$}(YN)
(YN-1)edge node[below] {$\id$}(YN)
;
}
\end{equation*}

\end{remark}

\begin{remark}
Since we will be working with lifting problems of morphisms with source objects of the form $\mu_j^i(R)$, we will take the shorthand notation $1^k$ for the element $1\in R=\mu_j^i(R)^k$ for $i-j+1 \le k \le i$.
\end{remark}

Now we can prove statements which are analogous to the case of the projective model structure for ordinary chain complexes. This will give us a classification of the class $J \boxslash$.

\begin{proposition}\label{7.4}
The class $J \boxslash$ is the class of levelwise epimorphisms.
\begin{proof}We will prove that both classes mutually include each other:
\begin{itemize}

\item[``$\subseteq$'':]
Let $p:X \to Y$ be in $J \boxslash$. So $p$ has the right lifting property with respect to all morphisms in $J$. Furthermore let $y \in Y^i$. To show that there exists a preimage of $y$ under $p$ consider the following lifting problem:
\begin{equation*}
\tikz[heighttwo,xscale=3,yscale=2,baseline]{
\node (Y) at (1,0) {$Y,$};
\node (X) at (1,1) {$X$};
\node (A) at (0,0) {$\mu_N^{i+N-1}(R)$};
\node (0) at (0,1) {$0$};
    
\draw[->]
(X) edge node[right] {$p$}(Y)
(0)edge(X)
(0)edge(A)
(A)edge node[above] {$\varphi_y$} (Y);
}
\end{equation*}
here $\varphi_y$ is the map uniquely defined by $1^i \mapsto y$, this is well defined by the adjunction between $\mu_N^{i+N-1}(\_)$ and $(\_)^i$ from Proposition~\ref{functoriality mu}. Now there exists, since $p \in J \boxslash$, a lift $h:\mu_N^{i+N-1}(R) \to X$, such that the following diagram commutes:
\begin{equation*}
\tikz[heighttwo,xscale=3,yscale=2,baseline]{
\node (Y) at (1,0) {$Y.$};
\node (X) at (1,1) {$X$};
\node (A) at (0,0) {$\mu_N^{i+N-1}(R)$};
\node (0) at (0,1) {$0$};
    
\draw[->]
(X) edge node[right] {$p$}(Y)
(A) edge node[above left] {$h$}(X)
(0)edge(X)
(0)edge(A)
(A)edge node[above] {$\varphi_y$} (Y);
}
\end{equation*}
Now $p \circ h = \varphi_y$. This means if we define $x:=h(1^i)$, that $y=\varphi_y(1^i)=p(h(1^i))=p(x)$ and hence we have a preimage of $y$, so $p$ is surjective.
\item[``$\supseteq$'':] Let $p:X\to Y$ be a levelwise epimorphism. Consider the following lifting problem: 
\begin{equation*}
\tikz[heighttwo,xscale=3,yscale=2,baseline]{
\node (Y) at (1,0) {$Y.$};
\node (X) at (1,1) {$X$};
\node (A) at (0,0) {$\mu_N^{i+N-1}(R)$};
\node (0) at (0,1) {$0$};
    
\draw[->]
(X) edge node[right] {$p$}(Y)
(0)edge(X)
(0)edge(A)
(A)edge node[above] {$\varphi$} (Y);
}
\end{equation*}
Now by the adjunction $  \mu_N^{i+N-1}(\_)\dashv (\_)^{i}$ we get that the morphism $\varphi$ comes from a map $\psi : R \to Y^i$, which is uniquely determined by the image of $1 \in R$, let $y:=\psi(1)$. Now, since $p$ is a levelwise epimorphism, there exists $x \in X^i$ such that $p(x)=y$. Consider $h':R \to X^i$, defined by $h'(1)=x$. Now we get an adjoint map $h$ to $h'$, and so the following diagram commutes:
\begin{equation*}
\tikz[heighttwo,xscale=3,yscale=2,baseline]{
\node (Y) at (1,0) {$Y.$};
\node (X) at (1,1) {$X$};
\node (R) at (0,0) {$\mu_N^{i+N-1}(R)$};
\node (0) at (0,1) {$0$};
    
\draw[->]
(X) edge node[right] {$p$}(Y)
(0)edge(X)
(0)edge(A)
(A)edge node[above] {$\varphi$} (Y)
(A) edge node[above left] {h}(X);
}
\end{equation*}
Since $i$ and $\varphi$ were arbitrary, we get $p \in J \boxslash$.
\end{itemize}
\end{proof}
\end{proposition}

\begin{proposition}\label{7.7}
The class $I\boxslash$ is precisely the class of levelwise surjective $N$-quasi-isomorphisms.
\begin{proof}
We will prove this by showing that each class is included in the other:
\begin{itemize}
\item[``$\subseteq$'':]Let $(p:X \to Y) \in I\boxslash$. First observe that $J \subset I$ directly gives $I \boxslash \subset J \boxslash $. So by Proposition~\ref{7.4} the morphism $p$ has to be a levelwise epimorphism.

Hence it remains to show that $p$ is an $N$-quasi-isomorphism. First we show that $p$ induces an epimorphism on cohomology. For this consider $y \in \ker(d_Y^r)^{i-N+1}$ and the lifting problem:
\begin{equation*}
\tikz[heighttwo,xscale=2,yscale=2,baseline]{
\node (Y) at (1,0) {$Y,$};
\node (X) at (1,1) {$X$};
\node (A2) at (0,0) {$\mu_N^{i}(R)$};
\node (A1) at (0,1) {$\mu_{N-r}^{i}(R)$};
    
\draw[->]
(X) edge node[right] {$p$}(Y)
(A1)edge node[above] {$0$}(X)
(A1)edge node[left]{$\iota_{i,N-r}$}(A2)
(A2)edge node[below] {$\varphi_y$} (Y);
}
\end{equation*}
where the bottom map is the map induced by $1^{i-N+1} \mapsto y$. This diagram commutes since $y \in \ker(d_Y^r)$. Since $p \in I \boxslash$ there is a lift $h: \mu_N^{i}(R) \to X$ such that the following square commutes:
\begin{equation*}
\tikz[heighttwo,xscale=2,yscale=2,baseline]{
\node (Y) at (1,0) {$Y.$};
\node (X) at (1,1) {$X$};
\node (A2) at (0,0) {$\mu_N^{i}(R)$};
\node (A1) at (0,1) {$\mu_{N-r}^{i}(R)$};
    
\draw[->]
(X) edge node[right] {$p$}(Y)
(A1)edge node[above] {$0$}(X)
(A1)edge node[left]{$\iota_{i,N-r}$}(A2)
(A2)edge node[below] {$\varphi_y$} (Y)
(A2) edge node[above left] {$h$} (X);
}
\end{equation*}
But the commutativity of the top triangle just means that $x:=h(1^{i-N+1}) \in \ker(d^r)$ and the bottom triangle gives $p(x)=y$. So $p$ induces a surjection on $\ker(d^r)$ for all $1 \le r \le N-1$ and hence $H_{(r)}(p)$ is an epimorphism as well.

In order to show that $H_{(r)}(p)$ is an monomorphism it suffices to show, that if for an $x\in \ker(d_X^r)^{i-r+1}$ there exists an $y'$ such that $y:=p(x)=d_Y^{N-r}(y')$, then there exists an $x'\in X$ such that $d_X^{N-r}(x')=x$.

To prove this, consider the lifting problem:
\begin{equation*}
\tikz[heighttwo,xscale=2,yscale=2,baseline]{
\node (Y) at (1,0) {$Y,$};
\node (X) at (1,1) {$X$};
\node (A2) at (0,0) {$\mu_N^{i}(R)$};
\node (A1) at (0,1) {$\mu_{r}^{i}(R)$};
    
\draw[->]
(X) edge node[right] {$p$}(Y)
(A1)edge node[above] {$\varphi_x$}(X)
(A1)edge node[left]{$\iota_{i,r}$}(A2)
(A2)edge node[below] {$\varphi_{y'}$} (Y);
}
\end{equation*}
where $\varphi_{y'}$ is again defined by $\varphi_{y'}(1^{i-N+1})=y'$ and $\varphi_x$ is defined by $\varphi_x(1^{i-r+1})=x$. Hence the diagram indeed commutes, since $p(x)=d_Y^{N-r}(y')$. Now we again get a morphism $h:\mu_N^{i}(R) \to X$ such that the square
\begin{equation*}
\tikz[heighttwo,xscale=2,yscale=2,baseline]{
\node (Y) at (1,0) {$Y.$};
\node (X) at (1,1) {$X$};
\node (A2) at (0,0) {$\mu_N^{i}(R)$};
\node (A1) at (0,1) {$\mu_{r}^{i}(R)$};
    
\draw[->]
(X) edge node[right] {$p$}(Y)
(A1)edge node[above] {$\varphi_x$}(X)
(A1)edge node[left]{$\iota_{i,r}$}(A2)
(A2)edge node[below] {$\varphi_{y'}$} (Y)
(A2) edge node[above left] {$h$} (X);
}
\end{equation*}
 commutes. The commutativity of the top triangle now yields that
\begin{align*}
d_X^{N-r}\circ h(1^{i-N+1})=h \circ d_{\mu_{N}^{i}(R)}^{N-r}(1^{i-N+1})=h(1^{i-r+1})= h \circ \iota_{i,r} (1^{i-r+1})=\varphi_x(1^{i-r+1})=x.
\end{align*}
Hence for $x':= h(1^{i-N+1})$ we get $d^{N-r}(x')=x$, and we have that $H_{(r)}(p)$ is injective. Which proves the claim.

\item[``$\supseteq$'':] Let $p:X \to Y$ be a surjective $N$-quasi-isomorphism and consider a lifting problem:
\begin{equation*}
\tikz[heighttwo,xscale=2,yscale=2,baseline]{
\node (Y) at (1,0) {$Y.$};
\node (X) at (1,1) {$X$};
\node (A2) at (0,0) {$\mu_N^{i}(R)$};
\node (A1) at (0,1) {$\mu_{N-r}^{i}(R)$};
    
\draw[->]
(X) edge node[right] {$p$}(Y)
(A1)edge node[above] {$\psi$}(X)
(A1)edge node[left]{$\iota_{i,N-r}$}(A2)
(A2)edge node[below] {$\varphi$} (Y);
}
\end{equation*}
Now since $p$ is surjective there is a short exact sequence in $\CC_N(\Mod_R)$ of the form
\begin{equation*}
\tikz[heighttwo,xscale=2,yscale=2,baseline]{
\node (K) at (-1,0.3) {$\ker(p)$};
\node (Y) at (1,0.3) {$Y$};
\node (X) at (0,0.3) {$X$};

\draw[->>]
(X) edge node[above] {$p$}(Y);
\draw[right hook ->]
(K) edge (X);
}.
\end{equation*}
Hence we can interpret $\ker(p)$ as a subcomplex of $X$, and since $p$ is an $N$-quasi-isomorphism we have that $H_{(r)}(\ker(p))=0$ by Lemma~\ref{5.7}. Let $x'\in X^{i-N+1}$ be a preimage of $y:=\varphi(1^{i-N+1})$ under $p$ and let $\widehat{h}:\mu_N^{i}(R) \to X$ be the map defined by $\widehat{h}(1^{i-N+1})=x'$. Now this morphism makes the bottom triangle in the following diagram commute,
\begin{equation*}
\tikz[heighttwo,xscale=1,yscale=1,baseline]{
\node (Y) at (2,0) {$Y,$};
\node (X) at (2,2) {$X$};
\node (A2) at (0,0) {$\mu_N^{i}(R)$};
\node (A1) at (0,2) {$\mu_{N-r}^{i}(R)$};
\node (L) at (0.75,1.25){\Huge{\Lightning}};

\draw[->]
(X) edge node[right] {$p$}(Y)
(A1)edge node[above] {$\psi$}(X)
(A1)edge node[left]{$\iota_{i,N-r}$}(A2)
(A2)edge node[below] {$\varphi$} (Y)
(A2) edge node[below right] {$\widehat{h}$} (X);
}
\end{equation*}
but in general not the top triangle. To solve this problem, let $x:=\psi(1^{i-N+r+1})$. now consider $d^r_X(x')-x$. This is an element of $\ker(p)$, since the commutativity of the outer square yields
\begin{align*}
p(d^r_X(x')-x)&=d^r_Y(p(x'))-p(x)\\
&=d^r_Y(y)-p(\psi(1^{i-N+r+1}))\\
&=d^r_Y(y)-\varphi(\iota_{i,r} (1^{i-N+r+1}))\\
&=d^r_Y(y)-\varphi(1^{i-N+r+1})\\ 
&=d^r_Y(y)-\varphi(d^r (1^{i-N+1}))\\
&=d^r_Y(y)-d^r(\varphi (1^{i-N+1}))\\
&=d^r_Y(y)-d_Y^r(y)=0.
\end{align*}
Furthermore we have $d^r(x'),x \in \ker(d^{N-r})$, since on one hand $d_{\ker(p)}^N=0$ and on the other hand $x= \psi(1^{i-N+r+1})$ and $\mu^i_{N-r}(R) \in \CC_{N-r}(\Mod_R)$.
 
So $d^r_X(x')-x=0$ in $H_{(N-r)}(\ker(p))$. Therefore there exists an element $\widehat{x}\in \ker(p)$ such that $d^r(\widehat{x})=d^r_X(x')-x$. Now we may replace $\widehat{h}$ by the unique map $h:\mu_N^{i}(R) \to X$ which satisfies $h(1^{i-N+1})=x'- \widehat{x}$. Then the diagram 
\begin{equation*}
\tikz[heighttwo,xscale=1,yscale=1,baseline]{
\node (Y) at (2,0) {$Y$};
\node (X) at (2,2) {$X$};
\node (A2) at (0,0) {$\mu_N^{i}(R)$};
\node (A1) at (0,2) {$\mu_{N-r}^{i}(R)$};

\draw[->]
(X) edge node[right] {$p$}(Y)
(A1)edge node[above] {$\psi$}(X)
(A1)edge node[left]{$\iota_{i,N-r}$}(A2)
(A2)edge node[below] {$\varphi$} (Y)
(A2) edge node[below right] {$h$} (X);
}
\end{equation*}
indeed commutes since 
\begin{align*}
h \circ \iota_{i,N-r}(1^{i-N+r+1})&=h(1^{i-N+r+1})=d^r_X \circ h (1^{i-N+1})\\
&= d^r_X(x'- \widehat{x})=d^r_X(x')-d^r_X(x')+x=x=\psi(1^{i-N+r+1})
\end{align*} 
and
\begin{align*}
p\circ h (1^{i-N+1})=p(x'- \widehat{x})=p(x')-0=y=\varphi(1^{i-N+1}).
\end{align*}
Hence the lifting problem has a solution. Now $p \in J \boxslash$ remains to be shown, but this is clear by Proposition~\ref{7.4}. So we altogether get $p \in I \boxslash$.
\end{itemize}
\end{proof}
\end{proposition}

\begin{proposition}\label{trivial cofibrations are weak equivalences}
The class $\boxslash (J \boxslash)$ is precisely the class of split monomorphisms with projective cokernel. In particular $\boxslash (J \boxslash)\subset W$.
\begin{proof}
We will prove this by showing that each class is included in the other:
\begin{itemize}

\item[``$\supseteq$'':] Let $i: X \to Y$ be a split monomorphism with projective cokernel. Then $i$ is isomorphic in the morphism category to $\iota_X:=\begin{bmatrix}
\id \\
0
\end{bmatrix}: X \to X \oplus C$ with $C$ projective. Now consider the lifting problem
\begin{equation*}
\tikz[heighttwo,xscale=1,yscale=1,baseline]{
\node (X) at (0,2) {$X$};
\node (Y) at (0,0) {$X\oplus C$};
\node (A) at (2,2) {$A$};
\node (B) at (2,0) {$B$};

\draw[->]
(X) edge node[left] {$\begin{bmatrix}
\id \\
0
\end{bmatrix}$}(Y)
(A) edge node[right] {$p$} (B)
(X) edge node[above]{$\varphi$}(A)
(Y) edge node[above] {$\psi$}(B);

}
\end{equation*}
with $p: A \to B \in J\boxslash$. By Proposition~\ref{7.4} $p$ is a levelwise epimorphism, and hence an epimorphism in $\CC_N(\Mod_R)$. This means we get a morphism $\widehat{\psi}: C \to A$ such that $p \circ \widehat{\psi}=\psi \circ  \iota_C$, where $\iota_C$ is the canonical inclusion $C \to X \oplus C$. Now consider the morphism $\begin{bmatrix}
\varphi & \widehat{\psi}
\end{bmatrix}: X \oplus C \to A$. This morphism makes the diagram 
\begin{equation*}
\tikz[heighttwo,xscale=2,yscale=1.5,baseline]{
\node (X) at (0,2) {$X$};
\node (Y) at (0,0) {$X\oplus C$};
\node (A) at (2,2) {$A$};
\node (B) at (2,0) {$B$};

\draw[->]
(X) edge node[left] {$\begin{bmatrix}
\id \\
0
\end{bmatrix}$}(Y)
(A) edge node[right] {$p$} (B)
(X) edge node[above]{$\varphi$}(A)
(Y) edge node[above] {$\psi$}(B)
(Y) edge node[above left] {$\begin{bmatrix}
\varphi & \widehat{\psi}\end{bmatrix}$}(A);
}
\end{equation*}
commute in $\CC_N(\Mod_R)$, since $$\begin{bmatrix}
\varphi & \widehat{\psi}\end{bmatrix} \circ \begin{bmatrix}
\id \\
0
\end{bmatrix} = \varphi$$
and 
$$p \circ \begin{bmatrix}
\varphi & \widehat{\psi}\end{bmatrix} =\begin{bmatrix}
p \circ \varphi & p \circ \widehat{\psi}\end{bmatrix}=\begin{bmatrix}\psi\circ \iota_X & \psi \circ \iota_C \end{bmatrix}= \psi  .$$
\item[``$\subseteq$'':] Let $i:X \to Y \in \boxslash(J \boxslash)$. First consider for the splitting part the diagram
\begin{equation*}
\tikz[heighttwo,xscale=1,yscale=1,baseline]{
\node (X) at (0,2) {$X$};
\node (Y) at (0,0) {$Y$};
\node (X2) at (2,2) {$X$};
\node (0) at (2,0) {$0$};

\draw[->]
(X) edge node[left] {$i$}(Y)
(X2) edge (0)
(Y) edge (0)
(X) edge node[above] {$\id_X$}(X2);

}
\end{equation*}
in $\CC_N(\Mod_R)$. Since $i \in  \boxslash (J \boxslash)$ and $0: X \to 0 \in J \boxslash$ by Proposition~\ref{7.4} we get $h: Y \to X$ such that $h \circ i = \id_X$ and hence $i$ is a split monomorphism.
Now to prove that $C:=\coker(i)$ is projective in $\CC_N(\Mod_R)$ we need to show that for every levelwise epimorphism $g:A \to B$ and $f: C \to B$ there exists an $\widehat{f}:C \to A$ such that the following diagram commutes in $\CC_N(\Mod_R)$
\begin{equation*}
\tikz[heighttwo,xscale=1,yscale=1,baseline]{
\node (M) at (2,2) {$A$};
\node (N) at (2,0) {$B.$};
\node (X) at (0,0) {$C$};

\draw[->>]
(M) edge node[right] {$g$}(N);
\draw[->]
(X) edge node[above] {$f$}(N);
\draw[dashed, ->]
(X) edge node[above left]{$\widehat{f}$}(M);
}
\end{equation*}
To prove this consider the commutative diagram
\begin{equation*}
\tikz[heighttwo,xscale=1,yscale=1,baseline]{
\node (X) at (0,4) {$X$};
\node (Y) at (0,2) {$Y$};
\node (A) at (2,4) {$A$};
\node (B1) at (2,2) {$B$};
\node (C) at (0,0) {$C$};
\node (B2) at (2,0) {$B$};

\draw[->]
(X) edge node[left] {$i$}(Y)
(A) edge node[right]{$g$}(B1)
(Y) edge node[above] {$f \circ \pi$}(B1)
(X) edge node[above] {$0$}(A)
(Y) edge node[left] {$\pi$}(C)
(B1) edge node[right] {$\id_B$}(B2)
(C) edge node[above] {$f$}(B2);

}
\end{equation*}
in $\CC_N(\Mod_R)$. Here the left column is the short exact sequence exhibiting $C$ as cokernel of $i$. Now since $i$ is a split monomorphism we get that $\pi$ has a splitting $j:C \to Y$ as well. Furthermore we get by Proposition~\ref{7.4} that $g \in J\boxslash$ . So using $i \in \boxslash(J \boxslash)$ we get a morphism $h: Y \to A$ such that $g \circ h = f \circ \pi$. So altogether we get the following diagram
\begin{equation*}
\tikz[heighttwo,xscale=1,yscale=1,baseline]{
\node (X) at (0,4) {$X$};
\node (Y) at (0,2) {$Y$};
\node (A) at (2,4) {$A$};
\node (B1) at (2,2) {$B$};
\node (C) at (0,0) {$C$};
\node (B2) at (2,0) {$B.$};

\draw[->]
(X) edge node[left] {$i$}(Y)
(A) edge node[right]{$g$}(B1)
(Y) edge node[above] {$f \circ \pi$}(B1)
(X) edge node[above] {$0$}(A)
(Y) edge node[left] {$\pi$}(C)
(B1) edge node[right] {$\id_B$}(B2)
(C) edge node[above] {$f$}(B2)
(Y) edge node[above left] {$h$}(A);
\draw[dashed,->]
(C) edge [bend right=30] node[right] {$j$}(Y)
;}
\end{equation*}
in $\CC_N(\Mod_R)$, which commutes except for the dashed arrow.
In particular we get
$$g \circ h  \circ j= f \circ \pi \circ j = f$$
which gives with $\widehat{f}:=h  \circ j$ the desired lift.

To conclude now the ``in particular'' part observe that since $C$ is projective we have that $\lambda^{\PP}_C : \PP_N(C) \to C$ splits. So, as a summand of an acyclic $N$-complex, $C$ is  nullhomotopic and hence acyclic itself. Corollary~\ref{5.7} now gives that $i $ has to be an $N$-quasi-isomorphism, and so $i \in W$.
\end{itemize}
\end{proof}
\end{proposition}

\begin{remark}
We remind the reader that the notion of a relative $\Psi$-cell complex for $\Psi$ a set of morphisms is defined as a transfinite composition of pushouts of morphisms in the set $\Psi$. We will use the shorthand notation of $\Psi$-cell for the class of morphisms of this form for a given $\Psi$. 
\end{remark}

Now in order to apply the small object argument in the form presented in \cite{hovey} we need that the source objects of the morphisms in $I$ and $J$ are all compact, since in this case they are in particular small with respect to $\Psi$-cell, for any class of morphisms $\Psi$.

\begin{remark}
We remind here that for a category $\C$ which admits filtered colimits an object $X \in \C$ is  compact if and only if the functor $\C(X,\_)$ commutes with filtered colimits.
\end{remark}


\begin{proposition}\label{7.9}
Let $R$ be a ring. Then all $\mu_k^i(R)$ are compact in $\CC_N(\Mod_R)$ for $1 \le k\le N$. 
\begin{proof}
Consider a filtered colimit $\colim_{m \in M}X_m$ and a morphism $f:\mu_k^i(R) \to \colim_{m \in M}X_m$ in $\CC_N(\Mod_R)$. This morphism is now uniquely defined by $$f^{i-k+1}(1)\in  \left(\colim_{m \in M}X_m\right)^i=\colim_{m \in M}X_m^i.$$ Furthermore by pointwise calculation we have $$f^{i-k+1}(1) \in \ker(d_{\colim_{m \in M}X_m}^k)=\ker(\colim_{m \in M}d_{X_m}^k)=\colim_{m \in M}(\ker(d_m^k)).$$
Here the last equality comes from the fact that $\Mod_R$ is a Grothendieck category, and hence satisfies property (AB4).
Using these observations we get a bijection $$\CC_N(\Mod_R)(\mu_k^i(R), \colim_{m \in M}X_m)\cong \Mod_R(R,\colim_{m \in M}(\ker(d_m^k))).$$ Now $R$ is compact in $\Mod_R$ and using the above observation we get
\begin{align*}
\CC_N(\Mod_R)(\mu_k^i(R), \colim_{m \in M}X_m)&\cong \Mod_R(R,\colim_{m \in M}(\ker(d_m^k)))\\
&\cong \colim_{m \in M} \Mod_R(R,(\ker(d_m^k)))\\
&\cong  \colim_{m \in M} \CC_N(\Mod_R)(\mu_k^i(R),X_m),
\end{align*}
which finishes the proof.
\end{proof}
\end{proposition}

For the case of $J$, observe that the zero complex is the zero object in $\CC_N({\Mod_R})$ and hence compact.

After these statements we can finally prove that this indeed gives a model structure, and that its model theoretic homotopy category is the Gabriel-Zismann localization $\CC_N(\Mod_R)[W^{-1}]$.

\begin{theorem}
The classes $I$, $J$ and $W$ cofibrantly generate a model structure on $\CC_N(\Mod_R)$ such that:
\begin{enumerate}
\item The weak equivalences  are precisely the $N$-quasi-isomorphisms.
\item The fibrations are precisely the levelwise epimorphisms.
\end{enumerate}
\begin{proof}
We will use Theorem~2.1.19 in \cite{hovey}. To apply this we need to check the following:
\begin{itemize}
\item The class $W$ has the ``two out of three'' property and is closed under retracts,
\item $\CC_N(\Mod_R)$ is complete and cocomplete,
\item the domains of $I$ are small relative to $I$-cell,
\item the domains of $J$ are small relative to $J$-cell,
\item $J\text{-cell} \subset  W \cap (\boxslash (I \boxslash))$,
\item and $I \boxslash = W \cap (J \boxslash)$.
\end{itemize}
First of all the set $W$ satisfies ``two out of three'' and is closed under retracts since it is obtained by pulling back isomorphisms along a functor and isomorphisms have these two properties. Furthermore $\CC_N(\Mod_R)$ is complete and cocomplete by Proposition~\ref{complete/cocomplete} and we have by Proposition~\ref{7.9} that the source objects of $I$ and $J$ are compact, so they are in particular small with respect to $\Psi$-cell, for $\Psi$ any class of morphisms. Proposition~\ref{7.7} and Proposition~\ref{7.4} give $(I\boxslash) = W \cap(J \boxslash)$. Hence we just have to show 
\begin{equation*}
J \text{-cell} \subset W \cap (\boxslash (I \boxslash)).
\end{equation*} 
Lemma~2.1.10 in  \cite{hovey} in general shows $J\text{-cell} \subset \boxslash(J \boxslash)$. Proposition~\ref{trivial cofibrations are weak equivalences} gives $\boxslash(J \boxslash) \subset W$ and $J \subset I$ directly implies $ \boxslash(J \boxslash) \subset \boxslash (I \boxslash)$. This shows the above inclusion.

Now using Theorem~2.1.19 in \cite{hovey} we get a model structure generated by $I$,$J$ and $W$, so $W$ is precisely the class of weak equivalences, which shows property (1), and the fibrations are $J\boxslash$, which by Proposition~\ref{7.4} is precisely the class of levelwise epimorphisms and yields property (2).
\end{proof}
\end{theorem}

\begin{remark}
We remind here that by Definition~1.2.4 and Lemma~1.2.5 in \cite{hovey} two morphisms $f,g:X \to Y$ between fibrant-cofibrant objects in a category with a model structure are homotopy equivalent (written $f\sim g$) if one of the following equivalent statements hold:
\begin{enumerate}
\item There exists a cylinder object $C_X$ for $X$, i.\,e.\ a factorisation of the fold map $X \amalg X \to X$ into a cofibration $X \amalg X  \xrightarrow{\iota_0 + \iota_1} C_X$ followed by a weak equivalence $C_X \to X$, and a left homotopy $\eta:C_X \to Y$ such that $f= \eta \circ \iota_0$ and $g= \eta \circ \iota_1$.
\item There exists a path object $P_Y$ for $Y$, i.\,e.\ a factorisation of the diagonal map $Y  \to Y \times Y$ into a weak equivalence $Y \to P_Y$ followed by a fibration $P_Y \xrightarrow{(\pi_0, \pi_1)} Y \times Y$, and a right homotopy $\nu:X \to P_Y$ such that $f=  \pi_0 \circ \nu$ and $g= \pi_1 \circ \nu$.
\end{enumerate} 
Furthermore we will denote for a category $\C$ with a model structure the full subcategory of fibrant-cofibrant objects of $\C$ by $\C^{cf}$.
\end{remark}

\begin{corollary}\label{equiv homotopycat gabriel zisman}
The obvious functor $\CC_N(\Mod_R)^{cf}_{/\sim} \to \CC_N(\Mod_R)[W^{-1}]$ is an equivalence.
\begin{proof}
This is a direct application of Theorem~1.2.10 in \cite{hovey}.
\end{proof}
\end{corollary}

Now to conclude this section we will show that just as in the case of chain complexes we get a sequence of equivalences: $$\DD_N(\Mod_R)\xleftarrow{\sim} \hproj_{N}(\Mod_R) \xleftarrow{\sim}  \CC_N(\Mod_R)^{cf}_{/\sim} \xrightarrow{\sim} \CC_N(\Mod_R)[W^{-1}].$$
In order to establish these equivalences we need to inspect what it means to be a cofibrant object, which will give us essential surjectivity of the functor in the middle.

\begin{proposition}
Let $X$ be a cofibrant object in $\CC_N(\Mod_R)$ with respect to the above defined model structure. Then $X^n$ is projective for all $n \in \Z$.
\begin{proof}
Let $n \in \N$, $p:M \to Q$ an epimorphism in $\Mod_R$ and $f:X^n\to Q$ a morphism in $\Mod_R$. Consider the following diagram in $\Mod_R$, where we need to construct the lift $\widehat{f}$
\begin{equation*}
\tikz[heighttwo,xscale=1,yscale=1,baseline]{
\node (M) at (2,2) {$M$};
\node (N) at (2,0) {$Q.$};
\node (X) at (0,0) {$X^n$};

\draw[->>]
(M) edge node[right] {$p$}(N);
\draw[->]
(X) edge node[above] {$f$}(N);
\draw[dashed, ->]
(X) edge node[above left]{$\widehat{f}$}(M);
}
\end{equation*}
To solve this, consider the lifting problem
\begin{equation*}
\tikz[heighttwo,xscale=1,yscale=1,baseline]{
\node (M) at (2,2) {$\mu_N^{n}(M)$};
\node (N) at (2,0) {$\mu_N^{n}(Q),$};
\node (X) at (0,0) {$X$};
\node (0) at (0,2) {$0$};

\draw[->]
(M) edge node[right] {$\mu_N^n(p)$}(N)
(0) edge (M)
(0) edge (X)
(X) edge node[above] {$\widetilde{f}$}(N);

}
\end{equation*}
where $\widetilde{f}$ is the adjoint map to $f$ under the adjunction $(\_)^n \dashv \mu_{N}^{n}$. Now $\mu_{N}^n(p)$ is a levelwise epimorphism and hence a fibration. Since we have $\mu^n_{N}(Y)\cong \PP_N(\mu^n_1(Y))$ for all $Y \in \Mod_R$, it is even a trivial fibration. So, since $X$ is cofibrant, we get a lift $h: X \to \mu_N^{n}(M)$ such that the following diagram commutes:
\begin{equation*}
\tikz[heighttwo,xscale=1,yscale=1,baseline]{
\node (M) at (2,2) {$\mu_N^{n}(M)$};
\node (N) at (2,0) {$\mu_N^{n}(Q).$};
\node (X) at (0,0) {$X$};
\node (0) at (0,2) {$0$};

\draw[->]
(M) edge node[right] {$\mu_N^n(p)$}(N)
(0) edge (M)
(0) edge (X)
(X) edge node[above] {$\widetilde{f}$}(N)
(X) edge node[above left] {$h$}(M);
}
\end{equation*}
Now passing to adjoints yields 
\begin{equation*}
\tikz[heighttwo,xscale=1,yscale=1,baseline]{
\node (M) at (2,2) {$M$};
\node (N) at (2,0) {$Q.$};
\node (X) at (0,0) {$X^n$};
\node (0) at (0,2) {$0$};

\draw[->]
(M) edge node[right] {$p$}(N)
(0) edge (M)
(0) edge (X)
(X) edge node[above] {$f$}(N)
(X) edge node[above left] {$\widehat{h}$}(M);
}
\end{equation*}
So we get with $\widehat{f}:=\widehat{h}$ the desired lift, which shows that $X^n$ is projective for all $n \in \Z$.
\end{proof}
\end{proposition}

\begin{proposition}\label{7.12}
Let $X$ be a cofibrant object. Then $X$ is h-projective.
\begin{proof}
Let $f:X \to T$ be a morphism with $T$ an acyclic complex. Consider the following lifting problem:
\begin{equation*}
\tikz[heighttwo,xscale=1,yscale=1,baseline]{
\node (PA) at (2,2) {$\PP_N(T)$};
\node (A) at (2,0) {$T.$};
\node (X) at (0,0) {$X$};
\node (0) at (0,2) {$0$};

\draw[->]
(PA) edge node[right] {$\lambda^{\PP}_T$}(A)
(0) edge (PA)
(0) edge (X)
(X) edge node[above] {$f$}(A);
}
\end{equation*}
Now since $X$ is cofibrant and $\lambda^{\PP}_T$ is a levelwise epimorphism between acyclic objects we get a lift $h:X \to \PP_N(T)$ such that $f= \lambda^{\PP}_T \circ h$. So $f$ factors over an $\Se_\oplus$-injective-projective object which shows that $f=0$ in $\KK_N(\Mod_R)$. Hence for all acyclic objects $T$ we have $\KK_N(\Mod_R)(X,T)=0$ and so $X$ is h-projective.
\end{proof}
\end{proposition}

The only thing missing in the proof of the statement announced below Corollary~\ref{equiv homotopycat gabriel zisman} is that the two equivalence relations in the homotopy categories coincide.

\begin{proposition}\label{7.13}
Let $f,g: X \to Y$ be two morphisms in $\CC_N(\Mod_R)$ between two cofibrant objects. Then they are homotopic with respect to the above model structure if and only if $f-g:X \to Y$ is  zero in $\KK_N(\Mod_R)$.
\begin{proof}
First of all observe that by Proposition~\ref{7.4} $X$ and $Y$ are fibrant. We will show both implications separately, since they use quite different techniques.
\begin{description}
\item[\normalfont $\Rightarrow$] Let $C_X$ be a cylinder object for $X$, hence we have morphisms $\pi_X: X\oplus X \to C_X$ and $\varphi_X: C_X \to X$ such that $\pi_X$ is a cofibration, $\varphi_X$ is a weak equivalence and we have $\varphi_X \circ \pi_X = \begin{bmatrix}
\id_X & \id_X
\end{bmatrix}$. Furthermore by the assumptions $X$ and $X\oplus X$ are cofibrant and hence h-projective. So $C_X$ is cofibrant as well, since $\pi_X$ is a cofibration. 
Now consider the following commutative diagram in $\KK_N(\Mod_R)$:

\begin{equation*}
\tikz[heighttwo,xscale=1.5,yscale=1.5,baseline]{

\node (X2) at (0,0) {$X$};
\node (XX2) at (2,0) {$X \oplus X$};
\node (CX2) at (4,0) {$C_X$};
\node (SX) at (6,0) {$\Sigma(X).$};
\node (X3) at (0,2) {$X$};
\node (XX3) at (2,2) {$X \oplus X$};
\node (X4) at (4,2) {$X$};
\node (SX3) at (6,2) {$\Sigma(X)$};
    
\draw[->]
(X2) edge node[above]{$\begin{bmatrix} \id \\ -\id
\end{bmatrix}$} (XX2)
(XX2) edge node[above]{$\pi_X$} (CX2)
(X3) edge node[above]{$\begin{bmatrix} \id \\0
\end{bmatrix}$} (XX3)
(XX3) edge node[above] {$\begin{bmatrix}
0 & \id_X
\end{bmatrix}$} (X4)
(CX2) edge node[above]{$0$} (SX)
(X4) edge node[above]{$0$}(SX3)
(CX2) edge node[right] {$\varphi_X$} (X4)
(X2) edge node[left] {$\id_X$}(X3)
(XX2) edge node[right] {$\begin{bmatrix}
 \id_X & 0 \\
 \id_X & \id_X
\end{bmatrix}$}(XX3)
(SX) edge node[right] {$\id_{\Sigma X}$}(SX3)

}
\end{equation*}
Using this diagram one gets that bottom row is a triangle, since the top one is clearly one.
Now the fact that $f$ and $g$ are homotopic with respect to our model structure, yields a commuting diagram in $\CC_N(\Mod_R)$ and hence also in $\KK_N(\Mod_R)$ of the form:
\begin{equation*}
\tikz[heighttwo,xscale=1,yscale=1,baseline]{
\node (X) at (0,2) {$X$};
\node (XX) at (2,2) {$X \oplus X$};
\node (CX) at (4,2) {$C_X$};
\node (SX) at (6,2) {$\Sigma(X)$};
\node (01) at (0,0) {$0$};
\node (Y1) at (2,0) {$Y$};
\node (Y2) at (4,0) {$Y$};
\node (02) at (6,0) {$\Sigma(0),$};

\draw[->]
(01) edge (Y1)
(Y2) edge (02)
(X) edge node[above] {$ \begin{bmatrix} \id \\ -\id
\end{bmatrix}$} (XX)
(XX) edge node[above] {$\pi_X$}(CX)
(CX) edge (SX)
(XX) edge node[left] {$[f,g]$}(Y1)
(CX) edge node[right] {$\eta$}(Y2)
(Y1) edge node[above] {$\id_Y$}(Y2);
}
\end{equation*}
where $\eta$ is a homotopy between $f$ and $g$. By the filling axiom this yields a commutative diagram in $\KK_N(\Mod_R)$:
\begin{equation*}
\tikz[heighttwo,xscale=1,yscale=1,baseline]{
\node (X) at (0,2) {$X$};
\node (XX) at (2,2) {$X \oplus X$};
\node (CX) at (4,2) {$C_X$};
\node (SX) at (6,2) {$\Sigma(X)$};
\node (01) at (0,0) {$0$};
\node (Y1) at (2,0) {$Y$};
\node (Y2) at (4,0) {$Y$};
\node (02) at (6,0) {$\Sigma(0).$};

\draw[->]
(X) edge (01)
(SX) edge (02)
(01) edge (Y1)
(Y2) edge (02)
(X) edge node[above] {$ \begin{bmatrix} \id \\ -\id
\end{bmatrix}$} (XX)
(XX) edge (CX)
(CX) edge (SX)
(XX) edge node[left] {$[f,g]$}(Y1)
(CX) edge node[right] {$\eta$}(Y2)
(Y1) edge node[above] {$\id_Y$}(Y2)
;
}
\end{equation*}
Now the commutativity of the leftmost square just says that in $\KK_N(\Mod_R)$ we have 
\begin{align*}
f-g= \begin{bmatrix} f & g
\end{bmatrix} \circ \begin{bmatrix} 1 \\ -1
\end{bmatrix}=0
\end{align*}
and thus $f-g$ is nullhomotopic.
\item[\normalfont $\Leftarrow$] We know that $f-g$ is nullhomotopic, hence $f-g$ factors as $\lambda_Y^\PP \circ \varphi$ for some $\varphi: X \to \PP_N(Y)$. Now consider the commutative diagram
\begin{equation*}
\tikz[heighttwo,xscale=1,yscale=1,baseline]{
\node (XX) at (4,0) {$Y\oplus Y.$};
\node (PX) at (2,2) {$Y \oplus \PP_N(Y)$};
\node (X) at (0,0) {$Y$};

\draw[->]
(X) edge node[below]{$ \begin{bmatrix}
 \id_Y \\
  \id_Y
\end{bmatrix}$}(XX)
(X) edge node[above left]{$ \begin{bmatrix} 
\id_Y \\ 
0 
\end{bmatrix}$}(PX)
(PX) edge node[above right] {$ \begin{bmatrix} 
\id_Y & 0 \\
\id_Y & \lambda^\PP_X
\end{bmatrix}$} (XX)
}
\end{equation*}
We have that $\lambda^\PP_Y$ is levelwise surjective. So we get 
$$ \im \left( \begin{bmatrix} 
\id &0\\
 \id & \lambda_N^\PP
\end{bmatrix}\right)=\{(y,y)+(0,y')|y,y' \in Y\}=Y\oplus Y, $$ which means that $Y \oplus \PP_N(Y) \to Y \oplus Y$ is a fibration. Furthermore, since $\PP_N(Y)\cong 0$ in $\KK_N(\Mod_R)$ we have that $\begin{bmatrix} 
\id_Y \\ 
0 
\end{bmatrix}$ is a weak equivalence. So the above diagram exhibits $Y \oplus \PP_N(Y)$ as a path object.

 Now we may construct a left homotopy by defining $\eta: X \to Y \oplus \PP_N(Y)$ as $\begin{bmatrix} 
f \\
 -\varphi
\end{bmatrix}$, since we get 
\begin{align*}
 \begin{bmatrix} 
\id & 0 \\ 
\id & \lambda_Y^\PP
\end{bmatrix}
\circ \eta
&=
 \begin{bmatrix} 
\id & 0 \\ 
\id & \lambda_Y^\PP
\end{bmatrix}
\circ 
\begin{bmatrix} 
f \\
 -\varphi
\end{bmatrix} 
 = 
\begin{bmatrix} 
f \\ f -\lambda_Y^\PP \circ  \varphi
\end{bmatrix}\\
&=
\begin{bmatrix} 
f \\
 f-(f-g)
\end{bmatrix}=\begin{bmatrix} 
f \\ g
\end{bmatrix}
\end{align*}
So the following diagram commutes:
\begin{equation*}
\tikz[heighttwo,xscale=1,yscale=1,baseline]{
\node (X) at (0,2) {$X$};
\node (PY) at (-4,0) {$Y \oplus \PP_N(Y)$};
\node (YY) at (0,0) {$Y\oplus Y,$};

\draw[->]
(X) edge node[right]{$ \begin{bmatrix}
 f\\  g
\end{bmatrix}$}(YY)
(PY) edge node[below]{$ \begin{bmatrix} 
\id & 0 \\ 
\id & \lambda_X^\PP
\end{bmatrix}$}(YY)
(X) edge node[above left] {$ \eta$} (PY)
}
\end{equation*}
which proves the claim.
\end{description}
\end{proof}
\end{proposition}

Now we just need to put the above arguments together to get:

\begin{theorem}\label{equiv hproj homotopycat}
There is an equivalence of categories $\CC_N(\Mod_R)^{cf}_{/\sim} \xrightarrow{\sim}  \hproj_{N}(\Mod_R)$ induced by the identity on $\CC_N(\Mod_R)$.
\begin{proof}
First observe that by definition every $N$-complex is fibrant.
Consider the identity functor 
$$\CC_N(\Mod_R) \to \CC_N(\Mod_R).$$
This induces the inclusion functor 
$$\CC_N(\Mod_R)^{cf}\to \CC_N(\Mod_R).$$
Composing this with the projection $\CC_N(\Mod_R)\to \KK_N(\Mod_R)$ now yields 
$$\CC_N(\Mod_R)^{cf}\to \KK_N(\Mod_R).$$
 By Proposition~\ref{7.12} the image of this functor is contained in $\hproj_N(\Mod_R)$, so we get a functor 
 $$\CC_N(\Mod_R)^{cf}\to \hproj_N(\Mod_R).$$ In particular we get that this functor induces the projection $$\CC_N(\Mod_R)(X,Y) \to \KK_N(\Mod_R)(X,Y):=\CC_N(X,Y)/I$$ on morphisms, but by Proposition~\ref{7.13} this induces a bijection $$\CC_N(\Mod_R)(X,Y)^{cf}_{/\sim} \to \KK_N(\Mod_R)(X,Y).$$

To show that this functor is also essentially surjective consider $X \in  \hproj_{N}(\Mod_R)$ and its cofibrant replacement $Q(X)$. Now both $X$ and $Q(X)$ are h-projective by Proposition~\ref{7.12}, so we get by Proposition~\ref{quasi iso between hproj} that every weak equivalence $f:Q(X) \to X$ is already a homotopy equivalence. In particular this yields essential surjectivity, which finishes the claim.
\end{proof}
\end{theorem}
This gives us our desired sequence of equivalences.
\begin{corollary}
There is a sequence of canonical equivalences
$$\DD_N(\Mod_R)\xleftarrow{\sim} \hproj_{N}(\Mod_R)\xleftarrow{\sim} \CC_N(\Mod_R)^{cf}_{/\sim} \xrightarrow{\sim} \CC_N(\Mod_R)[W^{-1}].$$
\begin{proof}
The first equivalence is Corollary~\ref{equiv derived hproj}, the second is Theorem~\ref{equiv hproj homotopycat} and the third is Corollary~\ref{equiv homotopycat gabriel zisman}.
\end{proof}
\end{corollary}